 \def\ekv#1#2{\begeq\label{#1}#2\endeq}
\def\eekv#1#2#3{\begin{eqnarray}\label{#1}#2 \\ #3
\nonumber\end{eqnarray}}
\def\eeekv#1#2#3#4{\begin{eqnarray}\label{#1}#2 \\ #3
\nonumber\\#4\nonumber\end{eqnarray}}
\def\eeeekv#1#2#3#4#5{\begin{eqnarray}\label{#1}#2 \\ #3
\nonumber\\#4\nonumber\\#5\nonumber\end{eqnarray}}

\def\iint{\int\hskip -2mm\int} \def\iiint{\int\hskip -2mm\int\hskip
-2mm\int}
\def\iiiint{\int\hskip -2mm\int\hskip -2mm\int\hskip -2mm\int}
  
scaled \magstep1  

\font\lagom=cmr10 at 10pt
  \def\3{\vert \hskip -1pt\vert\hskip -1pt\vert }

\def\ably{arbitrarily}

\def\asy{asymptotic}

\def\bdd{bounded}
\def\bdy{boundary}

\def\coef{coefficient}

\def\ctf{canonical transformation}

\def\dop{differential operator}
\def\ef{eigenfunction}
\def\ev{eigenvalue} \def\e{equation}
\def\fu{function}
\def\fy{family}

\def\fop{Fourier integral operator}
\def\hol{holomorphic}

\def\indep{independent}
\def\inv{^{-1}}
\def\lhs{left hand side}
\def\KFP{Kramers-Fokker-Planck}
\def\mfld{manifold}

\def\neigh{neighborhood}
\def\no#1{(\ref{#1})}
\def\nondeg{non-degenerate}
\def\3{\vert\hskip -1pt\vert\hskip -1pt\vert }
\def\on{orthonormal}

\def\op{operator} 

\def\pb{problem}    \def\pert{perturbation}
 \def\pol{polynomial} 
 \def\pop{pseudodifferential operator}
  \def\res{resonance} \def\rhs{right hand side}

\def\sa{selfadjoint} 
\def\scl{semiclassical}
\def\sm{\setminus }
 \def\st{strictly}   \def\sufly{sufficiently}
  
  \def\trans{^t\hskip -2pt}
\def\tr{{\rm tr\,}}
 \def\uf{uniform} \def\ufly{uniformly}
\def\vf{vector field} \def\wrt{with respect to} 
\def\Re{{\mathrm Re\,}}
\def\Im{{\mathrm Im\,}}
\def\R{{\bf R}}

\def\OO{{\cal O}}
\def\dd{\partial }
\def\cS{{\cal S}}

\def\D{\Delta }
\documentclass[12pt]{article}
\usepackage{amssymb,latexsym}
\title{Tunnel effect for Kramers-Fokker-Planck
type \op{}s} \author{Fr\'ed\'eric H\'erau\\ \lagom
Laboratoire de Math{\'e}matiques \\ \lagom Universit{\'e} de Reims \\ \lagom
Moulin de
la Housse B.P. 1039 \\ \lagom 51687 Reims cedex 2, France \\ \lagom
herau@univ-reims.fr \\ \lagom UMR 6056--CNRS \and Michael Hitrik\\ \lagom Department of
Mathematics \\ \lagom University of California \\ \lagom Los Angeles
\\ \lagom CA 90095-1555, USA\\ \lagom hitrik@math.ucla.edu
\and  Johannes
Sj{\"o}strand\\ \lagom CMLS\\ \lagom Ecole Polytechnique\\ \lagom FR  91120
Palaiseau c\'edex, France\\ \lagom johannes@{}math.polytechnique.fr \\ \lagom
UMR7640--CNRS}  \date{}
\def\wrtext#1{\relax\ifmmode{\leavevmode\hbox{#1}}\else{#1}\fi}
\def\abs#1{\left|#1\right|} \def\begeq{\begin{equation}}
\def\endeq{\end{equation}}

\def\part#1{\frac{\partial}{\partial #1}}

\renewcommand{\exp}{\mbox{\rm exp\,}}
\newcommand{\supp}{\mbox{\rm supp\,}}

\newtheorem{dref}{Definition}[section]
\newtheorem{lemma}[dref]{Lemma}
\newtheorem{theo}[dref]{Theorem}
\newtheorem{prop}[dref]{Proposition}
\newtheorem{hyp}[dref]{Hypothesis}
\newtheorem{remark}[dref]{Remark}

\newtheorem{cor}[dref]{Corollary}

\newenvironment{proof}{\par\noindent{{\bf Proof}}}{\hfill$\Box$
\medskip}

\begin{document}
\maketitle

\begin{abstract}
\par We consider operators of Kramers-Fokker-Planck type in the
semi-classical limit such that the exponent of the associated Maxwellian
is a Morse function with two local minima and a saddle
point. Under suitable additional assumptions we establish the complete
asymptotics of the exponentially small splitting between the first two
eigenvalues.
\smallskip\par
\centerline{\bf R\'esum\'e}
\smallskip\par
On consid\`ere des op\'erateurs du type de Kramers-Fokker-Planck dans la
limite semi-classique tels que l'exposant du maxwellien associ\'e soit
une fonction de Morse avec deux minima et un point selle. Sous des
hypoth\`eses suppl\'ementaires convenables on \'etablit un
d\'eveloppement asymptotique complet de l'\'ecart exponentiellement petit
entre les deux premi\`eres valeurs propres.
\end{abstract}

\vskip 2mm
\noindent
{\bf Mathematics Subject Clas\-si\-fi\-ca\-tion 2000}: 35P15, 35P20, 47A10, 47B44, \\
81Q20, 81Q60, 82C05, 82C31
\tableofcontents
\section{Introduction}\label{int}
\setcounter{equation}{0}

This paper is a natural continuation of the work \cite{HeSjSt},
investigating the low lying eigenvalues of the Kramers-Fokker-Planck
operator
\ekv{int.1}
{P=y\cdot h\partial _x-V'(x)\cdot h\partial _y+{\gamma
\over 2}(-h\partial _y+y)\cdot (h\partial _y+y),\ x,y\in {\bf R}^n,
}
where $\gamma >0$. Physically the \scl{} limit $h\to 0$ corresponds
to the low temperature limit. As explained in \cite{HeSjSt}, the
original motivation for that work was to give more explicit versions
of some results in \cite{HeNi} and later in \cite{HelNi}, giving
estimates on the time of return to equilibrium, or more or less
equivalently, on the gap between the first eigenvalue 0 (when the
potential $V$ tends to plus infinity sufficiently fast at
infinity) and the second eigenvalue. See also \cite{Vi} for further
developments in that direction. The methods of those works as well as
the one of Eckmann and Hairer \cite{EcHa} are inspired by those of hypoellipicity for
H\"ormander type operators. We will not repeat here all the
motivations of \cite{HeSjSt}, coming also from the works \cite{DeVi, Ko}
and others, including recent developments in pseudospectral theory.

\par In \cite{HeSjSt} it was assumed that $V$ in (\ref{int.1}) is a smooth real Morse
function with finitely many critical points, $U_j$, $j=1,..,N$,
that $\partial ^\alpha
V(x)$ is bounded for all multiindices $\alpha $ of length $|\alpha
|\ge 2$ and that $|\nabla V(x)|\ge 1/C$ for
$|x|\ge C$ for some $C>0$. Under these assumptions it was shown that the
eigenvalues in any disc $D(0,Ch)$ are of the form
\ekv{int.2}
{
\lambda _{j,k}(h)\sim h(\mu _{j,k}+h^{1/N_{j,k}}\mu _{j,k,1}+h^{2/N_{j,k}}\mu
_{j,k,2}+...),\,\, h\to 0
}
where the index $j$ labels the critical points and the leading
coefficients $\mu _{j,k}$ can be given explicitly in terms of the
Hessian of $V$ at the corresponding critical point. The values
$\mu _{j,k}$ are confined to a sector $\{ |{\rm arg\,}(z-\mu _{j,0})|\le
\theta_j \}$ for some $\theta_j \in [0,{\pi \over 2}[$,
and $\lambda _{j,0}$ is the eigenvalue with the smallest real part of
all the $\lambda _{j,k}$. (This comes formally
from a harmonic oscillator approximation.) Further $\mu _{j,k}\ne \mu
_{j,0}$ for $k\ne 0$ and the  asymptotic expansion \no{int.2} for
$\lambda _{j,0}$ contains only integer powers of $h$. See
\cite{HeSjSt} for more details. Using this result, as well as control
over the resolvent along suitable contours in the right half plane, the authors
were able to give asymptotic expansions for large times of $\exp\left(-tP/h\right)$,
that emphasize the role of the eigenvalues close to $0$ given
in \no{int.2}. Indeed, there are no others in a certain parabolic
neighborhood of the imaginary axis.

Moreover, we have  $\mu _{j,0}\ge 0$ with equality precisely when $U_j$ is
a local minimum of the potential, and in the case of such a minimum it
follows from the above results that
$\lambda _{j,0}$ is actually exponentially small.

\par In this paper we address the question of determining more
precisely the size of these exponentially small eigenvalues. In the
case when $V$ has precisely one local minimum, say $U_1$, and tends to
$+\infty $, when $x\to \infty $, we know that the corresponding
eigenvalue $\lambda _{1,0}$ is equal to zero (with the Maxwellian
$\exp \left(-(y^2/2+V(x))/h\right)$ as the corresponding eigenfunction) and that
this eigenvalue is separated from the other ones by a gap of size
$h$. This means that we have return to equilibrium with a speed that is
roughly 1. The situation becomes more complicated when there is
more than one local minimum. We are then in the presence of a
tunneling problem which is much more complicated than the
corresponding ones for the \scl{} Schr\"odinger operators since our
operators are non-elliptic. In principle one should be able to follow
the general approach of earlier works in the Schr\"odinger case as
\cite{HeSj}. However it seems that one necessarily runs into a
tunneling problem where the wave functions have to be studied also in
a neighborhood of some intermediate saddle
points of $V$, and as known from \cite{HeSj3} that can indeed
be done in the Schr\"odinger case with techniques that are
very useful in a variety of
problems. To carry out such an approach in the case of
Kramers-Fokker-Planck would require one to accumulate the difficulties of
non-resonant wells with the ones coming from the lack of ellipticity.
This seems to lead to considerations of degenerate non-symmetric Finsler
distances (see for instance \cite{BaChSh, LiNi}).

For the Witten Laplacian (see \cite{HeSj4}) we are also in the
presence of a tunneling problem with intermediate non-resonant wells
and in that case one could avoid the detailed study near the
non-resonant wells by studying directly the Witten complex as a
tunneling problem between critical points of neighboring
indices. More recently M. Klein, B. Helffer, and F. Nier \cite{HeKlNi} have used that
approach to study the exponentially small
non-vanishing eigenvalues of the Witten Laplacian.
Also in \cite{HeNi}, explicit estimates relating such small eigenvalues for the Witten
Laplacian and the Kramers-Fokker-Planck operator were established.

This relation with the Witten complex was strengthened further in the
works of J.~Tailleur, S.~Tanase-Nicola, J.~Kurchan \cite{TaTaKu} and
J.~M.~Bismut \cite{Bi}, who showed using respectively the languages of
supersymmetry and differential forms, that the Kramers-Fokker-Planck
operator can be viewed as a Witten Laplacian in degree 0 associated to
a certain non-semidefinite scalar product in the spaces of
differential forms. See also \cite{Le} for a quick introduction to the
differential form version of Bismut and \cite{BiLe}.

In the present paper, we use this supersymmetric approach. Our main
result, valid also for a class of more general operators, is that if a
certain weight function $\phi$ (which in the KFP-case is the function
$y^2/2+V(x)$) has precisely two local minima $U_\pm$ and an
intermediate saddle point $U_0$ then we can get a
complete asymptotic expansion for the second eigenvalue of the
corresponding Witten Laplacian (reducing to the KFP operator in the
special case). The logarithm of this eigenvalue is equal to $-2h\inv(\min
(\phi (U_0)-\phi (U_1),\phi (U_0)-\phi (U_{-1}))+o(1))$, but actually we do
have a complete asymptotic expansion. See Theorem \ref{2w1} for a
complete statement. It seems clear that this result can be somewhat
generalized but a more complete result
might require
exponential estimates and asymptotics for eigenfunctions also
far from the critical points. In our present approach we are able to
get such information for the eigenfunctions of the degree 0 operator
in the basin of attraction of each minimum and for the degree 1
Laplacian in a small neighborhood of the saddle point.

In most of the paper we work with a scalar real second order
non-elliptic operator, which is also non-\sa{}, and we were led to
reconsider some steps in
\cite{HeSjSt}. The plan of the paper is the following:

\smallskip
\par\noindent
In Section \ref{ap} we do some very simple and elementary exponential
estimates mainly designed to get the appropriate control near
infinity.
\smallskip
\par\noindent
In Section \ref{macr} we establish the $m$-accretivity for our
operators so that the step from a priori estimates to spectral
information becomes possible.
\smallskip
\par\noindent
In Section \ref{geo} we study certain auxiliary weights, somehow related to
escape functions in resonance theory (see \cite{HeSj2} and a large number
of more recent works) in connection with some dynamical conditions.
\smallskip
\par\noindent
In Section \ref{qw} we use those weights together with a machinery of \fop{}s
with complex phase in order to get phase space a priori estimates away
from the critical points. This section is perhaps technically the most
complicated one,
but the underlying ideas are now quite standard. Alternative methods
are certainly possible and we might return to this step in future works.
This section and the subsequent one are quite technical and
should not be studied in detail in the first reading.
\smallskip
\par\noindent
In Section \ref{co} we study the conjugation of our original
operator under the \fop{}s of the preceding section and in Section
\ref{app} we finally obtain the a priori estimates that we need.
\smallskip
\par\noindent
In Section \ref{as} it is now quite easy to get detailed asymptotic
results by adapting the methods of \cite{HeSjSt}.
\smallskip
\par\noindent
In Section \ref{exp} we show that the eigenfunction associated to
$\lambda _{j,0}$ for each (nondegenerate) critical point $U_j$ has the
form
$a(x;h)e^{-\phi (x)/h}$ in a neighborhood of that point, where $a(x;h)$ has
an asymptotic expansion in integer powers of $h$ with coefficients in
$C^\infty ({\rm neigh\,}(U_j))$ and $\phi =\phi _j$ is a smooth
function of the order of magnitude $|x-U_j|^2$. It would be very
interesting to extend such descriptions further away ``beyond caustics''.
\smallskip
\par\noindent
In Section \ref{ss} we review the supersymmetric approach of
\cite{TaTaKu, Bi} (see also \cite{Le}) and  establish various interesting links
between the dynamical conditions of Section \ref{geo} and old
results for non-\sa{} operators with double
characteristics (\cite{Sj}. This sheds additional light on some related
computations in \cite{HeSjSt}.
\smallskip
\par\noindent
Finally in Section \ref{2w}, we can put the various results together
and establish the precise exponential asymptotics of the spectral gap
between the first and the second eigenvalue (both real and the first
one being zero.)
\smallskip
\par We expect that the spectral results of the present paper will give rise to
precise asymptotics for the associated heat-evolution problem in the
limit of large times and we plan to treat that problem in a separate paper.

\par\noindent
\medskip
\par\noindent \bf Acknowledgements. \rm We are grateful to B. Helffer who
pointed out the work \cite{TaTaKu} to us, and for pointing out that an
earlier version of our condition (11.2) was not directly applicable to
the Kramers-Fokker-Planck operator. The research of the second author
is supported in part by the National Science Foundation under grant
DMS--0304970 and by an Alfred P. Sloan Research Fellowship. He is happy to acknowledge
the hospitality of \'Ecole Polytechnique and Universit\'e de Reims, where part of this work was done.

\section{A priori estimates}\label{ap}
\setcounter{equation}{0}

\par In this section we establish some simple a priori estimates which will
be important in Section \ref{co} and at other places. They illustrate
the technique of gaining ellipticity by means of exponential weights that
we shall later employ also in a micro-local setting.
\par Let $M$ denote either the space ${\bf R}^n$, or a smooth compact
$n$-dimensional \mfld{} equipped with a strictly positive smooth density
of integration $dx$. On $M$ we consider a second order \dop{}
\eekv{ap.1}
{
\hskip -5truemm P&=&\sum_{j,k=1}^n hD_{x_j}\circ b_{j,k}(x)\circ hD_{x_k}+{1\over
2}\sum_{j=1}^n (c_j(x)h\partial _{x_j}+h\partial _{x_j}\circ c_j(x))+p_0(x)
}
{&=&P_2+iP_1+P_0,\quad D_{x_j}={1\over i}{\partial \over \partial x_j},}
where the \coef{}s $b_{j,k}$, $c_j$, $p_0$ are assumed to be smooth and real, with
$b_{j,k}=b_{k,j}$. In the \mfld{} case, we use local coordinates such
that $dx=dx_1...dx_n$. To $P$ we associate the symbol in the \scl{}
sense,
\ekv{ap.2}
{p(x,\xi )=p_2(x,\xi )+ip_1(x,\xi )+p_0(x),}
\ekv{ap.3}
{p_2(x,\xi )=\sum_{j,k=1}^n b_{j,k}(x)\xi _j\xi _k,\ p_1(x,\xi
)=\sum_{j=1}^n c_j(x)\xi _j,}
so that $p_j(x,\xi )$ is a real-valued \pol{} in $\xi $, positively
homogeneous of degree $j$. (It is well-defined on $T^*M$ and coincides
with the Weyl symbol mod ${\cal O}(h^2)$ locally \ufly{}.) We assume that
\ekv{ap.4}
{p_2(x,\xi )\ge 0,\ p_0(x)\ge 0.}

\par In the case $M={\bf R}^n$, we impose the following growth conditions
at infinity:
\ekv{ap.5}
{
\partial _x^\alpha b_{j,k}(x)={\cal O}(1),\ \vert \alpha \vert \ge 0,
}
\ekv{ap.6}
{
\partial _x^\alpha c_j(x)={\cal O}(1),\ \vert \alpha \vert \ge 1,
}
\ekv{ap.7}
{
\partial _x^\alpha p_0(x)={\cal O}(1),\ \vert \alpha \vert \ge 2.
}

\par When discussing $P$
in the \op{} theoretical sense we will assume that it is the closure of
$P:{\cal S}(M)\to {\cal S}(M)$ as an un\bdd{} \op{} in $L^2(M)$. (When $M$
is compact, we identify the Schwartz space ${\cal S}(M)$ with $C^\infty (M)$.) Let ${\cal
D}(P)\subset L^2(M)$ be the domain.

\begin{lemma}\label{ap1} We have
\ekv{ap.8}
{
\int p_0(x) \vert u(x)\vert ^2dx+\int
\sum_{j,k}b_{j,k}(x)(hD_{x_j}u)(\overline{hD_{x_k}u})dx=\Re (Pu\vert u),
}
for all $u\in{\cal S}(M)$.
\end{lemma}
\begin{proof} Immediate by integration by parts.
\end{proof}

\par In the \mfld{} case, we view $B(x)=(b_{j,k}(x))$ as a positive
semi-definite matrix $T_x^*M\to T_xM$ and if we choose some smooth
Riemannian metric on $M$, we can view $B(x)$ as a map $T_x^*M\to T_x^*M$
and define $B(x)^{1\over 2}$ similarly. \no{ap.8} then becomes
\ekv{ap.9}
{
\Vert p_0(x)^{1\over 2}u\Vert ^2+\Vert B(x)^{1\over 2}hDu\Vert ^2=\Re
(Pu\vert u),
}
implying
\ekv{ap.10}
{\Vert p_0^{1\over 2}u\Vert +\Vert B^{1\over 2}hDu\Vert \le C_0(\Vert
Pu\Vert +\Vert u\Vert ).}
In particular, $u\in{\cal D}(P)\Rightarrow \Vert p_0^{1\over 2}u\Vert
+\Vert B^{1\over 2}hDu\Vert <\infty $.

\par Using the anti-\sa{} part $iP_1$ we shall obtain a similar estimate
where the averages of $p_0$ along the trajectories of
\ekv{ap.10.5}
{\nu (x,\partial _x)=\sum_1^n c_j(x)\partial _{x_j}}
will play a role.

\par In general, if $\psi (x)$ is a smooth real-valued \fu{}, the \op{}
\ekv{ap.11}
{
P_\psi :=e^{\psi /h}\circ P\circ e^{-\psi /h}
}
is of the same form as \no{ap.1} with new \coef{}s $c_j$, $p_0$ and the new
symbol
\eeeekv{ap.12}
{p_\psi (x,\xi )&=&p_2(x,\xi +i\psi '(x))+ip_1(x,\xi +i\psi '(x))+p_0(x)}
{&=& p_2(x,\xi )+i(p_1(x,\xi )+\partial _\xi p_2(x,\psi '(x))\cdot \xi
)}
{&&+p_0(x)-p_1(x,\psi '(x))-p_2(x,\psi '(x))}
{&=:& p_2(x,\xi )+ip_{1,\psi }(x,\xi )+p_{0,\psi}(x) }
In this section we choose $\psi $
very small and treat $p_2(x,\psi '(x))$
as a \pert{}. Notice that
\ekv{ap.13}
{p_1(x,\psi '(x))=\nu (x,\partial _x)\psi .}

\par Let $f(t)\in C^\infty ([0,\infty [;[0,3/2])$ be an increasing \fu{}
with $f(t)=t$ on $[0,1]$, $f(t)=3/2$ on $[2,\infty [$, $f(t)\le t$. Put
$f_\epsilon (t)=\epsilon f(t/\epsilon )$ and consider for $T_0>0$ fixed,
\ekv{ap.14}
{
\psi _\epsilon =\int k({t\over T_0})f_\epsilon \circ p_0\circ \exp
(t\nu )dt,
}
where
$$
k(t)=\cases{0,\ \vert t\vert \ge 1/2,\cr t+{1\over 2},\ -{1\over 2}\le
t<0,\cr -k(-t), \ 0<t\le {1\over 2}.}
$$
Then
\ekv{ap.15}
{
\nu (\psi _\epsilon )=f_\epsilon \circ p_0-\langle f_\epsilon \circ
p_0\rangle _{T_0},
}
where
\ekv{ap.16}
{
\langle f_\epsilon \circ p_0\rangle _{T_0}={1\over
T_0}\int_{-T_0/2}^{T_0/2}f_\epsilon \circ p_0\circ \exp (t\nu )dt
}
is the time $T_0$ average of $f_\epsilon \circ p_0$ along the integral
curves of $\nu $.
Clearly,
\ekv{ap.17}
{
\vert \psi _\epsilon \vert \le {T_0\over 4}{3\over 2}\epsilon ={3T_0\over
8}\epsilon .
}

\par From \no{ap.6} it is easy to see that $\Phi _t(x):=\exp (t\nu )(x)$ is
well-defined for all $t\in{\bf R}$, $x\in M$ and that

\ekv{ap.18}
{
\vert \partial _x^\alpha \Phi _t(x)\vert \le C(\alpha )e^{C_0\vert
\alpha \vert \vert t\vert },\ \alpha \in{\bf N}^n,\, \vert \alpha
\vert \ge 1 .}
In particular, $\vert \partial _x^\alpha \Phi _t(x)\vert \le C_{\alpha
,T_0}$, $\vert t\vert \le T_0$. On the other hand, since $p_0\ge 0$,
$p_0''={\cal O}(1)$ (by \no{ap.7}), we know that $\vert p_0'\vert \le
{\cal O}(1) p_0^{1\over 2}$ and this quantity is ${\cal O}(\epsilon
^{1/2})$ in the region where $0\le p_0\le 2\epsilon $. It follows that
$\partial ^\alpha (f_\epsilon \circ p_0)={\cal O}(\epsilon ^{1-\vert
\alpha \vert /2})$, $\alpha \in{\bf N}^n$ and together with \no{ap.18},
we get

\ekv{ap.19}
{
\partial _x^\alpha \psi _\epsilon ={\cal O}(\epsilon ^{1-\vert \alpha
\vert /2}),\ \forall \alpha \in{\bf N}^n, }
for every fixed $T_0$.

\par Then from \no{ap.12}, \no{ap.13}, \no{ap.15}, \no{ap.19} and the \uf{} \bdd{}ness
of the $b_{j,k}$, we get
\eekv{ap.21}{p_{0,\delta \psi_\epsilon }&=&p_0(x)-\delta f_\epsilon \circ p_0+\delta
\langle f_\epsilon \circ p_0\rangle _{T_0}-\delta ^2p_2(x,\psi '_\epsilon )
}
{&=&p_0(x)-\delta f_\epsilon \circ p_0+\delta
\langle f_\epsilon \circ p_0\rangle _{T_0}-{\cal O}(\delta ^2\epsilon ),
}
\ufly{} on ${\bf R}^n$ for $0\le \delta \le 1$. Using the properties of
$f$, we notice that $p_0(x)-\delta f_\epsilon \circ p_0\ge (1-\delta
)p_0$, so
\ekv{ap.22}
{
p_{0,\delta \psi_\epsilon  }\ge (1-\delta )p_0+\delta \langle f_\epsilon \circ
p_0\rangle _{T_0}-{\cal O}(\delta ^2\epsilon ).}

\par Since the \coef{}s of $P_{\delta \psi _\epsilon }$ grow at most
polynomially, Lemma \ref{ap1} can be applied and gives
\ekv{ap.23}
{
\int (p_{0,\delta \psi _\epsilon }(x)-\Re z)\vert u\vert ^2dx\le \Re
((P_{\delta \psi _\epsilon }-z)u\vert u),\ z\in{\bf C},\, \Re z\le
{\epsilon \over Z},
}
where $Z\gg 1$ is \indep{} of $\epsilon $.
(For Lemma \ref{ap1} we do not need that $p_0\ge 0$.) Let $\mu >0$
and rewrite \no{ap.23} as
\begin{eqnarray*}
\int \max (\widetilde{p}_{0,\psi },\mu )\vert u\vert ^2dx \le
\Re (((P_\psi -z)+(\mu -\widetilde{p}_{0,\psi })_+)u\vert u)\\
= \Re (\max (\widetilde{p}_{0,\psi },\mu )^{-{1\over 2}}(P_\psi -z +(\mu
-\widetilde{p}_{0,\psi })_+)u\vert \max (\widetilde{p}_{0,\psi },\mu )^{1\over 2}u).
\end{eqnarray*}
Here we write $\psi $ for $\delta \psi _\epsilon $ and
$\widetilde{p}_{0,\psi }=p_{0,\psi }-\Re z$.

\par Then using Cauchy-Schwarz,
$$
\Vert \max (\widetilde{p}_{0,\psi },\mu )^{1\over 2}u\Vert \le \Vert
\max (\widetilde{p}_{0,\psi },\mu )^{-{1\over 2}}(P_\psi -z)u\Vert
+\Vert {(\mu -\widetilde{p}_{0,\psi })_+\over \mu }\max
(\widetilde{p}_{0,\psi },\mu )^{1\over 2}u\Vert ,
$$
leading to
\ekv{ap.24}
{
\Vert \max (\widetilde{p}_{0,\psi },\mu )^{1\over 2}u\Vert \le \Vert
\max (\widetilde{p}_{0,\psi },\mu )^{-{1\over 2}}(P_\psi -z)u\Vert +\Vert
1_{\{ \widetilde{p}_{0,\psi }\le \mu \}}\max (\widetilde{p}_{0,\psi
},\mu )^{1\over 2}u\Vert .}
Notice that $1_{\{ \widetilde{p}_{0,\psi }\le \mu \} }\max
(\widetilde{p}_{0,\psi },\mu )^{1\over 2}=1_{\{\widetilde{p}_{0,\psi }\le
\mu \}}\mu ^{1\over 2}$.

\par We have
\begin{eqnarray*}
\widetilde{p}_{0,\psi }&=&p_0(x)-\delta f_\epsilon \circ p_0+\delta \langle
f_\epsilon \circ p\rangle _{T_0}-{\cal O}(\delta ^2\epsilon )-\Re z\\
&\ge & p_0(x)-\delta f_\epsilon \circ p_0+\delta \langle f_\epsilon \circ
p\rangle _{T_0}-{\cal O}(\delta ^2\epsilon )-{\epsilon \over Z}.
\end{eqnarray*}
Choose $Z=\delta ^{-2}$, so that
$$
\widetilde{p}_{0,\psi }
\ge  p_0(x)-\delta f_\epsilon \circ p_0+\delta \langle f_\epsilon \circ
p\rangle _{T_0}-{\cal O}(\delta ^2\epsilon ).
$$
Choose $\delta >0$
small enough (but \indep{} of $\epsilon $) so that ${\cal O}(\delta
^2\epsilon )\le \delta \epsilon /C_0$, where we shall fix $C_0$ \sufly{}
large. Then
$$
\widetilde{p}_{0,\psi }\ge (1-\delta )p_0(x)+\delta (\langle f_\epsilon
\circ p_0\rangle _{T_0}-{\epsilon \over C_0}).
$$

When $\langle f_\epsilon \circ p_0\rangle _{T_0}\ge {2\epsilon \over 3}$,
we get
$$
\widetilde{p}_{0,\psi }\ge (1-\delta )p_0(x)+{\delta \epsilon \over 2}.
$$
When $\langle f_\epsilon \circ p_0\rangle _{T_0}<{2\epsilon \over 3}$ we
have $\widetilde{p}_{0,\psi }\ge (1-\delta )p_0-{\delta \epsilon \over
C_0}$.

\par Choose
$$
\mu ={\delta \epsilon \over C_0}.
$$
Then
\ekv{ap.24.5}
{(1-\delta )p_0+\mu \le \widetilde{p}_{0,\psi }+2\mu \le 3\max
(\widetilde{p}_{0,\psi },\mu ).
}
Moreover, if $\widetilde{p}_{0,\psi }\le \mu $, then $\delta \langle
f_\epsilon \circ p_0\rangle_{T_0} -{\delta \epsilon \over C_0}\le {\delta
\epsilon \over C_0}$, so $\langle f_\epsilon \circ p_0\rangle _{T_0}\le
{2\epsilon \over C_0}$.
From \no{ap.24} and \no{ap.24.5} we then infer that
\eekv{ap.25}
{&&
\Vert ((1-\delta )p_0+{\delta \epsilon \over C_0})^{1\over
2}u\Vert \le} {&& 3\Vert ((1-\delta )p_0+{\delta \epsilon \over
C_0})^{-{1\over 2}}(P_\psi -z)u\Vert +
\sqrt{3}\sqrt{\delta \epsilon \over C_0}\Vert u\Vert _{\{ \langle
f_\epsilon \circ p_0\rangle _{T_0}\le
{2\epsilon \over C_0}
\}}, }
for $\Re z\le {\epsilon \over Z}$.

\par Here we can take $\epsilon =Mh$ with $M\gg 1$. Then by \no{ap.17} we
have
$$
\vert \psi /h\vert =\vert \delta \psi _\epsilon /h\vert \le C(M,T_0)\delta ,
$$
so there is a constant $C$ \indep{} of $h$ (but depending on $M,\delta $)
such that
$$
{1\over C}\le e^{\psi /h}\le C,
$$
From the discussion above, in particular \no{ap.25}, we get
\begin{prop}\label{ap2}
Let $P$ be of the form {\rm \no{ap.1}}, where $b_{j,k},c_j,p_0$ are smooth and real
and satisfy {\rm \no{ap.2}}--{\rm \no{ap.7}}. Define $\langle f_\epsilon \circ p_0
\rangle_{T_0}$ as in {\rm \no{ap.16}} with $f_\epsilon $ defined after {\rm \no{ap.13}}.
Then for every $C>0$, there exists $\widetilde{C}>0$ such that
\ekv{ap.26}
{
\Vert (p_0+h)^{1\over 2}u\Vert \le \widetilde{C}(\Vert (p_0+h)^{-{1\over
2}}(P-z)u\Vert +h^{1\over 2}\Vert u\Vert _{\{ \langle
f_{\widetilde{C}\epsilon} \circ p_0\rangle _{T_0}\le \widetilde{C}h \}}),
}
for $u\in {\cal S}$, $\Re z\le Ch$.
\end{prop}

Notice that \no{ap.8} implies that
\ekv{ap.27}
{
\Vert B^{1\over 2}hDu\Vert ^2 \le \Vert (p_0+h)^{-{1\over 2}}(P-z)u\Vert
\Vert (p_0+h)^{1\over 2}u\Vert +C\Vert h^{1\over 2}u\Vert ^2.
}

\section{From injectivity to the resolvent}\label{macr}
\setcounter{equation}{0}

\par Let $P=P_2+iP_1+P_0$ with symbol $p=p_2+ip_1+p_0$ be as in section
 \ref{ap}, so that we have \no{ap.8}
$$
\int p_0(x)\vert u(x)\vert ^2dx+(P_2u\vert u)=\Re (Pu\vert u),\ u\in {\cal
S}(M),
$$
leading to
$$
\vert \Re z\vert \Vert u\Vert ^2+\Vert p_0^{1\over 2}u\Vert ^2+(P_2u\vert
u)=\Re ((P-z)u\vert u)
\le \Vert (P-z)u\Vert \Vert u\Vert ,\ \Re z<0 .
$$
We get
$$
\Vert (\vert \Re z\vert +p_0)^{1\over 2}u\Vert ^2\le \Vert (\vert \Re
z\vert +p_0)^{-{1\over 2}}(P-z)u\Vert \Vert (\vert \Re z\vert
+p_0)^{1\over 2}u\Vert ,
$$
$$
\Vert (\vert \Re z\vert +p_0)^{1\over 2}u\Vert \le \Vert (\vert \Re
z\vert +p_0)^{-{1\over 2}}(P-z)u\Vert ,
$$
$$
\vert \Re z\vert ^{1\over 2}\Vert (\vert \Re z\vert +p_0)^{1\over 2}u\Vert \le \Vert
(P-z)u\Vert ,
$$
$$
\vert \Re z\vert\Vert u\Vert \le \Vert (P-z)u\Vert
$$
From this we get
$$
(P_2u\vert u)\le \Vert (P-z)u\Vert \Vert u\Vert \le {1\over \vert \Re
z\vert }\Vert (P-z)u\Vert ^2,
$$
and putting some of the estimates together,
\ekv{macr.1}
{
\vert \Re z\vert^2 \Vert u\Vert ^2+
\vert \Re z\vert \Vert p_0^{1\over 2}u\Vert ^2
+\vert \Re z\vert (P_2u\vert u))\le 2 \Vert (P-z)u\Vert ^2.
}

\par By $P$ we also denote the graph closure of $P:{\cal S}\to {\cal S}$.
From the estimates above we see that the range ${\cal R}(P-z)$ is
closed in $L^2$
when $\Re z<0$.
\begin{prop}\label{macr1}
$${\cal R}(P-z)=L^2,\ \Re z<0.$$
\end{prop}
\begin{proof}
It suffices to prove that ${\cal R}(P-z)=L^2$ for some $z$ with $\Re
z <0$,
because the a priori estimate then implies that $\Vert (P-z)\inv\Vert \le
\vert \Re z\vert \inv$
 and this fact extends by standard arguments to the whole left half plane.
 For the same $z$ it suffices to show that if $u\in L^2$ and
 $(P^*-\overline{z})u=0$ in the sense of distributions, then $u=0$. Now
 the formal adjoint $P^*=P_2-iP_1+P_0$ has the same properties as $P$, so in
 order to simplify the notations, we may just as well prove the
 corresponding fact for $P-z$ in stead of $P^*-\overline{z}$: There exists
 a $z$ with $\Re z<0$ such that if $u\in L^2$ and $(P-z)u=0$ in the sense
 of distributions, then $u=0$.

\par When $M={\bf R}^n$, let ${\rm Op}_h(q)$ denote the Weyl
 quantization of $q(x,h\xi )$ and put
$$
\Lambda _\epsilon ={\rm Op}_h(\langle \epsilon (x,\xi )\rangle ^{-N})
$$
where $N\ge 2$ is fixed and $\epsilon > 0$ is small and fixed.
(When $M$ is a compact \mfld{}, choose a Riemannian metric and put
$\Lambda _\epsilon =(1-\epsilon h^2\Delta )^{-N/2}$.)
Consider the
\e{}
$$
(P-z)u=v,\ u,v\in L^2
$$
Then $\Lambda _\epsilon (P-z)u=\Lambda _\epsilon v$, so
\ekv{macr.2}
{
(P-z)\Lambda _\epsilon u=\Lambda _\epsilon v+[P,\Lambda _\epsilon ]u.
}
Since $N\ge 2$, we can find a sequence $u_j\in{\cal S}$ such that
$\Lambda _\epsilon u_j\to \Lambda _\epsilon u$,
$p_0^{1\over 2}\Lambda _\epsilon u_j\to p_0^{1\over 2}\Lambda_\epsilon
u$, $\Lambda _\epsilon (P-z)u_j\to \Lambda _\epsilon v$,
$[P,\Lambda _\epsilon ]u_j\to [P,\Lambda _\epsilon ]u$ in $L^2$,
$(P_2\Lambda _\epsilon u_j\vert\Lambda _\epsilon  u_j)\to (P_2\Lambda _\epsilon u\vert \Lambda _\epsilon u)$.
This means that \no{macr.1} is applicable to \no{macr.2}, and we get
\ekv{macr.3}
{
\vert \Re z\vert ^2\Vert \Lambda _\epsilon u\Vert ^2+\vert \Re z\vert \Vert
p_0^{1\over 2}\Lambda _\epsilon u\Vert ^2+\vert \Re z\vert (P_2\Lambda
_\epsilon u\vert \Lambda _\epsilon u)\le 4(\Vert \Lambda _\epsilon v\Vert
^2+\Vert [P,\Lambda _\epsilon ]u\Vert ^2).
}

\par Consider $[P,\Lambda _\epsilon ]u=([P,\Lambda _\epsilon ]\Lambda
_\epsilon \inv )\Lambda _\epsilon u$. We study the \op{} $[P,\Lambda
_\epsilon ]\Lambda _\epsilon \inv$ and assume that $M={\bf R}^n$ in
order to fix the ideas. If $\rho =(x,\xi )$, we have
$$
\partial _\rho ^\alpha \langle \rho \rangle ^{-N}={\cal O}(1)\langle \rho
\rangle ^{-N-\vert \alpha \vert },
$$
so
$$
\partial _\rho ^\alpha \langle \epsilon \rho \rangle ^{-N}={\cal
O}(1)\langle \epsilon \rho \rangle ^{-N}({\epsilon \over \langle \epsilon
\rho \rangle })^{\vert \alpha \vert }={\cal O}(1)\langle \epsilon \rho
\rangle ^{-N}\langle \rho \rangle ^{-\vert \alpha \vert },
$$
\ufly{} \wrt{} $\epsilon $. From this we deduce that the symbol of
$[P,\Lambda _\epsilon ]$ is ${h\over i}\{ p,\Lambda _\epsilon \}+{\cal
O}_0(h^2\Lambda _\epsilon )$, using the same letters for \op{}s and their
symbols (except for $P, P_j$ where we already introduced a distinction by
using lower case letters for the symbols), and using the notation
${\cal O}_0(h^2\Lambda _\epsilon )$
for a
symbol $q$ satisfying $\partial _\rho ^\alpha q={\cal O}_\alpha
(h^2\Lambda _\epsilon )$, \ufly{} in $\epsilon $. Here
$$
\{p,\Lambda _\epsilon \} ={\cal O}_0(\Lambda _\epsilon (\rho )\langle
\rho \rangle ),
$$
so the symbol of $[\Lambda _\epsilon ,P]\Lambda _\epsilon \inv $ is
$$
{h\over i}{\{ p,\Lambda _\epsilon \}\over \Lambda _\epsilon }+{\cal
O}_0(h^2).
$$
Recall that $p(x,\xi )=p_0(x)+ip_1(x,\xi )+p_2(x,\xi )$. We get
$$
{\{ p_0,\Lambda _\epsilon \}\over \Lambda _\epsilon }={\cal O}_0({\langle
x\rangle \over \langle \rho \rangle })={\cal O}_0(1),
$$
$$
{\{ p_1,\Lambda _\epsilon \}\over \Lambda _\epsilon }={\cal O}_0({\langle
x\rangle \over \langle \rho \rangle })+{\cal O}_0({\langle
\xi \rangle \over \langle \rho \rangle })={\cal O}_0(1),
$$
$$
{\{ p_2,\Lambda _\epsilon \}\over \Lambda _\epsilon }=Q+{\cal O}_0(1),\
Q=-{\nabla _\xi \Lambda _\epsilon \over \Lambda _\epsilon }\cdot \nabla
_xp_2={\cal O}_0({1\over \langle \rho \rangle })\cdot \nabla _xp_2(\rho ).
$$

\par From these computations, we retain that
\ekv{macr.4}
{
\Vert [P,\Lambda _\epsilon ]\Lambda _\epsilon \inv w\Vert \le {\cal
O}(h)\Vert w\Vert +{\cal O}(h)\Vert Qw\Vert .  } The symbol $Q$ is real
and $\partial _\xi ^\alpha Q={\cal O}_0(\langle \xi  \rangle^{2-\vert \alpha
 \vert}/\langle \rho  \rangle )$, so
$$
\Vert Qw\Vert ^2=(Q^2w\vert w),
$$
where $Q\circ Q$
has the symbol $Q^2+{\cal O}_0(h^2)$.

\par Since $p_2\ge 0$, $\partial _x^2p_2={\cal O}(\xi ^2)$ we know that
$$
\vert \partial _xp_2\vert \le{\cal O}(1)\vert \xi \vert p_2(x,\xi )^{1\over
2},
$$
and we conclude that on the symbol level
$$
Q^2\le Cp_2.
$$
Hence by the semi-classical Fefferman-Phong inequality for operators with
symbols in the H\"ormander class $S^2_{1,0}$ (see for example Subsection 7.2
in \cite{HeSjSt} for a short review of the Weyl--H\"ormander calculus):
$$
\Vert Qw\Vert ^2\le C(P_2w\vert w)+{\cal O}(h^2)\Vert w\Vert ^2.
$$
Using this in \no{macr.4}, we get
\ekv{macr.5}
{
\Vert [P,\Lambda _\epsilon ]\Lambda _\epsilon \inv w\Vert ^2\le
Ch^2((P_2w\vert w)+\Vert w\Vert ^2).
}
we use this in \no{macr.3} with $[P,\Lambda _\epsilon ]u=[P,\Lambda
_\epsilon ]\Lambda _\epsilon \inv \Lambda _\epsilon u$ and get
\eekv{macr.6}
{
\vert \Re z\vert \Vert p_0^{1\over 2}\Lambda _\epsilon u\Vert ^2+(\vert
\Re z\vert ^2-Ch^2)\Vert \Lambda _\epsilon u\Vert ^2+(\vert \Re z\vert
-Ch^2)(P_2\Lambda _\epsilon u\vert \Lambda _\epsilon u )}{\le 4\Vert \Lambda
_\epsilon (P-z)u\Vert ^2.
}
So if $\Re z<-\sqrt{C}h$ and $(P-z)u=0$, $u\in L^2$, we have $u=0$.
\end{proof}
\begin{cor}\label{macr2}
The maximal closed extension $P_{\rm max}$ of $P$ (with domain given by
$\{ u\in L^2;\, Pu\in L^2\}$ coincides with the graph closure
(the minimal closed extension), already introduced.
\end{cor}
\begin{proof}
Let $z$ be fixed with $\Re z < -\sqrt{C}h$. Let $u,v\in L^2$ with $(P_{\rm
max}-z)u=v$. Denote by $P_{\rm min}$ the graph closure. Since ${\cal
R}(P_{\rm min}-z)=L^2$, there exists $\widetilde{u}\in {\cal D}(P_{\rm min})$
such that $(P_{\rm min}-z)\widetilde{u}=v$. Hence
$(P-z)(u-\widetilde{u})=0$ in the sense of distributions. We saw in the
proof of the proposition  that this implies that $u-\widetilde{u}=0$.
\end{proof}

\par This result can be extended to the following auxiliary problem:
Let $R_-:{\bf C}^N\to L^2$, $R_+:L^2\to {\bf C}^N$ be \bdd{} \op{}s
and assume that
$$
{\cal P}(z)=\pmatrix{P-z &R_-\cr R_+ &0}:{\cal D}(P)\times {\bf C}^N\to
L^2\times {\bf C}^N,
$$
is injective with
\ekv{macr.7}
{
\Vert u\Vert +\vert u_-\vert \le C(h)(\Vert v\Vert +\vert v_+\vert ),\
z\in\Omega (h),
}
whenever
\ekv{macr.8}
{
{\cal P}(z)\pmatrix{u\cr u_-}=\pmatrix{v\cr v_+}.
}
Here we assume that $\Omega (h)$ is open and connected, intersecting the
resolvent set $\rho (P)$ of $P$.

\par From \no{macr.7}, we see that ${\cal P}(z)$ is injective and has
closed range for every $z\in \Omega (h)$. When $z$ is also in the
resolvent set of $P$, we can rewrite \no{macr.8} as
\ekv{macr.9}
{
\pmatrix{1 &(P-z)\inv R_-\cr R_+ &0}\pmatrix{u\cr
u_-}=\pmatrix{(P-z)\inv v\cr v_+},
}
and the matrix appearing here is injective. On the other hand it is a
finite rank perturbation of the identity in $L^2\times {\bf C}^N$ and the
injectivity implies the bijectivity. It follows that ${\cal P}(z)$ is
bijective for $z\in\Omega (h)\cap\rho (P)$. Combining this with \no{macr.7}
and recalling that $\Omega $ is connected we see by a standard argument
that ${\cal P}(z)$ is bijective for all $z\in\Omega $ and that $\Vert
{\cal P}(z)\inv\Vert \le C(h)$.

\section{Geometric preparations}\label{geo}
\setcounter{equation}{0}

\par In this section we construct a special bounded weight that we shall
implement in Section \ref{qw} with the help of Fourier integral
operators with complex phase.

Let $P$
be the \op{} introduced in the beginning of Section \ref{ap}, satisfying
\no{ap.1}-\no{ap.7}.
The symbol
\ekv{geo.1}
{
\widetilde{p}(x,\xi )=p_0(x)+{p_2(x,\xi )\over \langle \xi \rangle ^2}
}
is non-negative and satisfies
\ekv{geo.2}
{
(\partial _x,\langle \xi \rangle \partial _\xi )^\alpha
\widetilde{p}={\cal O}(1),\ \vert \alpha \vert \ge 2,
}
locally \ufly{} \wrt{} $x$.
It follows that
\ekv{geo.3}
{
(\partial _x,\langle \xi \rangle  \partial _\xi )^\alpha \widetilde{p}={\cal
O}(\widetilde{p}^{1\over 2}),\ \vert \alpha \vert =1.
}

\par We now introduce a critical set associated to $p$.
\begin{hyp}\label{ny0}\rm
Assume
\ekv{ny.1}
{\hbox{The set }\{x \in M;\, p_0(x)=0,\ \nu (x,\partial _x)=0\} \hbox{ is
finite }=\{ x_1,...,x_N\} .}
\end{hyp}

\par
Let $\rho _j=(x_j,0)$ and put
\ekv{ny.2}
{
{\cal C}=\{ \rho _1,...,\rho _N\}
.}
Notice that $p_1,p_0,p_2,\widetilde{p}$ vanish to second order at each
$\rho _j$. Our weight will be of the form
\ekv{ny.3}
{
\psi _\epsilon =-\int k({t\over T_0})\widetilde{p}_\epsilon \circ \exp
(tH_{p_1})dt,\ 0<\epsilon \ll 1,
}
where $\widetilde{p}_\epsilon $ will be specified below and $k$ is the
same \fu{} as in Section \ref{ap}. Notice that
\ekv{ny.4}
{
H_{p_1}\psi _\epsilon =\langle \widetilde{p}_\epsilon \rangle
_{T_0}-\widetilde{p}_\epsilon ,
}
where we now write in general,
$$
\langle q\rangle _{T_0}={1\over T_0}\int_{-T_0/2}^{T_0/2}q\circ \exp
(tH_{p_1}) dt.
$$

\par Let $g(t)\in C^\infty ([0,+\infty [;[0,1])$ be a smooth decreasing
\fu{} with
\ekv{ny.5}
{
g(t)=\cases{1,\ 0\le t\le 1,\cr t\inv,\ t\ge 2.}
}
Notice that
$g^{(k)}(t)={\cal O}(\langle t\rangle ^{-1-k})$.
In a domain $0\le \vert \rho -\rho _j\vert <1/C$, $C\gg 1$, we put
\ekv{ny.6}
{
\widetilde{p}_\epsilon (\rho )=g({\vert \rho -\rho _j\vert ^2\over
\epsilon })\widetilde{p},
}
so that $\widetilde{p}_\epsilon (\rho )=\widetilde{p}(\rho )$ for $\vert
\rho -\rho _j\vert \le \sqrt{\epsilon }$. It is easy to show the symbol
estimates
\ekv{ny.7}
{
\partial ^\alpha \widetilde{p}_\epsilon ={\cal O}({\epsilon \over
(\epsilon +\vert \rho -\rho _j\vert ^2)^{\vert \alpha \vert /2}}),\ 0\le
\vert \rho -\rho _j\vert <1/C.}

\par Away from any fixed \neigh{} of ${\cal C}$, we simply get
$$
\partial ^\alpha \widetilde{p}_\epsilon ={\cal O}(\epsilon ).
$$
Let $\chi _j\in C_0^\infty ({\bf R}^{2n})$ be equal to 1 near $\rho _j$ and
have its support close to ${\cal C}$. To define $\widetilde{p}_\epsilon $
further out from ${\cal C}$, we put
\ekv{ny.8}
{
\widetilde{p}_\epsilon (\rho )=\sum_1^N\chi _j(\rho )g({\vert \rho -\rho
_j\vert ^2\over \epsilon })\widetilde{p}+\epsilon (1-\sum \chi _j(\rho
))\widetilde{p},
}
so that \no{ny.7} remains valid in a \neigh{} of ${\cal C}$, while
\ekv{ny.9}
{
\partial _x^\alpha \partial _\xi ^\beta \widetilde{p}_\epsilon ={\cal
O}(\epsilon \langle \xi \rangle ^{-\vert \beta \vert })={\cal O}(\epsilon
\langle \xi \rangle ^{-\vert \beta \vert })
}
outside any such region where $\vert x\vert $ is \bdd{}. We also have
\ekv{ny.10}
{
\partial _x^\alpha \partial _\xi ^\beta (\xi \cdot \partial _\xi
\widetilde{p}_\epsilon )={\cal O}(\epsilon \langle \xi \rangle ^{-2-\vert
\beta \vert })
}
in the same region. In fact, \no{ny.10} follows if we write for $\vert \xi
\vert \gg 1$:
$$
\xi \cdot \partial _\xi \widetilde{p}_\epsilon =\epsilon \xi \cdot
\partial _\xi {p_2\over \langle \xi \rangle ^2}=\epsilon \xi \cdot
\partial _\xi [ {p_2\over \vert \xi \vert ^2}(1-{1\over \langle \xi
\rangle ^2})]=-\epsilon {p_2\over \vert \xi \vert ^2}\xi \cdot \partial
_\xi ({1\over \langle \xi \rangle ^2}).
$$

\par From \no{ny.8} we see that $\widetilde{p}_\epsilon (\rho )=\epsilon
\widetilde{p}$ when ${\rm dist\,}(\rho ,{\cal C})\ge 1/C$ and in
particular for $\rho =(x,\xi )$ with $\vert x-x_j\vert \ge 1/C$, $\forall
j$. Further away from $\{ x_1,...,x_N\} $ we want to make
$\widetilde{p}_\epsilon (x,\xi )$ \indep{} of $\xi $, so we replace
$\widetilde{p}=p_0(x)+p_2(x,\xi )/\langle \xi \rangle ^2$ in $\no{ny.8}$ by
\ekv{ny.11}
{
\widetilde{p}_{\rm new}(x,\xi )=p_0(x)+{\chi (x)\over \langle \xi \rangle
^2}p_2(x,\xi ),
}
where $\chi \in C_0^\infty (M;[0,1])$ is equal to 1 near $\{
x_1,...,x_N\}$.

\par Even further out (in the case when $M={\bf R}^n$) we want to avoid
\pb{}s caused by $p_0$ being large, so when $\vert x\vert \gg 1$ we want to
replace $\widetilde{p}_\epsilon =\epsilon p_0$ there by $f_\epsilon (p_0)$,
where $f_\epsilon (t)=\epsilon f(t/\epsilon )$ is the \fu{} introduced in
Section \ref{ap}. Thus with a new cutoff \fu{} $\chi _{\rm new}\in
C_0^\infty (M;[0,1])$ being equal to 1 in a large \neigh{} of ${\rm
supp\,}\chi $, we get the final choice of $\widetilde{p}_\epsilon $:
\eekv{ny.12}
{
\widetilde{p}_\epsilon (\rho )=\chi _{\rm new}(x)(\sum_1^N \chi _j(\rho
)g({\vert \rho -\rho _j\vert ^2\over \epsilon })\widetilde{p}+\epsilon
(1-\sum \chi _j(\rho ))\widetilde{p}_{\rm new})
}{+(1-\chi _{\rm
new}(x))f_\epsilon (p_0).}
Notice that by construction
\ekv{ny.12.5}
{\widetilde{p}_\epsilon \le \widetilde{p}\le p_0+p_2.}
We also get
\begin{prop}\label{ny1}
We have {\rm \no{ny.7}} near ${\cal C}$ and {\rm \no{ny.9}}, {\rm \no{ny.10}} in any closed
subset of $T^*M$ disjoint from ${\cal C}$ and ${\rm supp\,}(1-\chi _{\rm
new})\times {\bf R}^n$, where $\chi _{\rm new}$ is the wider cutoff in
{\rm \no{ny.12}}. Over ${\rm supp\,}(1-\chi _{\rm new})$ (where
$\widetilde{p}_\epsilon $ only depends on $x$) we have $\partial _x^\alpha
\widetilde{p}_\epsilon ={\cal O}(\epsilon ^{1-\vert \alpha \vert /2})$.
\end{prop}

From the definition \no{ny.3}, we see that $\psi _\epsilon $ satisfies the
same estimates as $\widetilde{p}_\epsilon $ and will only depend on $x$
for $x$
outside a \neigh{} $\pi _x({\cal C})$ (that we can choose as small as we
like) and that the region where we only have  $\partial _x^\alpha \psi
_\epsilon (x)={\cal O}(\epsilon ^{1-\vert \alpha \vert /2})$ can be any
\neigh{} of infinity. This follows from the fact that our various symbol
estimates (as well as the $\xi $-independence)  are conserved by the flow
of $H_{p_1}$, as for \no{ny.10}, we here also use that $H_{p_1}$
commutes with dilations in $\xi $.

\par It also follows from the construction that
\ekv{ny.13}
{
\widetilde{p}_\epsilon (\rho )=\widetilde{p}(\rho ),\hbox{ when }\vert
\rho -\rho _j\vert \le\sqrt{\epsilon },
}
\ekv{ny.14}
{
\widetilde{p}_\epsilon (\rho )\backsim {\epsilon \over \vert \rho -\rho
_j\vert ^2}\widetilde{p},\hbox{ when }\sqrt{\epsilon }\le \vert \rho -\rho
_j\vert \le 1/C,
}
\ekv{ny.15}
{
\widetilde{p}_\epsilon (\rho )\backsim \epsilon \widetilde{p}_{\rm new},
\hbox{ when }{\rm dist\,}(\rho ,{\cal C})\ge 1/C,\ \vert \rho
_x\vert \le C,
}
\ekv{ny.16}
{
\widetilde{p}_\epsilon (\rho )=\epsilon f_\epsilon (p_0),\hbox{ when
}\vert \rho _x\vert \ge \widetilde{C}.
}

\par We introduce the following dynamical conditions where $T_0>0$ is fixed:
\begin{hyp}\label{ny2}\rm
\ekv{ny.17}
{
\hbox{Near each }\rho _j\hbox{ we have }\langle \widetilde{p}\rangle
_{T_0}\ge {1\over C}\vert \rho -\rho _j\vert ^2,
}
\eekv{ny.18}
{
\hbox{In any set }\vert x \vert \le C,\ {\rm dist\,}(\rho ,{\cal C})\ge
{1\over C},
}{\hbox{we have } \langle \widetilde{p}\rangle _{T_0}(\rho )\ge {1\over
\widetilde{C}(C)},\ \widetilde{C}(C)>0.
}
When $M={\bf R}^n$ we also need a modified dynamical assumption
\eekv{ny.20}
{
\forall\hbox{ \neigh{} }U\hbox{ of }\pi _x{\cal C}, \hbox{ and }\forall \
x\in {\bf R}^n \sm U,\ \exists C>0,
}
{
{\rm meas\,}(\{ t\in [-{T_0\over 2},{T_0\over 2}];\, p_0(\exp t\nu
(x))\ge {1\over C}\})\ge {1\over C}.
}
When $M$ is compact, we just assume \no{ny.17}, \no{ny.18}, where it is
understood that the estimate in \no{ny.18} should hold for all $\rho
\in T^*M$ of distance $\ge 1/C$ from ${\cal C}$.
\end{hyp}
\par Notice that in the region ${\rm dist\,}(x,\pi _x({\cal C}))\ge 1/C$, $\vert
x\vert \le C$, \no{ny.18} is equivalent (up to a change of
$\widetilde{C}(C)$) to
\ekv{ny.19}
{
\langle p_0\rangle _{T_0,\nu }(x)\ge {1\over \widetilde{C}(C)},\
\widetilde{C}(C)>0,
}
where
$$
\langle q\rangle _{T_0,\nu }={1\over T_0}\int_{-T_0/2}^{T_0/2}q\circ \exp
(t\nu )dt.
$$
This follows from the fact that $\pi _x\circ \exp tH_{p_1}=\exp (t\nu
)\circ \pi _x$, where $\pi _x((x,\xi )):=x$ and the fact that $\{ \xi =0\}$ is invariant under the
$H_{p_1}$ flow.

\par Actually, if we only assume \no{ny.17}Ê, then we see that $\langle p_2\rangle_{T_0}(x_j,\xi )\backsim \xi ^2$, and if we also assume \no{ny.19}, we recover \no{ny.18}, first over a \neigh{} of each $x_j$ and then again by \no{ny.19} over any \bdd{} set in $M$.

\par Assuming the dynamical conditions \no{ny.17}, \no{ny.18}, (\no{ny.19},
\no{ny.20}), we see that
\ekv{ny.21}
{
\langle \widetilde{p}_\epsilon \rangle _{T_0}\backsim {\rm dist\,}(\rho
,{\cal C})^2,\ {\rm dist\,}(\rho ,{\cal C})\le \sqrt{\epsilon },
}
\ekv{ny.22}
{
\langle \widetilde{p}_\epsilon \rangle _{T_0}\backsim \epsilon ,\ {\rm
dist\,}(\rho ,{\cal C})>\sqrt{\epsilon }.
}

\par From the estimates on $\psi _\epsilon $ in the various regions that we
mentioned after Proposition \ref{ny1} we shall often only retain that
\ekv{geo.8}
{
\partial _x^\alpha \partial _\xi ^\beta \psi _\epsilon (x,\xi )={\cal
O}(\epsilon ^{1-\vert \alpha +\beta \vert /2}\langle \xi \rangle ^{-\vert
\beta \vert }),
}
and we write this for short as
$$
\psi _\epsilon =\widetilde{{\cal O}}(\epsilon ).
$$
Similarly, we have in view of \no{ny.10}, and \no{ny.7} also valid for
$\psi _\epsilon $, that
\ekv{geo.10.5}
{
\xi \cdot \partial _\xi \psi _\epsilon =\widetilde{{\cal O}}(\epsilon
\langle \xi \rangle ^{-2}).
}

\par From~\cite{MeSj2} let us recall that if $f$ is a $C^1$ function locally defined on complexified
phase space, then at every point where $\overline{\partial }f=0$, we have
\ekv{geo.16}
{
\widehat{H}_f=H_{\Re f}^{\Re \sigma }=H_{\Im f}^{\Im \sigma },\
\widehat{H}_{if}=J\widehat{H}_f=H_{-\Im f}^{\Re \sigma }=H_{\Re f}^{\Im
\sigma },
}
where
$$H_f=\sum ({\partial f\over \partial \xi _j}{\partial \over \partial
x_j} - {\partial f\over \partial x _j}{\partial \over \partial \xi _j})$$
is the complex (1,0) Hamilton field of $f$ \wrt{} the complex
symplectic form $\sigma =\sum d\xi _j\wedge dx_j$ and
$\widehat{H}_f=H_f+\overline{H_f}$ is the real \vf{} which acts as $H_f$
on \hol{} \fu{}s. $H_{\Im f}^{\Im \sigma }$ denotes the real Hamilton
field of $\Im f$ \wrt{} the real symplectic form $\Im \sigma $, and
similarly for the other Hamilton fields appearing in \no{geo.16}. As
usual, $J$ denotes the action on tangent vectors induced by multiplication
by $i$. (When $M$ is compact we may assume without loss of generality
that $M$ is real-analytic.)

\par Assume first that $M={\bf R}^n$ and put $\Lambda _0=T^*{\bf R}^n={\bf
R}^{2n}$,
\ekv{geo.17}
{
\Lambda _s=\{ \rho +isH_{\psi _\epsilon }(\rho );\, \rho \in\Lambda _0\},
0\le s\ll 1.
}
If we extend $\psi _\epsilon $ to be a \fu{} on the complexified phase
space ${\bf C}^{2n}$, by setting
\ekv{geo.18}
{
\psi _\epsilon (\rho )=\psi _\epsilon (\Re \rho ),
}
then we have the equivalent description of $\Lambda _s$ as
\ekv{geo.19}
{
\Lambda _s=\exp (sH_{\psi _\epsilon }^{\Im \sigma })(\Lambda _0).
}
It follows that $\Lambda _s$ is an I-Lagrangian \mfld{}, i.e. a \mfld{} which
is Lagrangian for $\Im \sigma $.

\par Now we would like to parametrize $\Lambda _s$ by means of a \ctf{}
$\kappa (s):\Lambda _0\to \Lambda _s$ for the symplectic form $\Re \sigma $
(implying that $\Lambda _s$ like $\Lambda _0$ is an IR-\mfld, i.e. a \mfld{}
which is Lagrangian for $\Im \sigma $ and symplectic for the restriction
of $\Re \sigma $). For $0\le s\ll 1$, let $\widehat{\psi }_{\epsilon
,s}=\widehat{\psi }(s,\rho )$ be the unique \fu{} which is affine in $\Im
\rho $ and satisfies
\ekv{geo.20}
{
\overline{\partial }_\rho \widehat{\psi }_{\epsilon ,s}=0\hbox{ on }\Lambda
_s,\ \widehat{\psi }_{\epsilon ,s}=\psi _\epsilon \hbox{ on }\Lambda _s.
}
Recalling that $\psi _\epsilon (\rho )=\psi _\epsilon (\Re \rho )$, we get
with the notation $(\psi _\epsilon )'_\xi ={\partial \over \partial \Re
\xi }\psi _\epsilon $, $(\psi _\epsilon )'_x ={\partial \over \partial \Re x }
\psi _\epsilon $ to indicate derivatives in the real directions,
\ekv{geo.21}
{
\widehat{\psi }_{\epsilon ,s}(x,\xi )=\psi _\epsilon (x,\xi )+
ia_x(\Re (x,\xi ))
\cdot (\Im x-s(\psi _\epsilon )'_\xi
)+ia_\xi (\Re (x,\xi ))\cdot (\Im \xi
+s( \psi _\epsilon )' _x),
}
where
\ekv{geo.21.1}
{
(1-is{\trans F}_{\psi _\epsilon }(\Re (x,\xi )))\pmatrix {a_x\cr a_\xi
}=\pmatrix {(\psi _\epsilon )'_x \cr (\psi _\epsilon )'_\xi  }
}
and
\ekv{geo.21.2}
{F_{\psi _\epsilon }=\pmatrix{(\psi _\epsilon )''_{\xi x}
&(\psi _\epsilon)''_{\xi \xi }\cr -( \psi
_\epsilon )''_{xx} & -(\psi _\epsilon)''_{x\xi }} }
is the fundamental matrix of $\psi _\epsilon $. It also follows from the
construction that
\ekv{geo.21.3}
{
a_x=\widetilde{{\cal O}}(\epsilon ^{1\over 2}),\ a_\xi =\widetilde{{\cal
O}}(\epsilon ^{1\over 2}\langle \xi \rangle ^{-1}).
}
{}

\par Since $\widehat{\psi }_{\epsilon ,s}=\psi _\epsilon $ on $\Lambda _s$,
we know that $H^{\Im \sigma }_{\Re\widehat{\psi }_{\epsilon ,s}}- H^{\Im
\sigma }_{\psi _\epsilon }$ is tangent to $\Lambda _s$ for every $s$.  We
can therefore define $\kappa (s):\Lambda _0\to \Lambda _s$, by
\ekv{geo.22}
{
{d\over ds}\kappa (s)(\rho )=H_{\Re\widehat{\psi }_{\epsilon ,s}}^{\Im
\sigma }(\kappa (s)(\rho )),\ \rho \in\Lambda _0.
}
The second relation in \no{geo.16} implies that $H_{\Re \widehat{\psi
}_{\epsilon ,s}}^{\Im \sigma }=H_{-\Im \widehat{\psi }_{\epsilon
,s}}^{\Re \sigma }$ on $\Lambda _s$ (given that $\overline{\partial
}\widehat{\psi }_{\epsilon ,s}=0$ there) and hence $\kappa (s)$ is
symplectic for $\Re \sigma $: $\kappa (s)^*(\Re \sigma )=\Re \sigma $. Notice
that \no{geo.22} again shows that $\Lambda _s$ is an I-Lagrangian manifold.
\medskip

\par A priori, $\kappa (s)(\rho )$ is well-defined only for $s\ge 0$ small
enough \it depending on \rm $\rho $, but we shall next derive symbol
estimates for $\widehat{\psi }_{\epsilon ,s}$ and $\kappa (s)$ that will
imply that $\kappa (s)(\rho )$ is indeed well-defined for $s$ small enough
\it \indep{}ly of \rm $\rho $. Assume that we work
near a point $(x_0,r\xi _0)$, $x_0\in{\bf R}^n$, $\xi _0\in S^{n-1}$, $r\ge
0$ and replace $(x,\xi )$ by $(\widetilde{x},\widetilde{\xi })$, where
\ekv{geo.23}
{
x=x_0+\sqrt{\epsilon }\widetilde{x},\ \xi =r\xi _0+\langle r\rangle
\sqrt{\epsilon }\widetilde{\xi }.
}
Define $\widetilde{\psi }=\widetilde{\psi }_{\epsilon ,x_0,r,\xi _0}$ by
$\psi _\epsilon (x,\xi )=\epsilon \widetilde{\psi
}(\widetilde{x},\widetilde{\xi })$, so that
\ekv{geo.24}
{
\widetilde{\psi }(\widetilde{x},\widetilde{\xi })={1\over \epsilon }\psi
_\epsilon (x_0+\sqrt{\epsilon }\widetilde{x},r\xi _0+\langle r \rangle
\sqrt{\epsilon }\widetilde{\xi }),
}
\ekv{geo.25}
{
\partial _{\widetilde{x}}^\alpha \partial _{\widetilde{\xi }}^\beta
\widetilde{\psi }={\cal O}(1), \forall \alpha ,\beta ,\hbox{ for }\vert
(\widetilde{x},\widetilde{\xi })\vert \le {\rm Const.}
}
In the new coordinates $\Lambda _s$ takes the form
\ekv{geo.25.5}
{
\Im (\widetilde{x},\widetilde{\xi })={s\over \langle r\rangle
}H_{\widetilde{\psi }}(\Re (\widetilde{x},\widetilde{\xi })).
}

\par The scaling and the construction of $\widehat{\psi }_{\epsilon ,s}$
commute so if we put
$$
\widetilde{\widehat{\psi }}_{\epsilon ,s}={1\over \epsilon }\widehat{\psi
}_{\epsilon ,s}(x_0+\sqrt{\epsilon }\widetilde{x},r\xi _0+\langle r\rangle
\sqrt{\epsilon }\widetilde{\xi }),
$$
then \no{geo.25} still holds. In particular
$\widehat{\psi }_{\epsilon
,s}$ satisfies the same estimates on $\Lambda _s$ as $\psi _\epsilon $ in
\no{geo.8}, now with $\partial _x,\partial _\xi $ replaced by $(\partial
_{\Re x},\partial _{\Im x}),(\partial _{\Re \xi },\partial _{\Im \xi })$.

\par The transformation $\kappa (s)$ can be scaled similarly and the
scaling commutes with \no{geo.22} up to a factor $\langle r\rangle $:
If we put
$$
\kappa (s) (\rho )=(x(s),\xi (s))=(x_0+\sqrt{\epsilon
}\widetilde{x}(s),r\xi _0+\langle r\rangle \sqrt{\epsilon
}\widetilde{\xi }(s)),
$$
then from \no{geo.22}, \no{geo.16} we get
\eekv{geo.26}
{\partial _s\widetilde{x}(s)&=&{i\over \langle r\rangle }\partial
_{\widetilde{\xi }}\widetilde{\widehat{\psi }}_{\epsilon
,s}(\widetilde{x}(s),\widetilde{\xi }(s)),}
{\partial _s\widetilde{\xi }(s)&=&-{i\over \langle r\rangle }\partial
_{\widetilde{x }}\widetilde{\widehat{\psi }}_{\epsilon
,s}(\widetilde{x}(s),\widetilde{\xi }(s)).}
We conclude that with $(x(0),\xi (0))=(y,\eta )=(x_0+\sqrt{\epsilon }\widetilde{y},r\xi _0+\langle
r\rangle \sqrt{\epsilon }\widetilde{\eta })$,
\ekv{geo.27}
{
\partial _s^k\partial _{\widetilde{y}}^\alpha \partial
_{\widetilde{\eta }}^\beta (\widetilde{x},\widetilde{\xi })={\cal O}(\langle
r\rangle ^{-k}),
}
and hence for $(x,\xi )=\kappa (s)(y,\eta )$:
\eekv{geo.28}
{
\partial _s^k\partial _y^\alpha \partial _\eta ^\beta x&=&{\cal
O}(\epsilon ^{{1\over 2}-{\vert \alpha +\beta \vert \over 2}}\langle \eta
\rangle ^{-k-\vert \beta \vert }),
}
{\partial _s^k\partial _y^\alpha \partial _\eta ^\beta \xi &=&{\cal
O}(\epsilon ^{{1\over 2}-{\vert \alpha +\beta \vert \over 2}}\langle \eta
\rangle ^{1-k-\vert \beta \vert }),}
when $k+\vert \alpha +\beta  \vert\ge 1$.

\par Notice that the \rhs{}s in \no{geo.26} reduce to $i\langle r\rangle
^{-1}\partial _{\widetilde{\xi }}\widetilde{\psi }$ and $-i\langle r\rangle
^{-1}\partial _{\widetilde{x}}\widetilde{\psi }$ when $s=0$, where the
derivatives are taken in the real directions. The flow in \no{geo.28} is
therefore tangent to the one in \no{geo.19} at $s=0$ and we get
\begin{eqnarray*}
\kappa (s)(y,\eta )&=&(y,\eta )+isH_{\psi _\epsilon }(y,\eta )+(z,\zeta ),
\\
\partial _y^\alpha \partial _\eta ^\beta z&=&{\cal O}(\epsilon ^{{1\over
2}-{\vert \alpha +\beta \vert \over 2}}s^2\langle \eta \rangle ^{-2-\vert
\beta \vert }),\\
\partial _y^\alpha \partial _\eta ^\beta \zeta &=&{\cal O}(\epsilon ^{{1\over
2}-{\vert \alpha +\beta \vert \over 2}}s^2\langle \eta \rangle ^{-1-\vert
\beta \vert }).
\end{eqnarray*}

\par From \no{geo.28}
and the subsequent remark, we have for $\kappa (\delta )(x,\xi )=
\kappa (\delta ,x,\xi )$.
\ekv{geo.43}
{
\kappa (\delta ,x,\xi )=(x,\xi )+i\delta H_{\psi _\epsilon }(x,\xi
)+\epsilon ^{1\over 2}\delta ^2(\widetilde{{\cal O}}({1\over \langle \xi \rangle
^2}),\widetilde{{\cal O}}({1\over \langle \xi \rangle })).
}
 Recalling how $\kappa (\delta )$
was constructed we conclude that
\ekv{geo.44}
{
\kappa (\delta ,x,\xi )=(\widetilde{x},\widetilde{\xi })+i\delta H_{\psi
_\epsilon }(\widetilde{x},\widetilde{\xi }), } where
$(\widetilde{x},\widetilde{\xi })$ is real and
$(\widetilde{x},\widetilde{\xi })-(x,\xi )=\epsilon ^{1\over 2}\delta
^2(\widetilde{{\cal O}}({1\over \langle \xi \rangle ^2}),
\widetilde{{\cal O}}({1\over \langle \xi
\rangle }))$. Put $\alpha _\delta
(x,\xi )=(\widetilde{x},\widetilde{\xi })$.

\par The essential part of the discussion above took part near the
points of ${\cal C}$. In that region the discussion is the same in the
case when $M$ is compact.

\section{Quantization of weights}\label{qw}
\setcounter{equation}{0}

\par We will follow \cite{HeSj2}, \cite{Sj2} with one modification;
instead of analyticity we will use that our weights are "moderate" allowing
us to use almost \hol{} extensions. Another minor difference is that
we shall not use FBI-transforms explicitly, but rather rely on
certain \fop{}s with complex phase. We will assume that $M={\bf R}^n$
for simplicity, but as in the preceding section the essential part of
the work will take place near ${\cal C}$ and here there is no
difference between the case $M={\bf R}^n$ and the case when $M$ is
compact.

\par As a first step towards introducing some \fop{}s we shall study the
\fu{} $h(y,\eta )$ on $T^*{\bf R}^n$, given by
\ekv{qw.21}
{
\kappa (\delta )^*(\xi \cdot dx)-\eta \cdot dy=dh.
}
Recall here that
\ekv{qw.22}
{
{d\over ds}\kappa (s)\rho =H_{i\widehat{\psi }_{\epsilon ,s}}(\kappa
(s)\rho ),\ \rho \in\Lambda _0=T^*{\bf R}^n,
}
where $\widehat{\psi }_{\epsilon ,s}$ is given by \no{geo.21}, satisfying
\no{geo.20}. Of course, \no{qw.22} remains unchanged if we replace
$\widehat{\psi }_{\epsilon ,s}$ by an almost \hol{} extension from $\Lambda
_s=\Lambda _{s\psi _\epsilon }$. Now using Cartan's formula, we get
\begin{eqnarray*}
{d\over ds}\kappa (s)^*(\xi \cdot dx)&=&\kappa (s)^*{\cal
L}_{iH_{\widehat{\psi }_{\epsilon ,s}}}(\xi \cdot dx)\\
&=&
i\kappa (s)^*(H_{\widehat{\psi }_{\epsilon ,s}}\rfloor d(\xi \cdot
dx)+d(H_{\widehat{\psi }_{\epsilon ,s}}\rfloor \xi \cdot dx))\\
&=& i\kappa (s)^*(H_{\widehat{\psi }_{\epsilon ,s}}\rfloor \sigma +d(\xi
\cdot {\partial \widehat{\psi }_{\epsilon ,s}\over \partial \xi }))\\
&=&i\kappa (s)^*(d(\xi \cdot {\partial \widehat{\psi }_{\epsilon ,s}\over
\partial \xi})-d\widehat{\psi }_{\epsilon ,s})\\
&=& id\kappa (s)^*(\xi \cdot {\partial \over \partial \xi }\widehat{\psi
}_{\epsilon ,s}-\widehat{\psi }_{\epsilon ,s}).
\end{eqnarray*}
Thus, we can take
\ekv{qw.23}
{
h=i\int_0^\delta  (\xi \cdot {\partial \over \partial \xi }\widehat{\psi
}_{\epsilon ,s}-\widehat{\psi }_{\epsilon ,s})\circ \kappa (s)ds.
}
On $\Lambda _s$ we have ${\partial \widehat{\psi }_{\epsilon ,s}
\over \partial \xi }={\partial \widehat{\psi }_{\epsilon ,s}
\over \partial \Re\xi }$ and we recall that $\psi _\epsilon (x,\xi )=
\psi _\epsilon (\Re (x,\xi )).$
Using \no{geo.21}, we get on $\Lambda _s$:
\ekv{qw.24}
{
\xi \cdot {\partial \over \partial \xi }\widehat{\psi }_{\epsilon ,s}=\xi
\cdot {\partial \over \partial \Re \xi }\psi _\epsilon (x,\xi )-
is\langle a_x(\Re
(x,\xi ))\vert (\psi _\epsilon )''_{\xi\xi} \xi
\rangle +is\langle a_\xi (\Re (x,\xi ))\vert (\psi _\epsilon )''_{
x \xi}\xi \rangle ,
}
where again we use the notation $(\cdot )'_\xi $, $(\cdot )''_{x\xi }$ etc to
indicate derivatives in the real directions.

\par We now estimate the terms in the \rhs{} of \no{qw.24}. On $\Lambda
_s$ we have $\Im \xi =s(\psi _\epsilon )'_{ \xi }$, so
$$
\xi \cdot {\partial \over \partial \Re\xi }\psi _\epsilon =\Re \xi \cdot
(\psi _\epsilon )'_\xi  -is (\psi _\epsilon )'_x\cdot (\psi _\epsilon )'_\xi
.
$$
Hence by \no{geo.8}, \no{geo.10.5},
\ekv{geo.25'}
{
\xi \cdot {\partial \over \partial\Re \xi }\psi _\epsilon =\widetilde{{\cal
O}}(\epsilon \langle \xi \rangle ^{-1}).
}
Next, look at
$$
\langle a_x\vert  (\psi _\epsilon )''_{\xi \xi}\xi
\rangle =\langle a_x\vert (\psi _\epsilon )''_{\xi \xi}
\Re \xi \rangle -is\langle a_x\vert (\psi _\epsilon )''
_{\xi \xi }( \psi _\epsilon )_x' \rangle .
$$
The last term is $\widetilde{{\cal O}}(\epsilon \langle \xi \rangle ^{-2})$
by \no{geo.8}, \no{geo.21.3}. The first term is equal to
$$
\langle a_x\vert (\psi _\epsilon )''_{\xi \xi }\Re \xi
\rangle =a_x\cdot {\partial \over \partial \Re\xi }((\psi _\epsilon )'
_\xi \cdot \Re \xi )-(\psi _\epsilon )'_\xi\cdot a_x.
$$ From \no{geo.10.5}, \no{geo.8}, \no{geo.21.3} this is $\widetilde{{\cal
O}}(\epsilon
\langle \xi \rangle \inv )$, so
\ekv{qw.26}
{
\langle a_x\vert (\psi _\epsilon )''_{\xi \xi }\xi
\rangle =\widetilde{{\cal O}}(\epsilon \langle \Re \xi \rangle \inv ).
}
Similarly,
\begin{eqnarray*}
\langle a_\xi \vert  (\psi _\epsilon )''_{x\xi}
\xi
\rangle &=&\langle a_\xi \vert (\psi _\epsilon )''_{x\xi }\Re \xi \rangle -
is\langle a_\xi \vert (\psi _\epsilon )''_{x\xi }
( \psi _\epsilon )'_x\rangle \\
&=&\langle a_\xi \vert (\psi _\epsilon )''_{x\xi}\Re \xi \rangle+\widetilde{{\cal O}}(\epsilon \langle \Re \xi
\rangle ^{-2})\\
&=&a_\xi \cdot {\partial \over \partial \Re x}({\partial \psi _\epsilon
\over \partial \Re \xi }\cdot \Re \xi )
+\widetilde{{\cal O}}(\epsilon \langle \Re \xi
\rangle ^{-2}),
\end{eqnarray*}
so
\ekv{qw.27}
{
\langle a_\xi
\vert (\psi _\epsilon )''_{x\xi}\xi \rangle =
\widetilde{{\cal O}}(\epsilon \langle \Re \xi \rangle ^{-2}).
}
Returning to \no{qw.24}, we get
\ekv{qw.28}
{
\xi \cdot {\partial \over \partial \xi }\widehat{\psi }_{\epsilon
,s}=\widetilde{{\cal O}}(\epsilon \langle \Re \xi \rangle ^{-1}) \hbox{ on
}\Lambda _s.
}
Also recall that $\widehat{\psi }_{\epsilon ,s}=\psi _\epsilon $ on
$\Lambda _s$, so
\ekv{qw.29}
{
\xi \cdot {\partial \over \partial \xi }\widehat{\psi }_{\epsilon
,s}-\widehat{\psi }_{\epsilon ,s}=\widetilde{{\cal O}}(\epsilon )\hbox{
on }\Lambda _s.
}
Combining this with \no{geo.28}, we get
\ekv{qw.30}
{
(\xi \cdot {\partial \over \partial \xi }\widehat{\psi }_{\epsilon
,s}-\widehat{\psi }_{\epsilon ,s} )\circ \kappa (s)=\widetilde{{\cal
O}}(\epsilon ),
}
and finally from \no{qw.23} we obtain
\begin{lemma}\label{qw0}
\ekv{qw.31}
{
h=\widetilde{{\cal O}}(\delta \epsilon ).
}
\end{lemma}

\par Following \cite{HeSj2} we shall now quantize $\kappa (\delta )$ by
means of a \fop{}
\ekv{qw.32}
{
Au(x;h)=\iint_{(y,\alpha )\in{\bf R}^n\times T^*M}e^{{i\over h}\phi
(x,y,\alpha )}a(x,y,\alpha ;h)\chi (x,y,\alpha _x)u(y)dyd\alpha ,
}
with $a\in S^{{3n\over 2}+m,{n\over 2}+k}$, and where $\chi $ is a standard cutoff to a \neigh{}
of the diagonal: $\vert x-\alpha _x\vert ,\, \vert y-\alpha _x\vert <1/C$, equal to 1 on a
smaller \neigh{} of the same type. Here we take with $\beta =\kappa
(\delta ,\alpha )=\kappa (\delta )(\alpha )$:
\ekv{qw.33}
{
\phi (x,y,\alpha )=(x-\beta _x)\cdot \beta _\xi +(\alpha _x-y)\cdot \alpha
_\xi +\psi (x-\beta _x,y-\alpha _x,\alpha )+h(\alpha ),
}
where on the real domain,
\ekv{qw.34}
{
\psi (x-\beta _x,y-\alpha _x,\alpha )={\cal O}(\langle \alpha _\xi
\rangle ((x-\beta _x)^2+(y-\alpha _x)^2),
}
in the symbol sense:
$$
\partial _x^k\partial _y^\ell \partial _{\alpha_x}^m \partial _{\alpha _\xi
}^p\psi (x,y,\alpha )={\cal O}(\langle \alpha _\xi \rangle ^{1-\vert
p\vert }(\vert x\vert +\vert y\vert )^{(2-\vert k+\ell \vert )_+}),
$$
and we take almost holomorphic extensions satisfying the same estimates
in the complex domain. Further, $a\in S^{{3n\over 2}+m,{n\over 2}+k}$ means that
$a={\cal O}(h^{-{3n\over 2}-m}\langle \alpha _\xi  \rangle ^{{n\over 2}+k})$ in
the symbol sense, $\partial _{x,y}^\ell \partial _{\alpha _x}^p
\partial _{\alpha_\xi }^q a(x,y,\alpha ;h)=
{\cal O}(h^{-{3n\over 2}-m}\langle \alpha _\xi  \rangle
^{{n\over 2}+k-\vert q \vert})$.

\par We also assume that
\ekv{qw.35}
{
\Im \psi \backsim \langle \alpha _\xi \rangle ((x-\beta _x)^2+(y-\alpha
_x)^2).}

\par Viewing $h$ as a \fu{} on the graph of $\kappa (\delta )$, we have
$dh=\beta _\xi \cdot d\beta _x-\alpha _\xi \cdot d\alpha _x$, so for
\ekv{qw.36}
{
x=\kappa (\delta ,\alpha )_x,\ y=\alpha _x,
}
we have $d_\alpha \phi =0$. Moreover, in a \neigh{} of that set, we have
with $\psi =\psi (x,y,\alpha )$,
\eekv{qw.37}
{
d_\alpha \phi &=&(x-\kappa (\alpha )_x)\cdot d_\alpha (\kappa (\alpha
)_\xi )+(\alpha _x-y)\cdot d\alpha _\xi
}
{&&
-{\partial \psi \over \partial x}(x-\kappa (\alpha )_x,y-\alpha
_x,\alpha )\cdot d_\alpha (\kappa (\alpha )_x)-{\partial \psi \over
\partial y}\cdot d\alpha _x+{\partial \psi \over \partial \alpha }\cdot
d\alpha .
}
So $\phi $ is a \nondeg{} phase \fu{} in the sense of H\"ormander (a part
from the homogeneity condition in the fiber variables) with a critical set
\no{qw.36}, the associated \ctf{} is $\kappa =\kappa (\delta )$.

\par Similarly, to $\kappa (\delta )\inv$ we can associate
\ekv{qw.38}
{
Bv(x;h)=\iint_{(y,\gamma )\in {\bf R}^n\times T^*M}e^{{i\over
h}\widehat{\phi }(x,y,\gamma )}b(x,y,\gamma ;h)\chi (x,y,\gamma
)v(y)dyd\gamma ,
}
$b\in S^{{3n\over 2}+\widehat{m},{n\over 2}+\widehat{k}},$
$\widehat{m}, \widehat{k}\in {\bf R}$,
where with $\beta =\kappa (\delta )(\gamma )$:
\ekv{qw.39}
{
\widehat{\phi }(x,y,\gamma )=(x-\gamma _x)\cdot \gamma _\xi +(\beta
_x-y)\cdot \beta _\xi +\widehat{\psi }(x-\gamma _x,y-\beta _x,\gamma )
-h(\gamma )
}
and $\widehat{\psi }$ satisfies \no{qw.34}. Again this is a \nondeg{} phase
and the critical set is given by
\ekv{qw.40}
{
x=\gamma _x,\ y=\beta _x, (\beta =\kappa (\delta )(\gamma )),
}
and the associated \ctf{} is $\kappa (\delta )\inv$.

\par For $x,y\in {\bf R}^n$, $\vert x-\alpha _x\vert ,\, \vert y-\beta
_x\vert \le 1/C$, we have with $\beta =\kappa (\alpha )=\kappa (\delta
,\alpha )$:
$$
\Im \phi (x,y,\alpha )=(x-\Re \beta _x)\cdot \Im \beta _\xi -\Im \beta
_x\cdot \Re \beta _\xi +\Im \psi (x-\beta _x,y-\alpha _x,\beta ,\alpha
)+\Im h(\alpha ).
$$
Here
$$
\Im \beta _x=\delta \widetilde{{\cal O}}(\epsilon ^{1\over 2}\langle
\Re \alpha _\xi \rangle \inv ),\ \Im \beta _\xi =\delta
\widetilde{{\cal O}}(\epsilon ^{1\over 2}),
$$
$$
\Im \beta _x\cdot \Re \beta _\xi =\delta \Re \beta _\xi \cdot {\partial
\psi _\epsilon \over \partial \Re \beta _\xi }=\widetilde{{\cal O}}
({\delta \epsilon
\over \langle \Re \beta _\xi \rangle ^2 }),
$$
and
$$
\Im h(\alpha )=\widetilde{{\cal O}}(\delta \epsilon ),
$$
$$
\Im \psi (x-\beta _x,y-\alpha _x,\alpha
)\ge \langle \alpha _\xi \rangle ({\vert x-\Re \beta _x\vert ^2\over
C}+{1\over C}\vert y-\alpha _x\vert ^2-C\vert \Im \beta _x\vert ^2),
$$
so
$$
\Im \phi (x,y,\alpha )\ge {\langle \alpha _\xi \rangle \over C}(\vert
x-\Re \beta _x\vert ^2+\vert y-\alpha _x\vert ^2)-C{\delta ^2\epsilon
\langle \alpha _\xi \rangle \over \langle \alpha _\xi \rangle ^2}-C\vert
x-\Re \beta _x\vert \delta \epsilon ^{1\over 2}-C\delta \epsilon ,
$$
\ekv{qw.41}
{
\Im \phi \ge {\langle \alpha _\xi \rangle \over 2C}(\vert x-\Re\beta _x\vert
^2+\vert y-\alpha _x\vert ^2)-\widetilde{C}\delta \epsilon .
}

\par Recalling that $\epsilon =Ah$, $A\gg 1$, we get
\ekv{qw.42}
{
\vert e^{{i\over h}\phi (x,y,\alpha )}\vert \le \exp (-{\langle \alpha _\xi
\rangle \over 2Ch}(\vert x-\Re \beta _x\vert ^2+\vert y-\alpha _x\vert
^2)+\widetilde{C}\delta A).
}
A similar estimate holds for $\widehat{\phi }$. It follows that $A$, $B$
are well-defined \op{}s: ${\cal S}\to {\cal S}$ with semi-norm estimates
that are uniform in powers of $h$. Also for every $s\in {\bf R}$,
there is an $s'\in {\bf R}$ such that $A,B:\,H^s\to H^{s'}$ with norms bounded
by some power of $h$.  Moreover, our \op{}s
are \indep{} of the choice of almost \hol{} extensions of the phase and
amplitude modulo \op{}s whose integral kernels are ${\cal O}_A(h^\infty )$
with all their derivatives and supported in a domain of the form $\vert
x-y\vert \le {\cal O}(1)$.

\begin{prop}\label{qw1}
If $m,k,\widehat{m},\widehat{k}=0$, then $A,B={\cal O}(\exp {\cal
O}(1)\delta A ):L^2({\bf R}^n)\to L^2({\bf R}^n)$
\end{prop}
The proof will be given later.
\begin{prop}\label{qw2}
We have $BA={\rm Op}_h(c)$, where $c=\widetilde{{\cal
O}}(h^{-m-\widehat{m}}\langle \alpha _\xi \rangle ^{k+\widehat{k}})$
\end{prop}
\begin{proof}

We have
$$
BAu(x)=\iiiint e^{{i\over h}(\widehat{\phi }(x,z,\gamma )+\phi
(z,y,\alpha ))}b(x,z,\gamma ;h)a(z,y,\alpha ;h)u(y)dyd\alpha dzd\gamma ,
$$
where the cutoffs $\chi $, $\widehat{\chi }$ have been incorporated in
$a$, $b$. Here
\eeekv{qw.47}
{
&&\hskip -1truecm \widehat{\phi }(x,z,\gamma )+\phi (z,y,\alpha )=}{&&(x-\gamma _x)\cdot \gamma
_\xi +(\kappa (\gamma )_x-z)\cdot \kappa (\gamma )_\xi +(z-\kappa (\alpha
)_x)\cdot \kappa (\alpha )_\xi }{&& +(\alpha _x-y)\cdot \alpha _\xi +
\widehat{\psi }(x-\gamma _x,z-\kappa (\gamma )_x,\gamma  )+
\psi (z-\kappa (\alpha )_x,y-\alpha _x,\alpha )+h(\alpha
)-h(\gamma ).}

\par The contribution to the distribution kernel from a region with
$\vert \gamma _\xi -\alpha _\xi \vert /\langle \alpha _\xi \rangle \ge
C^{-1}$ is ${\cal O}(h^\infty )$, as can be seen by integration by parts
\wrt{} $z$. More precisely the distribution kernel of $BA$ is
$$
(BA)(x,y)=\iiint e^{{i\over h}(\widehat{\phi }+\phi )}ba\chi ({\vert \gamma
_\xi -\alpha _\xi \vert \over \langle \alpha _\xi \rangle })d\alpha
dzd\gamma +R(x,y;h),
$$
where $\chi \in C_0^\infty ([0,\infty [)$ is equal to 1 near 0, and
$\partial_x^k\partial _y ^\ell R={\cal O}(h^\infty ).$ Here \no{qw.42} and the
analogous
estimate for $e^{{i\over h}\widehat{\phi }}$ are essential of course.

\par Next using \no{qw.41} and the similar estimate for $\Im
\widehat{\phi }$, we get a localization in $\vert \gamma_x -\alpha _x\vert $,
leading to
\ekv{qw.48}
{
(BA)(x,y)=\iiint e^{{i\over h}(\widehat{\phi }+\phi )}ba\chi ((\gamma
_x-\alpha _x)^2+{( \gamma
_\xi -\alpha _\xi )^2 \over \langle \alpha _\xi \rangle ^2 })d\alpha
dzd\gamma +\widetilde{R}(x,y;h),
}
where $\partial _x^k\partial _y^\ell\widetilde{R}={\cal O}(h^\infty )$.
Here $\chi \in C_0^\infty ([0,\epsilon _0[)$, $\chi =1$
near 0, and $\epsilon _0$ can be any fixed number.

\par In the integral \no{qw.48}, we may assume that $\vert x-\gamma
_x\vert ,\, \vert y-\alpha _x\vert ,\,\vert z-\kappa (\alpha )_x\vert
<1/C$, where $C$ is as large as we like, since the integral in the
complementary region is exponentially small. We now want to eliminate
integration variables by means of the method of stationary phase, and we
start by carrying out the $z$-integration, so we first look for the
critical point of
\ekv{qw.49}
{
z\mapsto \widehat{\phi }(x,z,\gamma )+\phi (z,y,\alpha ),
}
where $\phi ,\widehat{\phi }$ also denote almost \hol{} extensions. Let $F_\delta (x,\gamma ,y,\alpha )$ denote the corresponding critical value.

\par In order to understand this function, we first treat $\widetilde{\gamma }
=\kappa (\gamma )$, $\widetilde{\alpha }
=\kappa (\alpha )$ as independent variables (writing $\kappa (\alpha )$ instead of $\kappa (\delta ,\alpha )$ for short). Let
\begin{eqnarray*}
G(x,\gamma ,\widetilde{\gamma },y,\alpha ,\widetilde{\alpha })
&=&{\rm vc}_z[(x-\gamma _x)\cdot \gamma _\xi +(\widetilde{\gamma }_x-z)\cdot \widetilde{\gamma }_\xi +(z-\widetilde{\alpha }_x)\cdot \widetilde{\alpha }_\xi
+(\alpha _x-y)\cdot \alpha _\xi \\
&&\hskip 15 truemm +\widehat{\psi }(x-\gamma _x,z-
\widetilde{\gamma }_x,\gamma )+\psi (z-
\widetilde{\alpha}_x, y-\alpha _x,\alpha )].
\end{eqnarray*}
Here the critical point $z=z_c(x,\gamma ,\widetilde{\gamma },y,\alpha ,\widetilde{\alpha })$ satisfies
$$
z_c=\widetilde{\alpha }_x+{\cal O}(\widetilde{\alpha }_x-\widetilde{\gamma }_x,
{\widetilde{\alpha }_\xi-\widetilde{\gamma }_\xi\over \langle
\widetilde{\gamma }_\xi  \rangle},(x-\gamma _x)^2,(y-\alpha _x)^2)
$$
in the natural symbol sense. Notice that
$$
G(\gamma _x,\gamma ,\widetilde{\gamma };\alpha _x,\alpha ,\widetilde{\gamma })=0.
$$
Moreover,
\begin{eqnarray*}
dG&=&\gamma _\xi \cdot dx-\alpha _\xi \cdot dy+\widetilde{\gamma }_\xi \cdot
d\widetilde{\gamma }_x-\gamma _\xi \cdot d\gamma _x-(\widetilde{\alpha }_\xi
\cdot d\widetilde{\alpha }_x-\alpha _\xi \cdot d\alpha _x)\\
&& \hskip 4 truecm +{\cal O}(\langle \gamma _\xi  \rangle (x-\gamma _x,
y-\alpha _x, \widetilde{\gamma }_x-\widetilde{\alpha }_x,
{\widetilde{\gamma }_\xi -
\widetilde{\alpha }_\xi \over\langle \gamma _\xi  \rangle }))
\end{eqnarray*}
in the same symbol sense, with the convention that the remainder term is
expressed as linear combination of the ``normalized'' forms $dx$, $dy$,
$d\alpha _x$, $\widetilde{\alpha }_x$, $d\gamma _x$, $\widetilde{\gamma }_x$,
$\langle \gamma _\xi  \rangle\inv d\alpha _\xi $, ... .

\par Now,
$$
F_\delta ={{G}_\big{\vert }}_{\widetilde{\gamma }=\kappa (\gamma ),\atop
\widetilde{\alpha}=\kappa (\alpha ) }+h(\alpha )-h(\gamma ),
$$
and we get
\ekv{qw.49.5}
{
dF_\delta =\gamma _\xi \cdot dx-\alpha _\xi \cdot dy+
{{{\cal O}(\langle \gamma _\xi  \rangle(x-\gamma _x,y-\alpha _x,
\widetilde{\gamma }_x-\widetilde{\alpha }_x,{\widetilde{\gamma }_\xi
-\widetilde{\alpha }_\xi
\over \langle \gamma _\xi  \rangle}))}_\big{\vert}}_{\widetilde{\gamma }=
\kappa
(\gamma )\atop \widetilde{\alpha }=\kappa  (\alpha )  }.
}
Moreover, $F_\delta (\gamma _x,\gamma ;\gamma _x,\gamma )=0$.
From \no{geo.28}, we know that $\widetilde{\alpha }_x=\alpha _x+
\widetilde{{\cal O}}(\delta \epsilon ^{1/2}\langle \alpha _\xi  \rangle
\inv)$, $\widetilde{\alpha }_\xi =\alpha _\xi +\widetilde{{\cal O}}
(\delta \epsilon ^{1/2})$ when $\widetilde{\alpha }=\kappa (\alpha )$, so
$$\widetilde{\gamma }_x-\widetilde{\alpha }_x=\gamma _x-\alpha _x+
\widetilde{{\cal O}}({\delta\over \langle \gamma _\xi  \rangle })
(\gamma _x-\alpha _x,
{\gamma _\xi-\alpha _\xi  \over \langle \gamma _\xi  \rangle})$$ and
similarly for ${\widetilde{\gamma }_\xi -\widetilde{\alpha }_\xi \over
\langle \alpha _\xi  \rangle}$.

\par Hence,
\ekv{qw.49.6}
{
dF_\delta =\gamma _\xi \cdot dx-\alpha _\xi \cdot dy
+{\cal O}(\langle \gamma _\xi  \rangle)
(x-\gamma _x,y-\alpha _x,\gamma _x-\alpha _x,{\gamma _\xi -\alpha _\xi
\over \langle \gamma _\xi  \rangle})+
\widetilde{{\cal O}}(\delta )(\gamma _x -\alpha _x,{\gamma _\xi -\alpha _\xi
 \over \langle \gamma _\xi  \rangle}),
}
and integrating this, we get
\begin{eqnarray}
\label{qw.51}
F_\delta (x,\gamma ;y,\alpha )
&=& \gamma _\xi \cdot (x-\gamma _x)-\alpha _\xi \cdot (y-\alpha _x)\\
&&+({\cal O}(\langle \gamma _\xi  \rangle )+\widetilde{{\cal O}}(\delta ))
(x-\gamma _x,y-\alpha _x,\gamma _x-\alpha _x,{\gamma _\xi -\alpha _\xi \over
\langle \gamma _\xi  \rangle})^2,
\nonumber
\end{eqnarray}
where the loss of $\epsilon ^{1/2}$ for each differentiation appears in the
variables $\alpha ,\gamma $ only. When $\delta =0$, we have on the real domain
\ekv{qw.52}
{
\Im F_\delta \backsim \langle \alpha _\xi  \rangle
((x-\gamma _x)^2+(y-\alpha _x)^2+(\gamma _x-\alpha _x)^2+
\Big( {\gamma _\xi -\alpha _\xi \over\langle \gamma _\xi  \rangle}\Big)^2),
}
and in view of \no{qw.51} this persists for $0\le \delta \ll 1$.

\par When applying stationary phase to \no{qw.48} we also have to make a
deformation of the integration contour in order to pass through the
critical point $z_c$. Here we recall  from \cite{MeSj} and \no{qw.41} and
its analogue for $\widehat{\phi }$ that
$$
{1\over \langle \alpha _\xi \rangle }\Im F_\delta \ge {1\over C}(\Im
z_c)^2-{C\delta \epsilon \over \langle \alpha _\xi \rangle },
$$
so the error from $\overline{\partial } _z$ of the almost \hol{} extension,
appearing in Stokes' formula, is
${\cal O}((\delta \epsilon +\Im F_\delta )/\langle \alpha _\xi \rangle
)^\infty )$. Since $\epsilon ={\cal O}(h)$, we conclude that
\ekv{qw.53}
{
(BA)(x,y)=\iint e^{{i\over h}F_\delta (x,\gamma ,y,\alpha )}d(x,\gamma
,y,\alpha ;h)\chi ((\gamma _x-\alpha _x)^2+({\gamma _\xi -\alpha _\xi
\over \langle \alpha _\xi \rangle })^2)d\alpha d\gamma +\widehat{R},
}
where $\widehat{R}$ has the same properties as $\widetilde{R}$ in \no{qw.48}
and
\ekv{qw.54}
{
\partial _{x,y}^k\partial _{\gamma _x,\alpha _x}^\ell \partial _{\gamma
_\xi ,\alpha _\xi }^p d={\cal O}(h^{-m-\widehat{m}-3n+{n\over 2}}\langle
\alpha _\xi \rangle ^{k+\widehat{k}+{n\over 2}-\vert p\vert }\epsilon
^{-{\vert \ell +p\vert \over 2}}).}

\par We now compute the Weyl symbol $c$ of $BA$ by means of the formula
\ekv{qw.55}
{
c(x,\xi ;h)=\int (BA)(x+{w\over 2},x-{w\over 2})e^{-iw\xi /h}dw.
}
The contribution from $\widehat{R}$ in \no{qw.53} is ${\cal O}(h^\infty
\langle \xi \rangle ^{-\infty })$ with all its derivatives. The
contribution from the integral in \no{qw.53} is
$$
\iiint e^{{i\over h}(F_\delta (x+{w\over 2},\gamma ,x-{w\over 2},\alpha
)-w\cdot \xi )}d(x+{w\over 2},\gamma ,x-{w\over 2},\alpha ;h)\chi ((\gamma
_x-\alpha _x)^2+({\gamma _\xi -\alpha _\xi \over \langle \alpha _\xi
\rangle })^2)dwd\alpha d\gamma .
$$
The contribution from a region $\{ \vert \xi -\alpha _\xi \vert \ge {1\over
C}\langle \alpha _\xi \rangle \}$ is ${\cal O}(h^\infty \langle \xi
\rangle ^{-\infty }$
with all its derivatives and the remaining region can be treated with the
method of stationary phase by working in the dilated tilde
variables given by $\gamma _x=x_0+\sqrt{\epsilon }\widetilde{\gamma
}_x,$ ... as in the addendum below. The proposition follows.
\end{proof}

\par\noindent\bf Addendum: \it Stationary phase with
$\widetilde{{\cal O}}$-symbols. \rm Assume
\begin{eqnarray*}
&&\phi \in C^\infty ({\bf R}^n),\ \Im \phi \ge 0,\ \phi (0)=0, \\
&&\vert \phi '(x)\vert \backsim \vert x\vert ,\ \phi '=x\cdot
\widetilde{{\cal O}}(1),\
\det \phi ''(0)\ne 0\hbox{ \ufly{} in }\epsilon .
\end{eqnarray*}
Let $a=\widetilde{{\cal O}}(1)$. We shall establish a stationary phase
development for
$$
I(h)=h^{-{n\over 2}}\int e^{{i\over h}\phi (x)}a(x)dx
$$
in powers of $\widetilde{h}=h/\epsilon $. (Assume $h\ll \epsilon \ll 1$.)

\par Put $x=\epsilon ^{1\over 2}\widetilde{x}$,
\ekv{qw.55.5}
{
I(h)=({h\over \epsilon })^{-{n\over 2}}\int e^{{i\over
\widetilde{h}}\widetilde{\phi
}(\widetilde{x})}\widetilde{a}(\widetilde{x})d\widetilde{x},
}
where $\phi (\epsilon ^{1\over 2}\widetilde{x})=\epsilon
\widetilde{\phi }(\widetilde{x})$, $\widetilde{\phi
}(\widetilde{x})={1\over \epsilon }\phi (\epsilon ^{1\over
2}\widetilde{x})$, $\widetilde{a}(\widetilde{x})=a(\epsilon ^{1\over
2}\widetilde{x})$.

\par Then $\widetilde{a}={\cal O}(1)$ in the symbol sense; $\partial
^\alpha \widetilde{a}={\cal O}(1)$, and
$$
\vert \widetilde{\phi }'(\widetilde{x})\vert =\vert {1\over \sqrt{\epsilon
}}\phi '(\sqrt{\epsilon }\widetilde{x})\vert \backsim \vert
\widetilde{x}\vert ,\ \widetilde{\phi
}'(\widetilde{x})=\widetilde{x}\cdot {\cal O}(1)
$$
and $\widetilde{\phi }''(0)=\phi ''(0)$.
The contribution to $I(h)$ from $\vert
\widetilde{x}\vert \ge 1$ is ${\cal O}(\widetilde{h}^\infty )$ by repeated
integrations by parts, and the contribution from $\vert \widetilde{x}\vert
<1$ can be handled in the usual way since $\widetilde{\phi }={\cal O}(1)$ here. Thus
\ekv{qw.55.7}
{
I(h)\sim \sum_0^\infty  c_j\widetilde{h}^j,\quad c_0={(2\pi )^{n/2}\over
\sqrt{\det \phi ''(0)}}e^{i{\pi \over 4}{\rm sgn\,}\phi ''(0)}.
}
\hfill{$\Box$}
\medskip
\par It follows from the proof that the proposition remains valid if we
relax the symbol condition in $y$ and only assume
\ekv{qw.56}
{
\partial _x^\ell \partial _{\alpha _x,y}^q\partial _{\alpha _\xi }^pa={\cal
O}(h^{-{3n\over 2}-m}\langle \alpha _\xi \rangle ^{{n\over 2}+k-\vert
p\vert }\epsilon ^{-{1\over 2}(\vert q\vert +\vert p\vert )}), } i.e. we
also allow for a loss of $\epsilon ^{1/2}$ for each $y$-derivation.
Similarly for $B$ (cf.  \no{qw.38}) we can content ourselves with
\ekv{qw.57}
{
\partial _y^\ell \partial _{\gamma _x,x}^q\partial _{\gamma _\xi }^pb=
{\cal
O}(h^{-{3n\over 2}-\widehat{m}}\langle \alpha _\xi \rangle ^{{n\over
2}+\widehat{k}-\vert
p\vert }\epsilon ^{-{1\over 2}(\vert q\vert +\vert p\vert )}).
}
Moreover, if $b$ an $a$ are elliptic, then $c$ is elliptic.

\par We get by standard arguments,
\begin{prop}\label{qw3}
Let $A$ be an elliptic \fop{} of order $(m,k)$ with symbol as in
{\rm \no{qw.56}}. Then there exists an elliptic \fop{} $B$ of order $(-m,-k)$
with symbol as in {\rm \no{qw.57}}, such that
\ekv{qw.58}
{
BA=1+R,
}
where $R$ is $1$-negligible in the sense that its symbol $R$ is
$\widetilde{{\cal O}}(({h\over \epsilon })^\infty \langle \xi \rangle
^{-\infty })$. In particular $A$ has the left inverse $(1+R)\inv B$ when
$\epsilon /h\gg 1$.
\end{prop}

\par We notice that when $\delta =0$, then $A,B$ are elliptic \pop{}s and
hence $(1+R)\inv B$ is also a  left inverse. (By the Beals lemma we also
know that $(1+R)\inv$ is an $h$-\pop{} with symbol $1+\widetilde{{\cal
O}}(({h\over \epsilon \langle \xi \rangle })^\infty )$.) For general
small $\delta $, ${\cal R}(A_\delta )$ is closed. Using suitable
deformations of $A$ we can produce a continuous deformation of closed
subspaces in $L^2$ from $L^2$ to ${\cal R}(A_\delta )$. All the deformed
subspaces then have to be equal to $L^2$ and ${\cal R}(A_\delta )=L^2$,
so $(1+R_\delta )\inv B_\delta $ is also a left inverse of $A_\delta $.
\medskip

\par We next turn to Egorov's theorem and start with some preparations.
Recall that by \no{geo.28}
\ekv{qw.59}
{
\kappa (\alpha )_x=\alpha _x+\widetilde{{\cal O}}(\delta \epsilon ^{1\over
2}\langle \alpha _\xi \rangle \inv ),\ \kappa (\alpha )_\xi =\alpha _\xi
+\widetilde{{\cal O}}(\delta \epsilon ^{1\over 2}).
}
It follows that
\eekv{qw.60}
{
d_\alpha (\kappa (\alpha )_\xi )&=&d\alpha _\xi +\widetilde{{\cal
O}}(\delta )d\alpha _x+\widetilde{{\cal O}}({\delta
\over \langle \alpha _\xi \rangle })d\alpha _\xi
}
{
d_\alpha (\kappa (\alpha )_x )&=&d\alpha _x +\widetilde{{\cal
O}}({\delta \over \langle \alpha _\xi \rangle })
d\alpha _x+\widetilde{{\cal O}}({\delta
\over \langle \alpha _\xi \rangle ^2 })d\alpha _\xi .
}
Substituting this into \no{qw.37}, we get
\eeeekv{qw.61}
{
d_\alpha \phi &=&((x-\kappa (\alpha )_x)-(y-\alpha _x))\cdot d\alpha
_\xi}{&&- ({\partial \psi \over \partial x}+{\partial \psi \over \partial
y})(x-\kappa (\alpha )_x,y-\alpha _x,\alpha _\xi )\cdot d\alpha _x
}
{&&+(x-\kappa (\alpha )_x,y-\alpha _x)\cdot (\widetilde{{\cal O}}(\delta
)d\alpha _x+\widetilde{{\cal O}}({\delta \over \langle \alpha _\xi
\rangle })d\alpha _\xi )}{&&+{\partial
\psi \over \partial \alpha }(x-\kappa (\alpha )_x,y-\alpha _x,\alpha
)\cdot d\alpha ,
}
which we can write
\ekv{qw.62}
{
\pmatrix{{\partial _{\alpha _x}\phi / \langle \alpha _\xi \rangle } \cr
\partial _{\alpha _\xi} \phi }=\pmatrix{ ({{\partial \psi \over \partial
\alpha _x}-({\partial \psi \over \partial
x}+{\partial \psi \over \partial y}))/ \langle \alpha _\xi \rangle}
\cr (x-\kappa (\alpha )_x)-(y-\alpha _x)+{\partial \psi \over \partial
\alpha _\xi }}+\widetilde{{\cal O}}({\delta \over
\langle \alpha _\xi \rangle })\pmatrix{x-\kappa (\alpha )_x\cr y-\alpha _x}.
}
Here we write $a(x,y,\alpha )=\widetilde{{\cal O}}(m)$ if $\partial
_{x,y}^\nu \partial _{\alpha _x}^\mu \partial _{\alpha _\xi }^\rho
a={\cal O}(m\epsilon ^{-{1\over 2}(\vert \mu \vert +\vert \rho \vert
)}\langle \alpha _\xi \rangle ^{-\vert \rho \vert })$.
The differential of this vector at a point of the critical set is given by
the matrix
\ekv{qw.63}
{
\pmatrix{ -\langle \alpha _\xi \rangle \inv (\psi ''_{xx}+\psi ''_{yx})&
-\langle \alpha _\xi \rangle \inv (\psi ''_{xy}+\psi ''_{yy})\cr 1
&-1}
+\widetilde{{\cal O}}({\delta \over \langle \alpha
_\xi \rangle }).
}
If $\pmatrix{t_x\cr t_y}$ is in the kernel of the first term, we get $t_x=t_y$ and
$(\psi ''_{xx}+\psi ''_{yx}+\psi ''_{xy}+\psi ''_{yy})t_x=0$. Here we
recognize the Hessian of $z\mapsto \psi (z,z,\alpha )$ which is invertible
because of the assumption on the imaginary part. More precisely, the matrix
\no{qw.63} has a \ufly{} \bdd{} inverse. From \no{qw.62}, we get
\ekv{qw.64}
{
\pmatrix{x-\kappa (\alpha )_x\cr y-\alpha _x}=M_\delta (x,y,\alpha
)\pmatrix{\partial _{\alpha _x}\phi /\langle \alpha _\xi \rangle \cr
\partial _{\alpha _\xi }\phi },\ \vert x-\kappa (\alpha )_x\vert ,\vert
y-\alpha _x\vert \le 1/{\cal O}(1),
}
where
\ekv{qw.65}
{
M_\delta (x,y,\alpha )=M_0(x,y,\alpha )+\widetilde{{\cal O}}({\delta
\over \langle \alpha _\xi \rangle }),
}
and $\partial _{x,y,\alpha _x}^k\partial _{\alpha _\xi }^\ell M_0={\cal
O}(\langle \alpha _\xi \rangle ^{-\vert \ell \vert })$.

\begin{lemma}\label{qw4}
Let $\Phi (x,y,\alpha )=\widetilde{{\cal O}}(\langle \alpha _\xi
\rangle ^{\widetilde{m}})$ in the sense that
$$
\partial _x^k\partial _{y,\alpha _x}^\ell \partial _{\alpha _\xi }^p\Phi
={\cal O} (\epsilon ^{-{1\over 2}(\vert \ell\vert +\vert p\vert )}\langle
\alpha _\xi \rangle ^{\widetilde{m}-\vert p\vert }).
$$
Similarly, let $a=\widetilde{{\cal O}}(\langle \alpha _\xi \rangle
^{\widetilde{n}})$ be elliptic. Then
\eeekv{qw.66}
{
&&\hskip -1truecm\iint e^{{i\over h}\phi }\Phi (x,y,\alpha )a(x,y,\alpha ;h)\chi
(x,y,\alpha )u(y)dyd\alpha
}
{
&=&\iint e^{{i\over h}\phi }\Phi (\kappa (\alpha )_x,\alpha _x,\alpha )a(x,y,\alpha ;h)\chi
(x,y,\alpha )u(y)dyd\alpha
}
{&&+
\iint e^{{i\over h}\phi }\widetilde{{\cal O}}({h\over \epsilon }\langle
\alpha _\xi \rangle ^{\widetilde{m}-1})a(x,y,\alpha ;h)\widetilde{\chi }
(x,y,\alpha )u(y)dyd\alpha,
}
where $\widetilde{\chi }$
is a similar cutoff.
\end{lemma}
\begin{proof}
We have
\begin{eqnarray*}
\Phi (x,y,\alpha )&=&\Phi (\kappa (\alpha )_x,\alpha _x,\alpha
)+\widetilde{{\cal O}}(\epsilon ^{-{1\over 2}}\langle \alpha _\xi
\rangle ^{\widetilde{m}})\pmatrix{x-\kappa (\alpha )_x\cr y-\alpha _x}\\
&=&\Phi (\kappa (\alpha )_x,\alpha _x,\alpha
)+\widetilde{{\cal O}}(\epsilon ^{-{1\over 2}}\langle \alpha _\xi
\rangle ^{\widetilde{m}-1})\cdot \partial _{\alpha _x}\phi
+\widetilde{{\cal O}}(\epsilon ^{-{1\over 2}}\langle \alpha _\xi
\rangle ^{\widetilde{m}})\cdot \partial _{\alpha _\xi }\phi .
\end{eqnarray*}
The contribution from the remainders to the \lhs{} in \no{qw.66} is
therefore
\begin{eqnarray*}
&&h\iint {\partial \over \partial \alpha _x}(e^{{i\over h}\phi })\cdot
\widetilde{{\cal O}}(\epsilon ^{-{1\over 2}}\langle \alpha _\xi \rangle
^{\widetilde{m}-1})a\chi udyd\alpha \\
&&+h\iint {\partial \over \partial \alpha _\xi }(e^{{i\over h}\phi })\cdot
\widetilde{{\cal O}}(\epsilon ^{-{1\over 2}}\langle \alpha _\xi \rangle
^{\widetilde{m}})a\chi udyd\alpha,
\end{eqnarray*}
and it suffices to integrate by parts.
\end{proof}

\par Actually, we shall not use the lemma directly, only its proof. The next
result is closely related and could probably be obtained from Lemma
\ref{qw4}. We will give a different proof however.
\begin{lemma}\label{qw5}
Under the same assumptions as in the preceding lemma, we have
\ekv{qw.67}
{
\iint e^{{i\over h}\phi }\Phi a\chi u(y)dyd\alpha =AQ,
}
where $Q$ is an $h$-\pop{} with symbol
\ekv{qw.68}
{
Q=\Phi (\kappa (x,\xi )_x,x,x,\xi )+\widetilde{{\cal O}}({h\over \epsilon
}\langle \xi \rangle ^{\widetilde{m}-1}).
}
Here
\ekv{qw.69}
{
Au(x)=\iint e^{{i\over h}\phi }a\chi u(y)dyd\alpha .
}
\end{lemma}

\begin{proof}
In order to harmonize with Proposition \ref{qw2}, we may change the
assumption on $a$ to $a=\widetilde{{\cal O}}(\langle \alpha _\xi
\rangle ^{\widetilde{m}+{n\over 2}}h^{-{3n\over 2}})$ and assume that $a$
is elliptic in this class. Here $\widetilde{{\cal O}}$ indicates a loss
of $\epsilon ^{1/2}$ for each differentiation in $\alpha ,y$. Let $B$ be of
the form \no{qw.38} with $b=\widetilde{{\cal O}}(\langle \alpha _\xi
\rangle ^{-\widetilde{m}+{n\over 2}}h^{-{3n\over 2}})$ where the
$\epsilon ^{1/2}$ loss is now for each differentiation in $(x,\alpha )$.
We also assume that $b$ is elliptic. Then from Proposition \ref{qw2} and the
remark following its proof, we know that $BA={\rm Op}_h(c)$, where
$c=\widetilde{{\cal O}}(1)$ is elliptic. Recalling that the proof was by
stationary phase, we see that if $\widetilde{A}$ is the \op{} given by
the \lhs{} of \no{qw.67}, then $B\widetilde{A}={\rm Op}_h(\widetilde{c})$,
where
$$
\widetilde{c}=\Phi (\kappa (x,\xi ),x,x,\xi )c+\widetilde{{\cal
O}}({h\over \epsilon }\langle \xi \rangle ^{\widetilde{m}-1}).
$$
By pseudodifferential calculus, we get $B\widetilde{A}=BAQ$, with $Q$ as in
the lemma. Moreover, $B$ is invertible.
\end{proof}
\begin{remark}\label{qw5'}\rm
Later on, we shall meet the special situation when $\Phi =\Phi (\alpha
)=\widetilde{{\cal O}}(\langle \alpha _\xi \rangle ^{m_1})P(\alpha
)=:RP$, $P(\alpha )={\cal O}(\langle \xi \rangle ^{m_2})$. In this case we
have
\ekv{qw.69.5}
{
Q=\widetilde{{\cal O}}(\langle \xi \rangle ^{m_1})\# P+\widetilde{{\cal
O}}({h\over \sqrt{\epsilon }}\langle \xi \rangle ^{m_1+m_2-1}).
}
In fact (anticipating on a part of the proof of Proposition \ref{qw7}), let
$$
Bu=\iint e^{{i\over h}\phi }Ra\chi u(y)dyd\alpha .
$$
Then (with ${\rm Op}={\rm Op}_h$ when nothing else is specified)
\begin{eqnarray*}
B\circ {\rm Op}(P)u&=&\iint {\trans{P}}(y,hD_y)(e^{{i\over h}\phi }Ra\chi
)u(y)dyd\alpha \\
&=&
\iint (P(y,-\phi '_y)Ra\chi +\widetilde{{\cal O}}(h^2\langle \alpha _\xi
\rangle ^{m_1+m_2-2}))e^{{i\over h}\phi }u(y)dyd\alpha .
\end{eqnarray*}
As we shall see in the proof of Proposition \ref{qw7},
\begin{eqnarray*}
&&\hskip -1cm \iint P(y,-\phi '_y)Ra\chi e^{{i\over h}\phi }u(y)dyd\alpha =\\
&&\iint e^{{i\over h}\phi }\Phi a\chi u(y)dyd\alpha +\iint e^{{i\over
h}\phi }\widetilde{{\cal O}}({h\over \sqrt{\epsilon }}\langle \alpha _\xi
\rangle ^{m_1+m_2-1})u(y)dyd\alpha .
\end{eqnarray*}
Applying Lemma \ref{qw5} to $B$ and the various remainder terms, we get
\no{qw.69.5}.
\end{remark}

\begin{prop}\label{qw6}
Let $A_s$ be a \fop{} quantizing $\kappa (s)$ as in Proposition {\rm {\ref{qw2}}}.
Also, assume that $a=a_s$ is elliptic and depends smoothly on $s$ in the
sense that
$$
\partial _s^k a_s=\widetilde{{\cal O}}(h^{-m-{3n\over 2}}\langle \alpha
_\xi \rangle ^{k+{n\over 2}}),
$$
where $\widetilde{{\cal O}}$ indicates a loss of $\epsilon ^{1/2}$ for
differentiations in $\alpha ,y$. Then
\ekv{qw.70}
{
hD_sA_s+iA_s\check{\psi }_s=0, \ 0\le s\le \delta _0,\ 0<\delta _0\ll 1,
}
where $ D_s=i\inv \partial _s$ and
$\check{\psi }_s$ is an $h$-\pop{} with symbol
\ekv{qw.71}
{
\check{\psi }_s=\widetilde{{\cal O}}(h+{h\langle \xi \rangle\inv
})+\widehat{\psi }_{\epsilon ,s}\circ \kappa (s).
}
The term "$h$" in the remainder can be dropped if $a_s$ is
\indep{} of $s$.
\end{prop}

\begin{proof}
Assume for simplicity that $a_s$ is \indep{} of $s$. Then,
\ekv{qw.72}
{
hD_sA_su(x)=\iint e^{{i\over h}\phi }a\partial _s\phi \, udyd\alpha .
}
Recalling that,
$$
\phi _s(x,y,\alpha )=(x-\kappa (s,\alpha )_x)\cdot \kappa (s,\alpha
)_\xi +(\alpha _x-y)\cdot \alpha _\xi +\psi (x-\kappa (s,\alpha
)_x,y-\alpha _x,\alpha )+h_s(\alpha ),
$$
we get
\eekv{qw.73}
{ &&\hskip -2cm \partial _s\phi _s(x,y,\alpha )= (x-\kappa (s,\alpha
)_x)\cdot \partial _s\kappa (s,\alpha )_\xi -\partial _s\kappa (s,\alpha
)_x\cdot \kappa (s,\alpha )_\xi} {&&\hskip 2cm -(\partial _x\psi )(x-\kappa
(s,\alpha )_x,y-\alpha _x,\alpha )\cdot \partial _s\kappa (s,\alpha
)_x+\partial _sh_s(\alpha ).  } The restriction to the critical set
\no{qw.36} is
\eeeekv{qw.74}
{
(\partial _s\phi _s)(\kappa (s,\alpha )_x,\alpha _x,\alpha )&=&(\partial
_sh_s)(\alpha )-\partial _s\kappa (s,\alpha )_x\cdot \kappa (s,\alpha
)_\xi
}
{
&=& (\partial _sh_s)(\alpha )-i({\partial \over \partial \xi
}\widehat{\psi }_{\epsilon ,s}\cdot \xi )\circ \kappa (s)
}
{
&=& -i\widehat{\psi }_{\epsilon ,s}\circ \kappa (s)(\alpha )
}
{
&=&\widetilde{{\cal O}}(\epsilon ),
}
where we used \no{qw.23}. More generally, from \no{qw.73} we get
\ekv{qw.75}
{
\partial _s\phi _s=-i\widehat{\psi }_{\epsilon ,s}\circ \kappa (s)(\alpha
)+\widetilde{{\cal O}}(\epsilon ^{1\over 2})\cdot (x-\kappa (s,\alpha
)_x,y-\alpha _x).
}
As in the proof of Lemma \no{qw4} we can make integrations by parts and see
that the contribution from the remainder to \no{qw.72} becomes
$$
\iint e^{{i\over h}\phi }\widetilde{{\cal O}}(h^{-m-{3n\over 2}}\langle
\alpha _\xi \rangle ^{k+{n\over 2}}h\langle \alpha _\xi \rangle \inv )u(y)
dyd\alpha .
$$
Combining this with \no{qw.75} and Lemma \ref{qw5}, we get the
proposition.
\end{proof}

\par We can now prove Proposition \ref{qw1}.
\begin{proof}
We only consider $A$
and we may assume that $A=A(\delta )$ where $A(s)$ is a smooth \fy{} as in
the preceding proposition. The result for $B$ will be the same since $B$ is
like the adjoint of $A$. From \no{qw.70} we get
$$
h\partial _sA_s^*=\check{\psi }_s^*A_s^*.
$$
If $u\in{\cal S}$, then $A_s^*u\in{\cal S}$ and we get
$$
h\partial _s\Vert A_s^*u\Vert ^2=(\check{\psi }_s^*A_s^*u\vert
A_s^*u)+(A_s^*u\vert \check{\psi }_s^*A_s^*u),
$$
so
$$
h\partial _s\Vert A_s^*u\Vert ^2\le 2\Vert \check{\psi}_s^*\Vert \Vert
A_s^*u\Vert ^2.
$$
But $\Vert \check{\psi }_s^*\Vert ={\cal O}(\epsilon )$, so
$$
\Vert A_s^*u\Vert \le e^{{\cal O}(\epsilon )\vert s\vert /h}\Vert u\Vert
=e^{{\cal O}(1)A\vert s\vert }\Vert u\Vert \quad (\epsilon =Ah).
$$
Then
$\Vert A_\delta ^*\Vert \le e^{{\cal O}(1)A\delta }$ and the Proposition
follows since $A_\delta $ and $A_\delta ^*$ have the same norm.
\end{proof}

\par Let $P={\cal O}(\langle \xi \rangle ^M)$ in the symbol sense:
$\partial _x^\alpha \partial _\xi ^\beta P={\cal O}(\langle \xi \rangle
^{M-\vert \beta \vert })$. We study $PA_s-A_sP$ where $A_s$ is as in
Proposition \ref{qw6}, $\partial _s^\ell a_s$ \bdd{} in $S^{m+{3n\over
2},k+{n\over 2}}$ (so no $\sqrt{\epsilon }$-loss in the symbol for
simplicity). Working with the Weyl quantization we know that
\begin{eqnarray*}
P(x,hD_x;h)(e^{{i\over h}\phi (x,y,\alpha )}a(x,y,\alpha ;h))&=&e^{{i\over
h}\phi }b+\widetilde{{\cal O}}(h^\infty \langle \alpha _\xi \rangle
^{-\infty })\\
{\trans{P}}(y,hD_y;h)(e^{{i\over h}\phi (x,y,\alpha )}a(x,y,\alpha ;h))&=&e^{{i\over
h}\phi }c+\widetilde{{\cal O}}(h^\infty \langle \alpha _\xi \rangle
^{-\infty }),
\end{eqnarray*}
where
\begin{eqnarray*}
b&\equiv& P(x,\phi '_x(x,y,\alpha );h)a(x,y,\alpha ;h)+
\widetilde{{\cal O}}(h^{-{3n\over 2}-m}\langle
\alpha _\xi \rangle ^{k+{n\over 2}+M}h^2\langle \alpha _\xi \rangle ^{-2}),\\
c&\equiv& P(y,-\phi '_y(x,y,\alpha );h)a(x,y,\alpha ;h)+
\widetilde{{\cal O}}(h^{-{3n\over 2}-m}\langle
\alpha _\xi \rangle ^{k+{n\over 2}+M}h^2\langle \alpha _\xi \rangle ^{-2})
,
\end{eqnarray*}
and the $\widetilde{{\cal O}}$ refers to $\sqrt{\epsilon }$-loss only
\wrt{} differentiations in $\alpha $. (Recall the general fact that the
Weyl symbol of $e^{-i\phi /h}\circ p^w(x,hD_x)\circ e^{i\phi /h}$ is equal to
$p(x,\xi +\phi '_x)+{\cal O}(h^2)$.)

\par We conclude that $PA_s-A_sP=\widetilde{A}_s+A_s\circ {\rm
Op\,}(\widetilde{{\cal O}}(h^2\langle \alpha _\xi \rangle ^{-2}))$, where
\ekv{qw.76}
{
\widetilde{A}_su=\iint e^{{i\over h}\phi }(P(x,\phi '_x(x,y,\alpha
);h)-P(y,-\phi '_y;h))au(y)dyd\alpha .
}
On the critical set \no{qw.36}, we have
$$
P(x,\phi '_x;h)-P(y,-\phi '_y;h)=P(\kappa (s)(\alpha ))-P(\alpha ),
$$
and more generally,
\begin{eqnarray*}
P(x,\phi '_x)-P(y,-\phi '_y)&=&P(\kappa (s)(\alpha ))-P(\alpha )\\&&+(x-\kappa
(s,\alpha )_x)\cdot ({\cal O}(\langle \alpha _\xi \rangle
^M)+\widetilde{{\cal O}}(s\epsilon ^{1\over 2}\langle \alpha _\xi
\rangle ^{M-1}))\\
&&+(y-\alpha _x)\cdot ({\cal O}(\langle \alpha _\xi \rangle
^M)+\widetilde{{\cal O}}(s\epsilon ^{1\over 2}\langle \alpha _\xi
\rangle ^{M-1})).
\end{eqnarray*}
Using \no{qw.64}, \no{qw.65}, we get
\eekv{qw.77}
{
P(x,\phi '_x)-P(y,-\phi '_y)&=&P(\kappa (s,\alpha ))-P(\alpha )}{&&+(\phi
'_{\alpha _\xi },\langle \alpha _\xi \rangle \inv \phi '_{\alpha
_x})\cdot ({\cal O}(\langle \alpha _\xi \rangle ^M)+\widetilde{{\cal
O}}(s\epsilon ^{1\over 2}\langle \alpha _\xi \rangle ^{M-1}).
}

\par The contribution from $P(\kappa (s)(\alpha ))-P(\alpha )$ in \no{qw.76}
is of the form
$$
A_s\circ {\rm Op\,}(P(\kappa (s)(\alpha ))-P(\alpha )+\widetilde{{\cal
O}}(s{h\over \sqrt{\epsilon }}\langle \alpha _\xi \rangle ^{M-2}))
$$
in general, and the remainder estimate improves to
$\widetilde{{\cal O}}(sh\langle \alpha _\xi \rangle ^{M-2})$ if $P(\kappa
(s,\alpha ))-P(\alpha )=\widetilde{{\cal O}}(s\epsilon \langle \alpha
\rangle ^{M-1})$. By integration by parts in $\alpha $, we see that the
contribution from the two remainder terms in \no{qw.77} is
\begin{eqnarray*}
h\iint e^{{i\over h}\phi }\partial _{\alpha _\xi }(({\cal O}(\langle
\alpha _\xi \rangle ^M)+\widetilde{{\cal O}}(s\epsilon ^{1\over 2}\langle
\alpha _\xi \rangle ^{M-1})+\partial _{\alpha _x}({1\over \langle \alpha
_\xi \rangle }({\rm idem}))a)u(y)dyd\alpha \\
=h\iint e^{{i\over h}\phi }({\cal O}(\langle \alpha _\xi \rangle
^{M-1})+\widetilde{{\cal O}}(s \langle \alpha _\xi \rangle
^{M-2}))au(y)dyd\alpha .
\end{eqnarray*}
In conclusion,
\begin{prop}\label{qw7}
We have
$$
PA_s=A_sP_s,\ P_s=P+Q_s,
$$
where
$$
Q_s=P(\kappa (s)(\alpha ))-P(\alpha )+\widetilde{{\cal O}}(s{h\over
\epsilon ^{1\over 2}}\langle \alpha _\xi \rangle ^{M-2})+\widetilde{{\cal
O}}(h\langle \alpha _\xi \rangle ^{M-1})
$$
in the general case, when $P(\kappa (s)(\alpha ))-P(\alpha
)=\widetilde{{\cal O}}(s\epsilon^{1\over 2} \langle \alpha _\xi \rangle
^{M-1})$.

\par In the special case when $P(\kappa (s,\alpha ))-P(\alpha
)=\widetilde{{\cal O}}(s\epsilon \langle \alpha _\xi \rangle ^{M-1})$ the
first remainder term improves to $\widetilde{{\cal O}}(sh\langle \alpha
_\xi \rangle ^{M-2})$ so in that case,
$$
Q_s(\alpha )=P(\kappa (s,\alpha ))-P(\alpha )+\widetilde{{\cal
O}}(h\langle \alpha _\xi \rangle ^{M-1}).
$$
\end{prop}

\par In order to treat certain conjugations, we need a more precise
description of $Q_s$ in the general case in the last proposition. The
proof above shows that
\ekv{qw.78}
{
A_sQ_su=\iint e^{{i\over h}\phi }(P(\kappa (s,\alpha ))-P(\alpha
))au(y)dyd\alpha +A_s{\rm Op\,}(\widetilde{{\cal O}}(h\langle \xi
\rangle ^{M-1})),
}
and we need to take a closer look at the oscillatory integral.

\par By Taylor's formula and \no{geo.28},
$$
P(\kappa (\alpha ))-P(\alpha )=\widetilde{{\cal O}}(s\sqrt{\epsilon
})\cdot (P_{\alpha _\xi }'(\alpha ),{1\over \langle \alpha _\xi \rangle
}P'_{\alpha _x}(\alpha ))+\widetilde{{\cal O}}(s^2\epsilon \langle \alpha
_\xi \rangle ^{M-2}).
$$
When passing to Weyl composition of symbols, we notice that if
$r=\widetilde{{\cal O}}(s\sqrt{\epsilon })$, then
$$
r\# (P'_{\alpha _\xi },{1\over \langle \alpha _\xi \rangle }P'_{\alpha
_x})=r\cdot (P'_{\alpha _\xi },{1\over \langle \alpha _\xi \rangle }P'_{\alpha
_x})+\widetilde{{\cal O}}(sh\langle \alpha _\xi \rangle ^{M-2}).
$$
Lemma \ref{qw5} and Remark \ref{qw5'} then show that on the
symbol level
\ekv{qw.79}
{
Q_s=\widetilde{{\cal O}}(h\langle \xi \rangle ^{M-1})+\widetilde{{\cal
O}}(s^2\epsilon \langle \xi \rangle ^{M-2})+\widetilde{{\cal
O}}(s\sqrt{\epsilon })\# (P'_\xi ,{1\over \langle \xi \rangle }P'_x).
}
Here the first term to the right is too large; we would like to have
$h\langle \xi \rangle ^{M-2}$, so we take a closer look at $Q_s$ using
Proposition \ref{qw6}, where we now add the assumption
\ekv{qw.80}
{
m=k=0,\ A_0=1,\ a_s\hbox{ is \indep{} of }s.
}
Then as noticed, \no{qw.71} improves to
\ekv{qw.81}
{
\check{\psi }_s=\widehat{\psi }_s\circ \kappa (s)+\widetilde{{\cal
O}}({h\over \langle \xi \rangle })=\widetilde{{\cal O}}(\epsilon ).
}
We get on the symbol level, writing $[A,B]=A\# B-B\# A$:
\begin{eqnarray*}
h\partial _sQ_s&=&[P_s,\check{\psi }_s]=[P,\check{\psi }_s]+[Q_s,\check{\psi
}_s]\\
&=&[P,\check{\psi }_s]+\widetilde{{\cal O}}(h\langle \xi \rangle
^{M-1}\epsilon {h\over \epsilon \langle \xi\rangle  })+\widetilde{{\cal
O}}(s^2\epsilon \langle \xi \rangle ^{M-2}\epsilon {h\over \epsilon
\langle \xi \rangle })\\
&&+\widetilde{{\cal O}}(s\sqrt{\epsilon }\epsilon {h\over \epsilon
\langle \xi \rangle })\# (P'_\xi ,{1\over \langle \xi \rangle
}P'_x)+\widetilde{{\cal O}}(s\sqrt{\epsilon })\# \widetilde{{\cal
O}}(\epsilon \langle \xi \rangle ^{M-1}){h\over \sqrt{\epsilon }\langle
\xi \rangle })\\
&=&
[P,\check{\psi }_s]+\widetilde{{\cal O}}(h^2\langle \xi \rangle
^{M-2})+\widetilde{{\cal O}}(s^2\epsilon h\langle \xi \rangle
^{M-3})
\\&&+\widetilde{{\cal O}}({s\sqrt{\epsilon }h\over \langle \xi
\rangle })\# (P'_\xi ,{1\over \langle \xi \rangle }P'_x)+\widetilde{{\cal
O}}(s\epsilon h\langle \xi \rangle ^{M-2}).
\end{eqnarray*}
On the other hand,
\begin{eqnarray*}
[P,\check\psi
_s]&=&{h\over i}\{ P,\check{\psi }_s\} +\widetilde{{\cal O}}(\langle \xi
\rangle ^M\epsilon {h^3\over \langle \xi \rangle ^3\epsilon ^{3/2}})\\
&=&
ih\{ \widehat{\psi }_s\circ \kappa (s),P\} +h
(P'_\xi ,{P'_x\over \langle \xi \rangle })\cdot
\widetilde{{\cal O}}({h\over \langle \xi
\rangle \epsilon ^{1/2}})+\widetilde{{\cal O}}(\epsilon ^{-{1\over
2}}h^3\langle \xi \rangle ^{M-3})\\
&=&
ih\{ \widehat{\psi }_s\circ \kappa (s),P\} +h\widetilde{{\cal O}}({h\over
\langle \xi \rangle \epsilon ^{1/2}})\# (P'_\xi ,{P'_x\over \langle \xi
\rangle })+\widetilde{{\cal O}}({h^3\over \epsilon }\langle \xi \rangle
^{M-3}),
\end{eqnarray*}
so
$$
h\partial _sQ_s=ih\{ \widehat{\psi }_s\circ \kappa (s),P\}
+\widetilde{{\cal O}}({h^2\over \epsilon ^{1/2}\langle \xi \rangle
}+{s\epsilon ^{1/2}h\over \langle \xi \rangle })\# (P'_\xi ,{P'_x\over
\langle \xi \rangle })+\widetilde{{\cal O}}(h^2\langle \xi \rangle
^{M-2}+s\epsilon h\langle \xi \rangle ^{M-2}).
$$
Since $\kappa (s)$ is obtained from integrating $iH_{\widehat{\psi }_s}$, we
obtain, using also that $Q_0=0$,
\ekv{qw.82}
{
Q_\delta =P\circ \kappa (\delta )-P+\widetilde{{\cal O}}({\delta
^2\sqrt{\epsilon }\over \langle \xi \rangle }+{\delta h\epsilon
^{-1/2}\over \langle \xi \rangle })\#(P'_\xi ,{P'_x\over \langle \xi
\rangle })+\widetilde{{\cal O}}(h\delta \langle \xi \rangle ^{M-2}+\delta
^2\epsilon \langle \xi \rangle ^{M-2})
}
under the assumption \no{qw.80}. (Recall that here $M$ is the order of $P$.)

\section{The conjugated \pop{}}\label{co}
\setcounter{equation}{0}

Let $P$ be the \op{} introduced in the beginning of Section \ref{ap}
satisfying the assumptions  \no{ap.1}--\no{ap.7}. To start with, we assume
Hypothesis \ref{ny0} and
define ${\cal C},\psi _\epsilon $ as in \no{ny.2}, \no{ny.3}. Later, we
shall also use the dynamical Hypothesis \ref{ny2} (implying \no{ny.19}).

\par The \fu{} $\psi _\epsilon $ satisfies $\psi _\epsilon
=\widetilde{{\cal O}}(\epsilon )$, $\xi \cdot \partial _\xi \psi
_\epsilon =\widetilde{{\cal O}}(\epsilon \langle \xi \rangle ^{-2})$ in
the sense of \no{geo.8}, \no{geo.10.5}, but we also know that these bounds
improve considerably away from ${\cal C}$.

\par Let $A=A_\delta $ be an elliptic \fop{} of order 0,0 quantizing
$\kappa (\delta )$ as in Section \ref{qw}, and let $B_\delta $ be a
corresponding \op{} quantizing $\kappa (\delta )\inv$. Then we can apply
Proposition \ref{qw2} to see that $BA={\rm Op}_h(c)$, where
$c=\widetilde{{\cal O}}(1)$ is elliptic. Now restrict the attention to the
intermediate region, where $\psi _\epsilon =\epsilon \psi (x)$, $\psi
(x)\in C^\infty $.
Using that $\epsilon =Ah$ with $A$ large but fixed, we find here that
$$
e^{-\delta \psi _\epsilon /h}\circ A={\rm Op}_h(c)
$$
where $c={\cal O}(1)$ in the standard symbol sense, $\partial _x^\alpha
\partial _\xi ^\beta c={\cal O}(\langle \xi \rangle ^{-\vert \beta
\vert })$ and $c$ is elliptic in this class. Let $\widetilde{c}={\cal
O}(1)$ be an elliptic symbol which is equal to 1 in the interior region
and which changes from 1 to $c$ when going outwards in the intermediate
region. We now replace $A$ by $A\circ {\rm
Op}_h(\widetilde{c})\inv$. (Here we assume that $A_0=1$, to avoid a
topological difficulty.) Then
we have achieved that $Au=e^{\delta \psi _\epsilon /h}u$ for $u$ supported
in the exterior region including a shell where $\psi _\epsilon =\epsilon
\psi $, $\psi =\psi (x)\in C^\infty $ as above and the region further out
where $\psi _\epsilon =\psi _\epsilon (x)=\widetilde{{\cal O}}(\epsilon )$.

Next, we consider $P^\delta $ defined as in Proposition \ref{qw7} by
$PA_\delta =A_\delta P^\delta  $. Recall that $P=P_2+iP_1+P_0$ and write
$P_jA_\delta =A_\delta P_j^\delta $, $P^\delta =P+Q^\delta $, $P_j^\delta
=P_j+Q_j^\delta $.

\par For $P_j$, $j=0,2$, we apply \no{qw.82} and get on the symbol level
over any \bdd{} set in $x$-space
\begin{eqnarray}\label{co.1}
Q_j^\delta &=& p_j\circ \kappa (\delta )-p_j+\widetilde{{\cal O}}({\delta
^2\sqrt{\epsilon }\over \langle \xi \rangle }+{\delta h\epsilon ^{-{1\over
2}}\over \langle \xi \rangle })\#(\partial _\xi p_j,{\partial _xp_j\over
\langle \xi\rangle }) +\widetilde{{\cal O}}(h\delta \langle \xi \rangle
^{j-2}+\delta ^2\epsilon \langle \xi \rangle ^{j-2})\nonumber
\\
&=&
\widetilde{{\cal O}}(\delta \sqrt{\epsilon })\# (\partial _{\xi
}p_j,{\partial _xp_j\over \langle \xi \rangle })+\widetilde{{\cal
O}}((h\delta +\delta ^2\epsilon )\langle \xi \rangle ^{j-2}).
\end{eqnarray}

\par For $p_1$ we need the more precise information about $\psi _\epsilon $
given in \no{ny.7}, \no{ny.9}, \no{ny.10}, satisfied also by $\psi
_\epsilon $, as well as the fact that $\psi _\epsilon (\rho )$ is \indep{}
of $\epsilon $ in $\vert \rho -\rho _j\vert \le \sqrt{\epsilon }$
(actually in $\vert \rho -\rho _j\vert \le \sqrt{\epsilon }/C$ for some
$C>0$, but we can always dilate in $\epsilon $). From these estimates it follows
that
\ekv{co.2}
{
({\partial _xp_1\over \langle \xi \rangle },\partial _\xi p_1)\otimes
(\partial _x\psi _\epsilon ,\langle \xi \rangle \partial _\xi \psi
_\epsilon )=\widetilde{{\cal O}}(\epsilon )
}
over a \neigh{} of the set where $\psi _\epsilon $ also depends on $\xi $.
It also follows that
\ekv{co.3}
{
p_1(\kappa (\delta)(x,\xi ))-p_1(x,\xi )=\widetilde{{\cal O}}(\delta
\epsilon ).
}
From Proposition \ref{qw7}, we deduce that
\ekv{co.4}
{
Q_1^\delta =p_1(\kappa (\delta )(x,\xi ))-p_1(x,\xi )+\widetilde{{\cal
O}}(h).
}
By construction of $\kappa (\delta )$, especially \no{geo.44}, we see that
\ekv{co.4.5}
{p_1(\kappa (\delta )(x,\xi ))-p_1(x,\xi )=i\delta H_{\psi _\epsilon
}p_1(x,\xi )+\widetilde{{\cal O}}(\epsilon \delta ^2).}
In fact, this is quite obvious in the intermediate and far out
regions, and near a point $\rho_j\in {\cal C}$, we choose canonical
coordinates so that $\rho _j=(0,0)$ and put $\rho =\sqrt{\epsilon
}\widetilde{\rho }$. Then if $\psi _\epsilon (\rho )=\epsilon
\widetilde{\psi }_\epsilon (\widetilde{\rho })$, we get $H_{\psi
_\epsilon }=\widetilde{H}_{\widetilde{\psi }_\epsilon }$, where the
tilde on the $H$ indicates that we take the Hamilton field \wrt{} the
$\widetilde{\rho }$-variables. The \mfld{} $\Lambda _\delta :$ $\Im
\rho =\delta H_{\psi _\epsilon }(\Re \rho )$ now becomes
$\Im\widetilde{\rho }=\delta \widetilde{H}_{\widetilde{\psi }_\epsilon
}(\Re \widetilde{\rho })$ while $p_1={\cal O}(\epsilon \widetilde{\rho
}^2)$ and since $\psi _\epsilon $ satisfies \no{ny.7},
$$
\partial ^\alpha _{\widetilde{\rho }}\widetilde{\psi }_\epsilon
={\cal O}({1\over \langle \widetilde{\rho }\rangle ^{|\alpha |}}).
$$
Then
$\kappa (\delta )(\rho )$ becomes $\widetilde{\kappa }(\delta
)(\widetilde{\rho })=\widetilde{\rho }+i\delta
\widetilde{H}_{\widetilde{\psi }_\epsilon }(\widetilde{\rho })+{\cal
O}(\delta ^2\langle \widetilde{\rho }\rangle ^{-2})$, so
$p_1(\widetilde{\kappa }(\delta )\widetilde{\rho }))=p_1(\rho
)+i\delta \widetilde{H}_{\widetilde{\psi }_\epsilon }(p_1)+{\cal
O}(\epsilon \delta ^2)$, leading to \no{co.4.5}.

\par It follows that
\ekv{co.5}
{
Q_1^\delta =-i\delta H_{p_1}\psi _\epsilon +\widetilde{{\cal O}}(\delta
^2\epsilon +h)=\widetilde{{\cal O}}(\delta \epsilon +h).
}
Away from any \neigh{} ${\cal C}$ we have the improved estimates
$$
Q_1^\delta =-i\delta H_{p_1}\psi _\epsilon +{\cal O}(\delta ^2\epsilon +h)
$$
in the usual symbol sense, as long as we stay away from the outer region
where $\psi _\epsilon =\psi _\epsilon (x)$ only satisfies $\psi _\epsilon
=\widetilde{{\cal O}}(\epsilon )$, and where we can apply the analysis of
Section \ref{ap}, so
\ekv{co.6}
{
Q_1^\delta =-i\delta H_{p_1}\psi _\epsilon =-i\delta \nu (\psi _\epsilon )
}
there.

\section{Estimates for the conjugated \pop{} and localization of the
spectrum}
\label{app}
\setcounter{equation}{0}

\par Let $P^\delta $ be the conjugated \op{} of Section \ref{co}. We shall
study lower bounds for
\ekv{app.1}
{
\Re (P^\delta u\vert u)=\Re ((P_0^\delta +P_2^\delta )u\vert u)-\Im
(P_1^\delta u\vert u),\ u\in{\cal S}(M).
}
Using \no{co.1}, we get
\begin{eqnarray*}
\Re (P_2^\delta u\vert u)&\ge& (P_2u\vert u)-C\delta \sqrt{\epsilon }\Vert u\Vert
\Vert {\rm Op\,}((\partial _\xi p_2,{\partial _xp_2\over \langle \xi
\rangle }))u\Vert -C(h\delta +\delta ^2\epsilon )\Vert u\Vert ^2\\
&\ge & (P_2u\vert u)-{1\over C_1}\Vert {\rm Op\,}((\partial _\xi
p_2,{\partial _xp_2\over \langle \xi \rangle }))u\Vert ^2-C_2(h\delta
+\delta ^2\epsilon )\Vert u\Vert ^2,
\end{eqnarray*}
where $C_1$ can be chosen \ably{} large and $C_2$ depends on $C_1$. Here
$$
{1\over 2}p_2-{1\over C_1}(\partial _\xi p_2,{\partial_x p_2\over \langle
\xi \rangle })^2\ge 0,
$$
if $C_1$ is large enough, so
\begin{eqnarray*}
(P_2u\vert u)-{1\over C_1}\Vert {\rm Op\,}((\partial _\xi p_2,{\partial_x
p_2\over \langle \xi \rangle }))u\Vert ^2 & = &
({\rm Op\,}(p_2-{1\over C_1}(\partial _\xi p_2,{\partial_x
p_2\over \langle \xi \rangle })^2)u\vert u)+{\cal O}(h^2)\Vert u\Vert ^2\\
&\ge & {1\over 2}(P_2u\vert u)-{\cal O}(h^2)\Vert u\Vert ^2,
\end{eqnarray*}
where the estimate follows from the Fefferman-Phong inequality in the
semi-classical setting (see \cite{HeSjSt}). It follows that
\ekv{app.2}
{
\Re (P_2^\delta u\vert u)\ge {1\over 2}(P_2u\vert u)-C(h^2+h\delta +\delta
^2\epsilon )\Vert u\Vert ^2.
}
For $P_0$ we have the same conclusion,
\ekv{app.3}
{
\Re (P_0^\delta u\vert u)\ge {1\over 2}(P_0u\vert u)-C(h^2+h\delta +\delta
^2\epsilon )\Vert u\Vert ^2,
}
since the same analysis applies for $x$ in a \bdd{} region and further
out, we just have $P_0^\delta =P_0\ge 0$.

\par Combining this with \no{co.5}, we get
\ekv{app.4}
{
\Re (P^\delta u\vert u)\ge \delta ({\rm Op\,}(H_{p_1}\psi _\epsilon )u\vert
u)+{1\over 2}((P_0+P_2)u\vert u)-C(h+h\delta +\delta ^2\epsilon )\Vert
u\Vert ^2.
}
Here we recall that $H_{p_1}\psi _\epsilon =\langle \widetilde{p}_\epsilon
\rangle _{T_0}-\widetilde{p}_\epsilon $ by \no{ny.4}, that $\langle
\widetilde{p}_\epsilon \rangle _{T_0}$ satisfies \no{ny.21}, \no{ny.22} and
that $\widetilde{p}_\epsilon \le \widetilde{p} \le p_0+p_2$ by \no{ny.12.5}. Write,

\eekv{app.5}
{
\delta H_{p_1}\psi _\epsilon +{1\over 2}(p_0+p_2)&=&\delta \langle
\widetilde{p}_\epsilon \rangle _{T_0}-\delta \widetilde{p}_\epsilon
+{1\over 2}(p_0+p_2)
}
{&=&\delta \langle \widetilde{p}_\epsilon \rangle _{T_0}+\delta
(\widetilde{p}-\widetilde{p}_\epsilon )+\delta
(p_0+p_2-\widetilde{p})+({1\over 2}-\delta )(p_0+p_2).
}
Here we want a lower bound for ${\rm
Op\,}(\widetilde{p}-\widetilde{p}_\epsilon )$. This is quite straight
forward away from ${\cal C}$, so we concentrate on a \neigh{} of a point
$\rho _j\in{\cal C}$. Assume $\rho _j=0$  for simplicity. Then near $0$ we
have by \no{ny.12}
$$
\widetilde{p}_\epsilon =g({\vert \rho \vert ^2\over \epsilon
})\widetilde{p},\quad \widetilde{p}-\widetilde{p}_\epsilon =(1-g({\vert \rho
\vert ^2\over \epsilon }))\widetilde{p}=(1-\ell_\epsilon (\rho
))^2\widetilde{p},
$$
where
$$
\ell_\epsilon (\rho )=\ell ({\vert \rho \vert ^2\over \epsilon
})=\widetilde{{\cal O}}({\epsilon \over \epsilon +\vert \rho \vert ^2}),
$$
and we may assume that $g$ has been chosen so that $1-g=(1-\ell )^2$ with
$\ell$ smooth.
Here
\begin{eqnarray*}
(1-\ell_\epsilon )\# \widetilde{p}\# (1-\ell_\epsilon
)&=&\widetilde{p}-\ell_\epsilon
\#\widetilde{p}-\widetilde{p}\#\ell_\epsilon +\ell_\epsilon
\#\widetilde{p}\#\ell_\epsilon \\
&=&(1-\ell_\epsilon (\rho ))^2\widetilde{p}+\widetilde{{\cal O}}(h).
\end{eqnarray*}
From this we conclude that
\ekv{app.6}
{
({\rm Op\,}(\widetilde{p}-\widetilde{p}_\epsilon )u\vert u)\ge -Ch\Vert
u\Vert ^2.
}
Similarly,
\ekv{app.7}
{
({\rm Op\,}(p_0+p_2-\widetilde{p})u\vert u)\ge -Ch^2\Vert u\Vert ^2,
}
by the Fefferman-Phong inequality or by a direct argument.

\par Let $0\le k_\epsilon =\widetilde{{\cal O}}(\epsilon )$ be equal to
$\epsilon $
in ${\cal C}+B(0,\sqrt{\epsilon })$ and have its support in ${\cal
C}+B(0,\sqrt{2\epsilon })$. Let $K_\epsilon ={\rm Op\,}(k_\epsilon )$. By
\no{ny.21}, \no{ny.22} we have $\langle \widetilde{p}_\epsilon \rangle
_{T_0}\backsim {\rm min\,}(\epsilon ,{\rm dist\,}(\cdot ,{\cal C})^2)$, so
\ekv{app.8}
{
\epsilon \backsim k_\epsilon +\langle \widetilde{p}_\epsilon \rangle
_{T_0}=\widetilde{{\cal O}}(\epsilon ).
}
Hence,
\ekv{app.9}
{
((K_\epsilon +{\rm Op\,}(\langle \widetilde{p}_\epsilon \rangle
_{T_0}))u\vert u)\backsim \epsilon \Vert u\Vert ^2.
}
Combining \no{app.4}--\no{app.9}, we get
\begin{prop}\label{app1}
We have
\ekv{app.10}
{
\Re ((P^\delta +K_\epsilon )u\vert u)\ge {\delta \epsilon \over C}\Vert
u\Vert ^2+({1\over 2}-\delta )((P_2+P_0)u\vert u)-Ch\Vert u\Vert ^2,\
u\in {\cal S}(M).
}
\end{prop}

\par Here we recall that $\delta >0$ should be small enough, $\epsilon
=Ah$ with $A$ \ably{} large and fixed, and $C$ in \no{app.10} is
\indep{} of $\delta ,A$ while $h$ is small enough depending on these
two parameters. From \no{app.10} we get the a priori estimate
\ekv{app.11}
{
({\delta Ah\over C}-Ch-\Re z)\Vert u\Vert\le \Vert (P^\delta
+K_\epsilon -z)u\Vert,\ u\in {\cal S}(M),
}
when $\Re z<\delta Ah/C-Ch$.

\par From Section \ref{macr} we know that $P$ has no spectrum in the
open left half-plane. We shall next prove
\begin{prop}\label{app2}
For every constant $B>0$ there is a constant $D>0$ such that $P$ has
no spectrum in
\ekv{app.12}
{
\{ z\in {\bf C};\, \Re z<Bh,\, |\Im z|>Dh\}
}
when $h>0$ is small enough. Moreover $\Vert (P-z)\inv\Vert ={\cal O}_B(h\inv )$ for $z$ in the set {\rm \no{app.12}}.
\end{prop}
\begin{proof}
Choose $\delta >0$ small, then $A$ large enough, so that ${\delta
Ah\over C}-Ch-\Re z\ge \epsilon /{\rm Const.\,}$ when $\Re z<Bh$. Then \no{app.11} gives
\ekv{app.13}
{
{\epsilon \over C_0}\Vert u\Vert\le \Vert (P^\delta +K_\epsilon
-z)u\Vert,\ u\in{\cal S}(M).
}
Take $z$ in the set \no{app.12}. When $\Re z<0$, we already know that
$z\not\in \sigma (P)$, so we may assume that $0\le \Re z <Bh$.

\par Now recall that the symbol of $K_\epsilon $ is $\widetilde{{\cal
O}}(\epsilon )$ and supported in ${\cal C}+B(0,\sqrt{2\epsilon })$. On
that set we have $p^\delta =\widetilde{{\cal O}}(\epsilon
)=\widetilde{{\cal O}}(Ah)$ and hence
$$
|p^\delta -z|>{\epsilon +|z|\over C_0},\hbox{ when }|\Im z|>Dh,
$$
and $D$ is large enough, assuming still that $0\le \Re z<Bh$.

\par It follows that we can find $E=\widetilde{{\cal O}}(\epsilon
/(\epsilon +|z|))$ such that
$$
K_\epsilon =E\circ (P^\delta -z)+F,\quad F=\widetilde{{\cal O}}(h),
$$
where $E,F$ also denote the corresponding $h$-\pop{}s. In particular,
\ekv{app.14}
{
\Vert K_\epsilon u\Vert\le {\cal O}(1)\Vert (P^\delta -z)u\Vert +{\cal
O}(h)\Vert u\Vert,
}
where the ${\cal O}$s are uniform in $\epsilon $. Combining this with
\no{app.13} with $\epsilon \gg h$, we get
\ekv{app.15}
{
\epsilon \Vert u\Vert \le {\cal O}(1)\Vert (P^\delta -z)u\Vert,\ u\in
{\cal S}(M),
}
for $0\le \Re z< Bh$, $0\le \Im z<Dh$. Now recall that
$P^\delta =A_\delta \inv PA_\delta $, where $A_\delta ,\, A_\delta
\inv :{\cal S}\to {\cal S},\ L^2\to L^2$ and have $L^2$ norm ${\cal
O}_A(1)$. Then \no{app.15} gives
\ekv{app.16}
{
h\Vert u\Vert\le {\cal O}_B(1)\Vert (P-z)u\Vert, \ u\in {\cal S}(M),
}
for $z$ in the set \no{app.12}. From Section \ref{macr} we then know
that $z\not\in \sigma (P)$ and that $\Vert (P-z)\inv \Vert\le {\cal O}_B(1)/h$.
\end{proof}

\section{Asymptotics of eigenvalues}\label{as}
\setcounter{equation}{0}

\par \par Let $\rho _j\in{\cal C}$ and let $F_{\rho _j}$ be the matrix of the
linearization of $H_p$ at $\rho _j$ (the so called fundamental matrix of $p$
at the doubly characteristic point $\rho _j$). Thanks to the fact that the
quadratic approximation of $p^\delta =p\circ \kappa (\delta ) $ at
$\rho _j$ is elliptic on $T_{\rho
_j}(\Lambda _\delta )$ and takes its values in a closed angle contained in the
union of $\{ 0\}$ and the open right half plane, we know from \cite{Sj} that
the \ev{}s of $F_{\rho _j}$ are of the form $\pm \lambda _{j,k}$, $1\le
k\le n$, when repeated with their multiplicity, with $\Im \lambda
_{j,k}>0$. Let
\ekv{as.1}
{
{\widetilde{\rm tr}\,}(p,\rho _j)={1\over i}\sum_k \lambda _{j,k}.
}
In our case the subprincipal symbol of $P$ at $\rho _j$ is zero and will
not enter into the description of the \ev{}s.

\par Put
$$
q(x,\xi )=-p(x,i\xi )=p_2+p_1-p_0.
$$
Let $F_q$, $F_p$ be the fundamental matrices of $q,p$ at one of the
critical points $\rho _j\in{\cal C}$. Since
$$
H_q(x,\xi )={1\over i}(p'_\xi (x,\eta )\cdot {\partial \over \partial x}-
p'_x (x,\eta )\cdot {\partial \over \partial \eta  }),\hbox{ with }\eta
=i\xi ,
$$
we see that $F_q$ and ${1\over i}F_p$ have the same \ev{}s; $\pm {1\over
i}\lambda _k$, $k=1,...,n$ ($j$ being fixed) where $\Re ({1\over i}\lambda
_k)>0$. Now $q$ is real-valued and we can apply the stable \mfld{} theorem
as in \cite{HeSj} (and at many other places) to see that the $H_q$-flow has
a stable outgoing \mfld{} $\Lambda _+$ passing through $\rho _j$ such
that $T_{\rho _j}\Lambda _+^{\bf C}$ is spanned by the generalized
eigenvectors corresponding to $+{1\over i}\lambda _k$, $k=1,...,n$. We also
know that $\Lambda _+$ is a Lagrangian \mfld{} and that $q$ vanishes on
$\Lambda _+$.

\begin{lemma}\label{as3}
Assume for simplicity that $\rho _j=(0,0)$. Then $T_{\rho _j}\Lambda _+$ is
transversal both to $\{ x=0\}$ and to $\{ \xi =0\}$.
\end{lemma}
\begin{proof}
Since we consider the linearized situation we may assume right away that $q$
is a quadratic form (the second order Taylor polynomial at $\rho _j$), so
that $p_j$
are quadratic forms as well. The dynamical condition \no{ny.17} implies that
\ekv{as.10.5}
{
\langle p_0+p_2\rangle _{T_0}>0 \hbox{ as a quadratic form}.
}

\par Let $L=\Lambda _+\cap \{ x=0\}$. Since $p_0=0$, $p_1=0$ on $\{
x=0\}$, we know that $p_2=0$ on $L$ and hence $H_q=H_{p_1}$ on $L$. Now $H_q$
is tangent to $\Lambda _+$
and $H_{p_1}$ is tangent to $\{ x=0\}$, so $H_q=H_{p_1}$
is tangent to $L$. Thus $L$ is an $H_q$- and $H_{p_1}$ invariant subspace
on which $p_2+p_0=0$, and since $\langle p_2+p_0\rangle _{T_0}>0$ away from
$0$, we
necessarily have $L=0$.

\par The proof of the fact that $\Lambda _+\cap \{ \xi =0\}=0$ is the same
after permuting the roles of $p_0$ and $p_2$.
\end{proof}

\par It follows from the lemma that
\ekv{as.11}
{
\Lambda _+:\ \xi =\phi _+'(x),\ x\in{\rm neigh\,}(0),
}
where $\phi _+\in C^\infty ({\rm neigh\,}(0);M)$, $\phi_+(0)=0$,
$\phi _+'(0)=0$, $\det \phi _+''(0)\ne 0$.

\par Let $\Lambda _-$ be the stable incoming $H_q$-invariant \mfld{} such
that $T_{\rho _j}\Lambda _-^{\bf C}$ is spanned by the generalized
eigenvectors of $F_q$ corresponding to $-{1\over i}\lambda _k$, $1\le k\le
n$. The lemma is valid also for $\Lambda _-$ and \no{as.11} has an obvious
analogue for $\Lambda _-$ where we let $\phi _-$ denote the corresponding
generating \fu{}.
\begin{prop}\label{as4}
We have $\phi _+''(0)>0$, $\phi _-''(0)<0$.
\end{prop}
\begin{proof}
Again we can consider the linearized quadratic case. If we make a smooth
deformation of $q$, then $\Lambda _+,\,\Lambda _-,\,\phi _+''(0),\,\phi
_-''(0)$ vary smoothly with the deformation parameter and $\det \phi
''_\pm (0)\ne 0$, provided of course that we maintain the condition
\no{as.10.5}. Consider the deformation from $q=q_0$ to to $\xi
^2-x^2=q_1$:
$$
q_t(x,\xi )=(1-t)q(x,\xi )+t(\xi ^2-x^2)=p_2^t+p_1^t-p_0^t,
$$
with
$$
p_2(\xi )=(1-t)p_2(\xi )+t\xi ^2,\ p_0^t(x)=(1-t)p_0(x)+tx^2.
$$
$p_2^t$ and $p_0^t$ are positive definite for $t>0$, so \no{as.10.5} is
maintained. For $t=1$ we have $\phi _\pm^1(x)=\pm x^2/2$ so $\pm \phi
_\pm ''(0)$ is positive definite. Since the signatures of $\phi _\pm ''(0)$
are \indep{} of $t$, we get the lemma.
\end{proof}

\par From Section \ref{co} we recall that the conjugated
operator $P_\delta =P_{\delta ,\epsilon }$ has the symbol
$$
p+\delta H_{p_1}\psi _\epsilon +\sum_{j=0,2}\widetilde{{\cal
O}}(\delta \sqrt{\epsilon })\#(\partial _{\xi
}p_j,{\partial _xp_j\over \langle \xi \rangle })+\widetilde{{\cal
O}}(h+\delta ^2\epsilon )
$$
and that we have the a priori estimate \no{app.10} expressing that the
real part of $P^\delta $ is $\ge {\delta \epsilon \over C}-Ch$ outside
${\cal C}+B(0,\sqrt{\epsilon })$. In the set ${\cal
C}+B(0,\sqrt{\epsilon })$ the symbol $P^\delta $ is independent of
$\epsilon $ modulo ${\cal O}(\epsilon \big({h\over \epsilon
}\big)^\infty )$ and is of the form $\widetilde{P}_\delta
\sim p^\delta +hr_1+h^2r_2+...$,
where $|p^\delta |\backsim {\rm dist\,}(\cdot ,{\cal C})^2$, $\Re
p^\delta \backsim {\rm dist\,}(\cdot ,{\cal C})^2$. In the following we shall
assume for a while (in order to simplify the notations) that ${\cal C}$ is
reduced to a single point $\rho _0=(0,0)$.

\par We fix $B>0$ and restrict the spectral parameter $z$ to the
disc $D(0,Bh)$. We shall take $\epsilon =Ah$ with $A\gg B$ sufficiently large.
Assume that $M={\bf R}^n$ for simplicity.
Following basically \cite{HeSjSt} we recall the construction of a well-posed
Grushin problem, first for $\widetilde{P}_\delta -z$ and then for $P_\delta
-z$. Let $\Lambda =(x^2+(hD_x)^2)^{1/2}$ so that $\widetilde{P}_\delta $ is
equipped with the natural domain ${\cal D}(\widetilde{P}_\delta )=
\{ u\in L^2;\, \Lambda ^2u\in L^2({\bf R}^n)\}$. In Section 11 of\cite{HeSjSt}
the authors
constructed operators
\ekv{gr.1}
{
R_-:{\bf C}^{N_0}\to L^2,\ R_+:L^2\to {\bf C}^{N_0}
}
of the form
\ekv{gr.2}
{
R_-u_-=\sum_{j=1}^{N_0}u_-(j)e_j^h(x),\ R_+u(j)=(u|f_j^h(x)),
}
with the following properties:
\ekv{gr.3}
{
e_j^h(x)=h^{-{n\over 4}}e_j({x\over\sqrt {h}}),\quad f_j^h(x)=
h^{-{n\over 4}}f_j({x\over\sqrt {h}}),
}
\ekv{gr.4}
{
e_j(x)=p_j(x)e^{i\Phi _0(x)},\quad f_j(x)=q_j(x)e^{i\Psi _0(x)},
}
where $p_j,q_j$ are polynomials and $\Phi _0,\Psi _0$ are quadratic forms
with $\Im \Phi _0,\, \Im \Psi _0>0$, $\Phi _0''=\phi _+''(0)$. If $\delta >0$
is small and fixed, $A$ sufficiently large, the problem
\ekv{gr.5}
{
(\widetilde{P}_\delta -z)u+R_-^\delta u_-=v,\ R_+^\delta u=v_+,
}
for $v\in L^2({\bf R}^n)$, $v_+\in {\bf C}^{N_0}$ has a unique solution
$u\in {\cal D}(\widetilde{P}_\delta )$, $u_-\in {\bf C}^{N_0}$, where
$R_+^\delta =R_+A_\delta $, $R_-^\delta =A_\delta \inv R_-$, Moreover, for the solution, we have the a priori estimate
\ekv{gr.6}
{
\Vert \Lambda ^2u\Vert +|u_-|\le C(\Vert v\Vert +h\vert v_+ \vert ).
}
Notice here that $(h^{-1/2}\Lambda )^NR_-={\cal O}(1):{\bf C}^{N_0}\to L^2$,
for every $N\in {\bf R}$ and similarly for $R_+$, $R_\pm^\delta $. From this,
 it follows that $R_-^\delta=R_-^{\delta,\epsilon } $ depends weakly on $\epsilon $ in the sense
that if $\widetilde{\epsilon }=\widetilde{A}h$, $\widetilde{A}\ge A$, then
$$
\Vert R_-^{\delta ,\widetilde{\epsilon }}-R_-^{\delta ,\epsilon }\Vert
_{{\cal L}({\bf C}^{N_0},L^2)}={\cal O}({1\over A}),
$$
and similarly for $R_+^\delta $.

\par We shall derive an a priori estimate for the problem
\ekv{gr.7}
{
(P_\delta -z)u+R_-^\delta u_-=v,\ R_+^\delta u=v_+,
}
when $u\in {\cal S}$. Let $\chi \in C_0^\infty (B(0,2))$ be equal to one
on $B(0,1)$, and put $\chi_{\sqrt{\epsilon }}(x,\xi )=\chi (\epsilon ^{-1/2}
(x,\xi ))$. We use the same notation for the corresponding $h$-quantization.
We may assume that $\widetilde{P}_\delta -P_\delta =\widetilde{{\cal O}}
(\epsilon (h/\epsilon )^\infty )$ on ${\rm supp\,}(\chi_{\sqrt{\epsilon }})$.
From the first equation in \no{gr.7}, we get
\ekv{gr.8}
{
\cases{\displaystyle
(\widetilde{P}_\delta -z)\chi_{\sqrt{\epsilon}}u+R_-^\delta u_-=
\chi_{\sqrt{\epsilon}}v+\chi_{\sqrt{\epsilon}}(\widetilde{P}_\delta -P_\delta )u
+[\widetilde{P}_\delta ,\chi_{\sqrt{\epsilon}}]u+
(1-\chi_{\sqrt{\epsilon}})R_-u_-,
\cr\displaystyle
R_+^\delta \chi_{\sqrt{\epsilon}}u
=v_+-R_+^\delta (1-\chi_{\sqrt{\epsilon}})u.
}
}
Here
\begin{eqnarray*}
\chi_{\sqrt{\epsilon}}\# (\widetilde{P}_\delta -P_\delta )&=&
\widetilde{{\cal O}}(h{h\over \epsilon })\\
{[} \widetilde{P}_\delta ,\chi_{\sqrt{\epsilon}} {]}&=&
\widetilde{{\cal O}}(h{h\over \epsilon })
\end{eqnarray*}
\ekv{gr.9}
{
\Vert (1-\chi_{\sqrt{\epsilon}})R_-^\delta u_-
\Vert= \Vert (1-\chi_{\sqrt{\epsilon}})(h^{-{1\over 2}}\Lambda )^{-N}
(h^{-{1\over 2}}\Lambda )^{N}R_-^\delta u_-\Vert={\cal O}(({h\over \epsilon} )^
{N\over 2})|u_-|,
}
\ekv{gr.10}
{
h|R_+^\delta (1-\chi_{\sqrt{\epsilon}})u|=h|R_+^\delta (h^{-{1\over 2}}
\Lambda )^N (h^{-{1\over 2}}
\Lambda )^{-N}(1-\chi_{\sqrt{\epsilon}})u|={\cal O}(({h\over\epsilon })^
{{N\over2}}h\Vert u\Vert ).
}
Thus, applying the a priori estimate \no{gr.6} to \no{gr.8}, we get
\ekv{gr.11}
{
\Vert \Lambda ^2\chi_{\sqrt{\epsilon}}u\Vert+\vert u_- \vert\le
C(\Vert \chi_{\sqrt{\epsilon}}v\Vert +{\cal O}(h{h\over \epsilon })\Vert u\Vert
+{\cal O}({h\over \epsilon })\vert u_- \vert +h|v_+|).
}

\par We next look for an a priori estimate for $(1-\chi_{\sqrt{\epsilon}})u$.
Apply $1-\chi_{\sqrt{\epsilon}}$ to \no{gr.7}:
\ekv{gr.12}
{
(P_\delta -z)(1-\chi_{\sqrt{\epsilon}})u=(1-\chi_{\sqrt{\epsilon}})v+
[P_\delta, \chi_{\sqrt{\epsilon}}]u-(1-\chi_{\sqrt{\epsilon}})R_-u_-.
}
As before,
$$[P_\delta ,\chi_{\sqrt{\epsilon}}]=[\widetilde{P}_\delta ,
\chi_{\sqrt{\epsilon}}]+[P_\delta -\widetilde{P}_\delta ,\chi_{\sqrt{\epsilon}}]
={\cal O}(h{h\over \epsilon }):\, L^2\to L^2,
$$
and using also \no{gr.9}, we get from \no{ap.11}, \no{gr.12} that
\ekv{gr.13}
{
{\epsilon \over C}\Vert (1-\chi_{\sqrt{\epsilon}})u\Vert
\le \Vert (1-\chi_{\sqrt{\epsilon}})v\Vert +\Vert K_\epsilon
(1-\chi_{\sqrt{\epsilon}})u\Vert +{\cal O}({h\over \epsilon }h)
\Vert u\Vert +{\cal O}({h\over \epsilon })|u_-|,
}
where we can take $K_\epsilon (x,\xi )=\epsilon \chi ({3(x,\xi )\over
\sqrt{\epsilon }})$. Since
$$
K_\epsilon \# (1-\chi_{\sqrt{\epsilon}})={\cal O}(\epsilon
({h\over \epsilon })^N):
\, L^2\to L^2,
$$
we get
\ekv{gr.14}
{
{\epsilon \over C} \Vert (1-\chi_{\sqrt{\epsilon}})u\Vert
\le \Vert (1-\chi_{\sqrt{\epsilon}})v\Vert+{\cal O}(h{h\over \epsilon })
\Vert u\Vert +{\cal O}({h\over \epsilon })|u_-|.
}

\par Let $L_\epsilon =\widetilde{{\cal O}}(\min (h+{\rm dist\,}
(x,\xi ;0,0)^2,\epsilon ))$ be an elliptic symbol in the class defined
by the right hand side, so that
$$
\Vert \Lambda ^2\chi_{\sqrt{\epsilon}}u\Vert+\epsilon \Vert
(1-\chi_{\sqrt{\epsilon}})u\Vert\backsim \Vert L_\epsilon u\Vert .
$$
Then summing \no{gr.11}, \no{gr.14}, we can absorb the various remainder
terms to the right, and obtain
\ekv{gr.15}
{
{1\over C}(\Vert L_\epsilon u\Vert+|u_-|)\le \Vert v\Vert +h\vert v_+ \vert,
}
when \no{gr.7} holds.

\par Now recall that $P_\delta $, $R_+^\delta $, $R_-^\delta $ have been
defined from $P$, $R_+$, $R_-$ by conjugation with $A_\delta $, and use
also that
$$
{1\over C}h\Vert u\Vert \le \Vert L_\epsilon u\Vert\le CAh\Vert u\Vert,
$$
to see that if $u,v\in {\cal S}(M)$ and
\ekv{gr.16}
{
\cases{
(P-z)u+R_-u_-=v\cr R_+u=v_+
}\ ,\ \ z\in D(0,Bh),
}
then
\ekv{gr.17}
{
h\Vert u\Vert+\vert u_- \vert \le C(\Vert v\Vert +h\vert v_+ \vert ).
}

\par From the discussion after the proof of Corollary \ref{macr2}, we conclude
that
\ekv{gr.18}
{
{\cal P}(z)=\pmatrix{P-z &R_-\cr R_+ &0}:\ {\cal D}(P)\times {\bf C}^{N_0}
\to L^2\times {\bf C}^{N_0}
}
is bijective with a bounded inverse
\ekv{gr.19}
{
{\cal E}=\pmatrix{E(z) &E_+(z)\cr E_-(z) &E_{-+}(z)}:\ L^2\times {\bf C}^{N_0}
\to {\cal D}(P)\times {\bf C}^{N_0},
}
for $z\in D(0,Bh)$, and \no{gr.17} shows that
\ekv{gr.20}
{
\Vert E(z)\Vert ={\cal O}({1\over h}), \Vert E_-(z)\Vert ={\cal O}(1),
\Vert E_+(z)\Vert ={\cal O}(1), \Vert E_{-+}(z)\Vert ={\cal O}(h).
}

\par In Section 11 of \cite{HeSjSt} the authors studied the action of
${\cal P}(z)$ on spaces of functions of the form
$(a(x;h)e^{i\Phi _0(x)/h},u_-)$,
where $a$ is a symbol, and deduced that $E_{-+}(z;h)$ has an asymptotic
expansion in half powers of $h$ with a certain additional structure.
From that was obtained the asymptotic expansion of the zeros of
$\det E_{-+}$, i.e. of the eigenvalues of $P$ in $D(0,Bh)$. That
discussion goes through without any changes in the present situation, so
we get the asymptotics for the eigenvalues in any disc $D(0,Bh)$, when $h\to 0$.

\begin{theo}\label{as1}
We make the assumptions {\rm \no{ap.1}}--{\rm \no{ap.7}}, {\rm \no{ny.1}}, {\rm \no{ny.17}},
{\rm \no{ny.18}}, and recall the definition of ${\cal C}$ in {\rm \no{ny.2}}.
Let $B>0$. Then there exists $h_0>0$ such that
for $0<h\le h_0$, the spectrum of $P$ in $D(0,Bh)$ is discrete and the
\ev{}s are of the form
\ekv{as.8}
{
\lambda _{j,k}(h)\sim h(\mu _{j,k}+h^{1/N_{j,k}}\mu _{j,k,1}+h^{2/N_{j,k}}\mu
_{j,k,2}+...),
}
where the $\mu _{j,k}$ are all the numbers in $D(0,B)$ of the form
\ekv{as.9}
{\mu _{j,k}={1\over i}\sum_{\ell =1}^n \nu _{j,k,\ell}\lambda _{j,\ell }
+{1\over 2}{\widetilde{\rm tr}\,}(p,\rho _j), \hbox{ with }\nu _{j,k,\ell}\in{\bf N},
}
for some $j\in\{ 1,...,N\}$, $N=\# {\cal C}$. (Possibly after
changing $B$, we may assume that $|\mu _{j,k}|\ne B$, $\forall j,k$.) Recall
here that  $\pm\lambda _j$ are the eigenvalues of $F_p$.
This description also
takes into account the multiplicities in the natural way. If the
\coef{}s $\nu_{j,k,\ell} $ in {\rm \no{as.9}} are unique, then $N_{j,k}=1$ and we
have only integer powers of $h$ in the \asy{} expansion {\rm \no{as.8}}.
\end{theo}

\begin{theo}\label{as2}
We make the same assumptions as in Theorem {\rm \ref{as1}}. For every $B$, $C>0$
there is a constant $D>0$ such that
\ekv{as.10}
{
\Vert (z-P)\inv\Vert \le {D\over h},\hbox{ for }z\in D(0,Bh)\hbox{ with
}{\rm dist\,}(z,\sigma (P))\ge {h\over C}.
}
\end{theo}
The last result follows from the formula
$$(z-P)\inv =-E(z)+E_+(z)E_{-+}(z)\inv E_-(z),$$
\no{gr.20} and the fact that $\Vert E_{-+}\inv (z)\Vert ={\cal O}(h\inv )$ when
\no{as.10} holds.

\par Still with $j=j_0$ fixed, let
\ekv{as.12}
{
\mu ={1\over i}\sum_{\ell =1}^n\nu _\ell \lambda _\ell +{1\over
2}{\widetilde{\rm tr}\,}(p,j_0),\ \nu _\ell\in {\bf N}
}
be a value as in \no{as.9} and assume that $\mu $ is simple in the
sense that $(\nu _1,...,\nu _n)\in{\bf N}^n$ is uniquely determined by
$\mu $. In
particular, every $\lambda _\ell$ for which $\nu _\ell\ne 0$ is a simple
\ev{} of $F_p$. Then as in \cite{HeSj} (see also Chapter 3 in \cite{DiSj})
we can construct
\ekv{as.13}
{
\lambda (h)\sim h(\mu +h\mu _1+h^2\mu _2+...)
}
with uniquely determined \coef{}s $\mu _1$, $\mu _2$, ... and
\ekv{as.14}
{
a(x;h)\sim a_0(x)+ha_1(x)+...\hbox{ in }C^\infty ({\rm neigh\,}(x_{j_0})),
}
where $a_j(x)={\cal O}(\vert x-x_{j_0}\vert ^{(m-2j)_+})$, $m=\sum \nu _\ell$
and $a_0$
has a non-vanishing Taylor polynomial of order $m$, such that
\ekv{as.15}
{
(P-\lambda (h))(a(x;h)e^{-\phi _+(x)/h})={\cal O}(h^\infty )e^{-\phi _+(x)/h}
}
in a \neigh{} of $x_{j_0}$. Actually any \neigh{} $\Omega \subset\subset
{\bf R}^n$ will do, provided that
\smallskip
\par\noindent 1) $\phi _+$ is well-defined in a \neigh{} of
$\overline{\Omega }$.\smallskip
\par\noindent 2) ${{H_q}_\vert}_{\Lambda _+}\ne 0$ on $\overline{\Omega
}\sm \{ x_{j_0}\}$. \smallskip
\par\noindent 3) $\Omega $ is star-shaped \wrt{} the point $x_{j_0}$
and the integral curves of the  \vf{}
$\nu _+:=\left(\pi_x\right)_*({{H_q}_\vert}_{\Lambda _+})$, where $\pi_x((x,\xi))=x$.\smallskip

\par We also know that $\lambda (h)$ is equal ${\rm mod\,}{\cal
O}(h^\infty )$ to the corresponding value in \no{as.8}.

\par As in \cite{HeSj} we notice that if $\gamma \subset D(0,B)$ is a
closed $h$-\indep{} contour avoiding all the values $\mu _{j,k}$ in
\no{as.9}, and
\ekv{as.16}
{
\pi _{h\gamma} ={1\over 2\pi i}\int_{h\gamma }(z-P)\inv dz
}
the corresponding spectral projection, then, using also Theorem \ref{as2},
\ekv{as.17}
{
\Vert \pi _{h\gamma }(\chi ae^{-\phi _+/h})-\chi ae^{-\phi _+/h}\Vert
_{L^2}={\cal O}(h^\infty )
}
if $\chi \in C_0^\infty (\Omega )$ is equal to one near $x_{j_0}$. It
follows that $\chi ae^{-\phi _+/h}$ is a linear combination of
generalized \ef{}s of $P$ with \ev{}s inside $h\gamma $ up to an error
${\cal O}(h^\infty )$ in $L^2$-norm.

\section{Exponentially weighted estimates}\label{exp}
\setcounter{equation}{0}

\par In this section we keep the general assumptions on $P$ and assume
for simplicity that ${\cal C}$ is reduced to a single point:
\ekv{exp.1}
{{\cal C}=\{ (0,0)\} .}
If $\psi \in C^\infty ({\rm neigh\,}(0,M);{\bf R})$, we have
\ekv{exp.2}
{
e^{\psi /h}\circ P\circ e^{-\psi /h}=P_\psi ,
}
with the symbol (cf \no{ap.12})
\ekv{exp.3}
{p_\psi (x,\xi )=p_2(x,\xi )-q(x,\psi '(x))+i(\partial _\xi q)(x,\psi
'(x))\cdot \xi ,}
where we recall that
\ekv{exp.4}
{q(x,\xi )=p_2(x,\xi )+p_1(x,\xi )-p_0(x).}
Notice that $\xi \mapsto q(x,\xi )$ is a convex \fu{} for every $x$.

\par Let $\phi =\phi _+(x)\in C^\infty ({\rm neigh\,}(0;{\bf R}))$ be the
\fu{}
introduced in Section \ref{as} so that $\Lambda _+=\Lambda _\phi $ is the
stable outgoing \mfld{} through $(0,0)$ for the $H_q$-flow. Recall that by
Proposition \ref{as4}
\ekv{exp.5}
{
\phi ''(0)>0.
}
We have the eikonal \e{}
\ekv{exp.6}
{q(x,\phi '(x))=0,}
so $p_\phi (x,\xi )=p_2(x,\xi )+i(\partial _\xi q)(x,\phi '(x))\cdot \xi $.
The \vf{} $(\partial_\xi q)(x,\phi '(x))\cdot \partial _x $ is the
$x$-space projection of ${{H_q}_\vert}_{\Lambda _\phi }$, so its
linearization at $x=0$ has all its \ev{}s with real part $>0$.
Consequently (as we shall see in more detail in the proof of Lemma
\ref{ss1} below), there exists $G\in C^\infty ({\rm neigh\,}(0,M);{\bf R})$
such that
\ekv{exp.7}
{
(\partial _\xi q)(x,\phi '(x))\cdot \partial _xG\backsim x^2,\ G(x)\backsim x^2.
}

\par Let $\Omega _G(r)=\{ x\in{\rm neigh\,}(0);\, G(x)\le r\}$ for $0<r\ll
1$. Outside the set $\Omega _G(C_0\epsilon )$, we put
\ekv{exp.8}
{\widehat{\psi }=\phi -\epsilon g(G),}
for a suitable \fu{} $g$ and for $0<\epsilon \ll 1$. (Eventually $\epsilon $
will be proportional to $h$.) Using \no{exp.6}, we get
\ekv{exp.9}
{
q(x,\widehat{\psi }'(x))=-\epsilon g'(G)(\partial _\xi q)(x,\phi
'(x))\cdot \partial _xG(x)+\epsilon ^2g'(G)^2p_2(x,G'(x)).
}
Choose $g(G)=\ln G$ for $G\ge C_0 \epsilon $, so that $g'(G)=1/G$. Then
\ekv{exp.10}
{
q(x,\widehat{\psi }'(x))=-\epsilon {1 \over G(x)}(\partial _\xi q)(x,\phi
'(x))\cdot \partial _xG+\epsilon {\epsilon \over G(x)}{p_2(x,G'(x))\over
G(x)}.
}
Here
$$
{1\over G}(\partial _\xi q)(x,\phi ')\cdot \partial _xG\backsim 1,\quad
{p_2(x,G')\over
G}={\cal O}(1).
$$
We conclude that
\ekv{exp.11}
{
q(x,\widehat{\psi }'(x))\le -{\epsilon \over C_0},\ x\in {\rm
neigh\,}(0,M)\sm \Omega _G(C_0\epsilon ),}
if $C_0>0$ is large enough.

\par Outside a small fixed \neigh{} of 0 we want to flatten out the
weight. Let $f_\delta (t)=\delta f({t\over \delta })$ be the \fu{}
introduced in
Section \ref{ap}. For some small and fixed $\delta _0>0$, we
put
\ekv{exp.12}
{
\psi =f_{\delta _0}(\widehat{\psi })=f_{\delta _0}(\phi -\epsilon g(G))
}
which is also well-defined as the constant $3\delta _0/2$ for large $x$.
From \no{exp.11}, the fact that $q(x,0)\le 0$ and the convexity of $q$, we
get
\ekv{exp.13}
{
q(x,\psi '(x))\le -{\epsilon \over C}f'_{\delta _0}(\widehat{\psi }),
}
where we keep in mind that $0\le f_{\delta _0}'\le 1$.

\par We extend the definition of $\psi $ to a full
\neigh{} of $x=0$, by putting
\ekv{exp.14}
{
g(G)=\ln (C_0\epsilon )+{1\over C_0\epsilon }(G-C_0\epsilon ),\hbox{ for }
0\le G\le C_0\epsilon .
}
Then in $\Omega _G(C_0\epsilon )$, we have
$\widehat{\psi }'=\psi '=\phi '-G'/C_0$,
so $q(x,\psi ')= {\cal O}(\epsilon )$ by \no{exp.9}.

\par In the exterior region where $\psi '=f_{\delta _0}'(\widehat{\psi
})\widehat{\psi }'$ is small, the conjugated \op{}
$\widetilde{P}=P_\psi $ is close to the unperturbed \op{} $P$ and we can
apply the method of Section \ref{ap}. Write
\ekv{exp.15}
{
\widetilde{p}=p_2+i\widetilde{p}_1+\widetilde{p}_0
}
for the symbol of $\widetilde{P}$, so that by \no{exp.3}, we have
\ekv{exp.16}
{
\widetilde{p}_1=(\partial _\xi q)(x,\psi '(x))\cdot \xi ,\
\widetilde{p}_0=-q(x,\psi '),
}
to be compared with $p_1=(\partial _\xi q(x,0))\cdot \xi $, $p_0=-q(x,0)$.
Let $\widetilde{\nu }(x,\partial _x)=\sum \widetilde{c}_j(x)\partial
_{x_j}$, where $\widetilde{p}_1=\sum \widetilde{c}_j(x)\xi _j$ and
define $\widetilde{\psi }_\epsilon $ as in \no{ap.14}, now for
$\widetilde{P}$ instead of $P$ and with an additional cut-off:
\ekv{exp.17}
{
\widetilde{\psi }_\epsilon =(1-\chi (x))\int k({t\over T_0})f_\epsilon
\circ \widetilde{p}_0\circ \exp (t\widetilde{\nu })dt.}
Here $\chi \in C_0^\infty (M)$ is equal to 1 near $x=0$ and has its
support in a small \neigh{} of that point.

\par Then for $\widetilde{P}_{\delta \widetilde{\psi }_\epsilon
}=e^{\delta \widetilde{\psi }_\epsilon /h}\circ \widetilde{P}\circ
e^{-\delta \widetilde{\psi }_\epsilon /h}$ we get (cf \no{ap.22}),
\ekv{exp.18}
{
\widetilde{p}_{\delta \widetilde{\psi }_\epsilon }\ge (1-\delta
)\widetilde{p}_0+\delta \langle f_\epsilon \circ \widetilde{p}_0\rangle
_{T_0}-{\cal O}(\delta ^2\epsilon )\hbox{ outside }{\rm supp\,}\chi ,
}
where the time average in the second term to the right is taken along the
trajectories of $\widetilde{\nu }$. In the region where $\chi =1$, we get
$\widetilde{P}_{\delta \widetilde{\psi }_\epsilon }=\widetilde{P}$ of
course, and in the intermediate region, ${\rm supp\,}(\nabla \chi )$, we
have
\ekv{exp.19}
{
\widetilde{p}_{0,\delta \widetilde{\psi }_\epsilon }=\widetilde{p}_0+{\cal
O}(\delta \epsilon ).
}

\par If $\delta _1>0$ is small enough, we know that
$$
f_{\delta _0}'(\widehat{\psi }(x))\le \delta _1\Rightarrow \langle
f_\epsilon \circ \widetilde{p}_0\rangle _{T_0}(x)\ge {\epsilon \over C_1},
$$
for some constant $C_1>0$.
In fact, the $\nu $ and $\widetilde{\nu }$ trajectories through a given
point with $f_{\delta _0}'(\widehat{\psi }(x))\le \delta $ stay close for
some fixed time $>0$ , so the conditions \no{ny.19}, \no{ny.20} imply that the
$\widetilde{\nu }$-trajectory will encounter points with $p_0\ge {1/{\rm
Const}}$ during a non-trivial interval of time.

\par Then
\begin{itemize}
\item In the region where $f_{\delta _0}'(\widehat{\psi }(x))\le \delta
_1$, we have $\widetilde{p}_{0,\delta \widetilde{\psi }_\epsilon }\ge
(1-\delta )\widetilde{p}_0+{\delta \epsilon \over C_1}$ (and we recall
that $\widetilde{p}_0\ge 0$).
\item In the region where $f_{\delta _0}'(\widehat{\psi }(x))>\delta _1$,
and $G(x)\ge C_0\epsilon $, we have $\widetilde{p}_{0,\delta \widetilde{\psi
}_\epsilon }=\widetilde{p}_0+{\cal O}(\delta \epsilon )\ge {\epsilon \delta
_1\over C}+{\cal O}(\delta \epsilon )$, by \no{exp.13}.
\item In $\Omega _G(C_0\epsilon )$, we have $\widetilde{p}_{0,\delta
\widetilde{\psi }_\epsilon }=\widetilde{p}_0= {\cal O}(\epsilon )$.
\end{itemize}

\par Choosing first $\delta _1>0$ small enough, then $\delta >0$ small
enough, we conclude that
\ekv{exp.20}
{
\widetilde{p}_{0,\delta \widetilde{\psi }_\epsilon }\cases{\displaystyle
\ge {\epsilon
\over C}\hbox{ outside }\Omega _G(C_0\epsilon ),\cr
\displaystyle
={\cal O}(\epsilon )
\hbox{ in }\Omega
_G(C_0\epsilon ).}}

\par Now $\widetilde{P}_{\delta \widetilde{\psi }_\epsilon }=e^{\psi
_\epsilon /h}\circ P\circ e^{-\psi _\epsilon /h}=P_{\psi _\epsilon }$, where
\ekv{exp.21}
{
\psi _\epsilon =\psi +\delta \widetilde{\psi }_\epsilon ,
}
and where we recall that $\psi $
also depends on $\epsilon $. Combining Lemma \ref{ap1} for $P_{\psi
_\epsilon }$ with \no{exp.20}, we get
\eekv{exp.22}
{\int_{M\sm \Omega_ G(C_0\epsilon )}({\epsilon \over C}-\Re z)\vert u\vert
^2dx+\int_M (P_2u)\overline{u}dx\le }
{\Re ((P_{\psi _\epsilon }-z)u\vert u)+\int_{\Omega _G(C_0\epsilon )}(
{\cal O}(\epsilon )+\Re
z ) \vert u\vert ^2dx.}
If $\Re z={\cal O}(h)$, we choose $\epsilon =Ah$ with $A$ so large that
${\epsilon \over C}-\Re z\ge h$ and deduce that
\ekv{exp.23}
{
h\Vert u\Vert ^2+(P_2u\vert u)\le \Re ((P_{\psi _\epsilon }-z)u\vert
u)+{\cal O}(h)\Vert u\Vert _{\Omega _G(C_0\epsilon )}^2,
}
leading to the a priori estimate
\ekv{exp.24}
{
h\Vert u\Vert \le \Vert (P_{\psi _\epsilon} -z)u\Vert +{\cal O}(h)\Vert
u\Vert _{\Omega _G(C_0\epsilon )}.}
Re-injecting this estimate in \no{exp.23}, we get
\ekv{exp.25}
{
h^2\Vert u\Vert ^2+h\Vert B^{1\over 2}hDu\Vert ^2\le (1+\alpha )\Vert
(P_{\psi _\epsilon }-z)u\Vert ^2+{\cal O}_\alpha (h^2)\Vert u\Vert
^2_{\Omega _G(C_0\epsilon ),}
}
for every fixed $\alpha >0$.

\par Here $(P_{\psi _\epsilon }-z)u=e^{\psi _\epsilon /h}(P-z)v$,
$u=e^{\psi _\epsilon /h}v$, so from \no{exp.25} we get the a priori
estimate for the original \op{}
\ekv{exp.26}
{
h\Vert e^{\psi _\epsilon /h}v\Vert +h^{1\over 2}\Vert B^{1\over
2}hD(e^{\psi _\epsilon /h}v)\Vert \le {\cal O}(1)\Vert e^{\psi _\epsilon
/h}(P-z)v\Vert +{\cal O}(h)\Vert e^{\psi _\epsilon /h}v\Vert _{\Omega
_G(C_0\epsilon )},
}
\ufly{}, for $\vert \Re z\vert \le Ch$ provided that $\epsilon =Ah$ for
$A$ large enough depending on $C$.

\par Now let $\lambda (h)=\lambda _{1,k}(h)$ be an \ev{} of $P$ as in
\no{as.8}, \no{as.13} and assume that $\mu $ is given by \no{as.12} and is
simple, as explained after that \e{}. Then $\lambda (h)$ is a simple \ev{}
of $P$ and is the only \ev{} in some disc $D(\lambda (h),h/C_0)$. Let
$u_{\rm WKB}(x;h)$ be the approximate solution given in \no{as.14},
\no{as.15} and let $u=\pi _{h\gamma }(\chi u_{\rm BKW})$ be the
corresponding exact \ef{}, where $\gamma =\partial D(\mu ,{1\over 2C_0})$.
\begin{theo}\label{exp1}
{\rm a)} Outside any $h$-\indep{} \neigh{} of 0, we have
$$
u,\, B^{1\over 2}hDu={\cal O}(e^{-1/(Ch)})
$$
in $L^2$-norm.
\smallskip
\par\noindent {\rm b)} There exists a \neigh{} $\Omega $ of $0$, where
\eekv{exp.27}
{
&&u(x;h)=(a+r)e^{-\phi _+(x)/h},}{&&\Vert r\Vert _{L^2(\Omega )},\, \Vert
B^{1\over 2}hDr\Vert _{L^2(\Omega )}={\cal O}(h^\infty ).
}
\end{theo}
\begin{proof}
Apply \no{exp.26} with $v=u$, $z=\lambda (h)$, $\epsilon =Ah$, $A\gg 1$, to
get
$$
h\Vert e^{\psi _\epsilon /h}u\Vert +h^{1\over 2}\Vert B^{1\over
2}hD(e^{\psi _\epsilon /h}u)\Vert \le {\cal O}(h^{-N}),\ N=N(A),
$$
where we also used that $e^{\psi _\epsilon /h}={\cal O}(h^{-N})$ in $\Omega
_G(C_0\epsilon )$. Here, $\psi _\epsilon =f_{\delta _0}(\phi -\epsilon
g(G))+{\cal O}(\delta \epsilon )$ is larger than a positive constant
outside any fixed \neigh{} of $0$, so $u={\cal O}(e^{-1/(Ch)})$ in
$L^2$-norm there. Moreover, since $\Vert e^{\psi _\epsilon /h}u\Vert ={\cal
O}(h^{-N})$, we have
$$
{\cal O}(h^{-N})\ge h^{1\over 2}\Vert B^{1\over 2}hD(e^{\psi _\epsilon
/h}u)\Vert \ge h^{1\over 2}\Vert e^{\psi _\epsilon /h}B^{1\over
2}hDu\Vert -{\cal O}(h^{1\over 2})\Vert \vert \nabla \psi _\epsilon \vert
e^{\psi _\epsilon /h}u\Vert .
$$
Here $\nabla \psi _\epsilon ={\cal O}(1)$ (as we shall see more in detail
below), so $\Vert e^{\psi _\epsilon /h}B^{1\over 2}hDu\Vert ={\cal
O}(h^{-N-{1\over 2}})$, so $B^{1\over 2}hDu={\cal O}(e^{-{1\over Ch}})$ in
$L^2$-norm away from any given fixed \neigh{} of $0$. The proof of a) is
complete.

\par To prove b), we apply \no{exp.26} to $v=u-\chi u_{\rm WKB}$,
$z=\lambda (h)$ with $\epsilon =Ah$, $A\gg 1$.  Since $\Vert u-\chi u_{\rm
WKB}\Vert ={\cal O}(h^\infty )$ by \no{as.17},
and $\Vert e^{\psi _\epsilon /h}(P-z)\chi u_{{\rm WKB\,}}\Vert
={\cal O}(h^\infty )$ if we arrange so that $\psi _\epsilon <\phi $ on
${\rm supp\,}\nabla \chi$,
we conclude that $\Vert
e^{\psi _\epsilon /h}(u-\chi u_{\rm WKB})\Vert _{\Omega _G(C_0\epsilon
)}={\cal O}(h^\infty )$ and hence
\ekv{exp.28}
{
h\Vert e^{\psi _\epsilon /h}(u-\chi u_{\rm WKB})\Vert +h^{1\over 2}\Vert
B^{1\over 2}hD(e^{\psi _\epsilon /h}(u-\chi u_{\rm WKB})\Vert \le {\cal
O}(h^\infty ).
}

\par In a small \neigh{} of $0$, we have
$$
e^{\psi _\epsilon /h}=e^{{\phi \over h}-{\epsilon g(G)\over
h}}=e^{-Ag(G)}e^{\phi \over h},
$$
where by \no{exp.14}
\smallskip
\par\noindent\hskip 2truecm $e^{-Ag(G)}={\cal O}(1)\epsilon ^{-A}$
for $G(x)\le C_0\epsilon $,
\smallskip
\par\noindent\hskip 2truecm $e^{-Ag(G)}=G(x)^{-A}$ for $G(x)\ge C_0\epsilon $.
\smallskip
\par\noindent It then follows from \no{exp.28}, that $u=(a+r)e^{-\phi
_+/h}$ in a \neigh{} $\Omega $ of 0, with $\Vert r\Vert _{L^2(\Omega
)}={\cal O}(h^\infty )$.

\par To get the corresponding bound on $B^{1\over 2}hDr$, we just have to
proceed as in the proof of a) and use that $\nabla \psi _\epsilon
=\nabla \phi -\epsilon g'(G)\nabla G$, where
$$
g'(G)\nabla G={\cal O}(\epsilon ^{-{1\over 2}})\hbox{ for }G(x)<C_0\epsilon ,\ g'(G)\nabla G={\nabla G\over
G}={\cal O}(\epsilon ^{-{1\over 2}})\hbox{ for }G(x)\ge C_0\epsilon .
$$
Thus $\nabla \psi _\epsilon =\nabla \phi +{\cal O}(\sqrt{\epsilon })={\cal
O}(1)$ and we conclude that $\Vert B^{1\over 2}hDr\Vert _{L^2(\Omega )}={\cal
O}(h^\infty )$.
\end{proof}
\begin{remark}\label{exp2} \rm
If we drop the assumption \no{exp.1} and allow $N-1$ more points $\rho
_2,...,\rho _N$ in ${\cal C}$, then Theorem \ref{exp1} is still valid,
provided that all $\mu _{j,k}$ in \ref{as.8} with $j\ge 2$ are
different from the value $\mu $, associated to $\rho _1=(0,0)$.
\end{remark}
\section{Supersymmetric approach}\label{ss}
\setcounter{equation}0

\par The Witten approach has been \indep{}ly extended to the case of
non-elliptic \op{}s like the \KFP{} \op{} in \cite{TaTaKu} (in
supersymmetric language) and in \cite{Bi} (in terms of differential
forms). See also \cite{Le}.

\par We start by a quick review of that in the semiclassical case,
then we establish some basic facts about the principal and
subprincipal symbols, especially at the critical points of the given
weight \fu{.}

\subsection{Generalities}\label{ssa}
Let
\ekv{ss.1}
{
A(x):T_x^*M\to T_xM,\ x\in M,
}
be an invertible map depending smoothly on $x\in M$.
Then we have the real
\nondeg{} bilinear form
\ekv{ss.2}
{
\langle u\vert v\rangle _{A(x)}=\langle \wedge ^kA(x)(u)\vert v\rangle ,\
u,v\in \wedge ^kT_x^*M.
}
If $a:\wedge ^kT_x^*M\to \wedge ^\ell T_x^*M$ is a linear map, we define
the "adjoint" $a^{A,*}:\wedge ^\ell T_x^*M\to \wedge ^kT_x^*M$ by
\ekv{ss.3}
{
\langle au\vert v\rangle _{A(x)}=\langle u\vert a^{A,*}v\rangle _{A(x)}
}
(In the complexified case, we use the sesquilinear scalar product $(u\vert
v)_A=\langle u\vert \overline{v}\rangle _A$ and define $a^{A,*}$ the same
way.)

\par If $\omega $
is a one form and $u$ and $v$ are $k-1$ and $k$ forms respectively, we get
at $x$:
\ekv{ss.4}
{
\langle \omega \wedge u\vert v\rangle _A=\langle A\omega \wedge (\wedge
^{k-1}A)u\vert v\rangle =\langle \wedge ^{k-1}Au\vert (A\omega )^\rfloor
v\rangle =\langle u\vert (A\omega )^\rfloor v\rangle _A,
}
so
\ekv{ss.5}
{ (\omega ^\wedge )^{A,*}=(A\omega)^\rfloor,  } where ${}^\rfloor$ denotes
the usual \op{} of contraction. Let $\mu (dx)$ be a locally
finite measure with a smooth positive density. When $M={\bf
R}^n$, $\mu $
will be the Lebesgue measure. Sometimes we also use the symbol $\mu $ for
the corresponding density.
\par If $u,v$ are smooth $k$ forms with $\supp u\cap\supp v$ compact, we
define
$$
\langle u\vert v\rangle _A=\int \langle u(x)\vert v(x)\rangle _{A(x)}\mu
(dx),\ ( u\vert v) _A=\int (u(x)\vert v(x)) _{A(x)}\mu (dx) $$ and denote
by $a^{A,*}$ the formal adjoint of an \op{} $a:C_0^\infty (M;\wedge
^kT^*M)\to {\cal D}'(M;\wedge ^\ell T^*M)$.  If we fix some local
coordinates $x_1,..,x_n$ and write $\mu (dx)=\mu (x)dx$ (by slight abuse of
notation), we can consider
$$
\partial _{x_j}:C_0^\infty (M;\wedge ^kT^*M)\to
C_0^\infty (M;\wedge ^kT^*M),
$$
acting coefficient-wise, and a straightforward computation shows that
\ekv{ss.6}
{ (\partial _{x_j})^{A,*}=(\mu\, {\trans ( \wedge ^kA}))^{-1}(-\partial
_{x_j})\circ (\mu \, {\trans (\wedge ^kA})) =-\partial _{x_j}-{\partial _{x_j}\mu
\over \mu }-{\trans (\wedge ^kA)\inv }{\trans (\partial _{x_j}(\wedge
^kA)).} } We only retain that
\ekv{ss.7}
{
(h\partial _{x_j})^{A,*}= -h\partial _{x_j}+{\cal O}(h),
}
where ${\cal O}(h)$ stands for multiplication by a smooth matrix, which
is ${\cal O}(h)$ with all its derivatives, \ufly{} on all of $M$ when
$M={\bf R}^n$, and which is $=0$ when $A(x)$, $\mu (x)$ are constant.

\par Let $\phi \in C^\infty (M;{\bf R})$
and introduce the Witten (de Rham) complex
\ekv{ss.8}
{
d_\phi =e^{-\phi /h}\circ hd\circ e^{\phi /h}=hd+(d\phi )^\wedge :\,
C_0^\infty (M;\wedge ^kT^*M)\to C_0^\infty (M;\wedge ^{k+1}T^*M),
}
with $d_\phi ^2=0$.

\par In local coordinates (always the canonical ones when $M={\bf R}^n$) we have
\ekv{ss.9}
{d_\phi =\sum_1^n(h\partial _{x_j}+\partial _{x_j}\phi )\circ dx_j^\wedge ,}
where $h\partial _{x_j}+\partial _{x_j}\phi $ acts coefficient-wise and
commutes with $dx_j^\wedge $,
so
\ekv{ss.10}
{
d_\phi ^{A,*}=\sum_1^n (-h\partial _{x_k}+\partial _{x_k}\phi +{\cal
O}(h))\circ A(dx_k)^\rfloor ,
}
where from now on, ${\cal O}(h)$ and ${\cal O}(h^2)$ will
have the same meaning as after \no{ss.7}.
\par The corresponding Witten-Hodge Laplacian is given by
\ekv{ss.11}
{-\Delta _A=d_\phi ^{A,*}d_\phi +d_\phi d_\phi ^{A,*}.}
Since
$$
(d_\phi ^{A,*})^2=(d_\phi ^2)^{A,*}=0,
$$
we also have
$$
-\Delta _A=(d_\phi +d_\phi ^{A,*})^2,
$$
and $-\Delta _A$ conserves the degree of differential forms.

\par Choose local coordinates $x_1,...,x_n$ (to be the standard ones
when $M={\bf R}^n$)
and write
$$
A(dx_k)=\sum_j A_{j,k}(x)\partial _{x_j},\ A_{j,k}(x)=\langle A(dx_k)\vert
dx_j\rangle .
$$

Let $Z_j=h\partial _{x_j}+\partial _{x_j}\phi$, $Z_k^{A,*}=h\partial
_{x_k}^{A,*}+\partial _{x_k}\phi $. Notice that $[Z_j,dx_k^\wedge ]=0$,
so
$$
[Z_k^{A,*},A(dx_j)^\rfloor ]=0.
$$
Writing $d_\phi =\sum_1^nZ_j\circ
dx_j^\wedge $, $d_\phi ^{A,*}=\sum_1^nZ_k^{A,*}\circ A(dx_k)^\rfloor$, we
get
$$
-\Delta _A=\sum_{j,k}(Z_k^{A,*}\circ A(dx_k)^\rfloor dx_j^\wedge
Z_j+Z_jdx_j^\wedge A(dx_k)^\rfloor Z_k^{A,*}).
$$
Here, we use the general identity $\nu ^\wedge \mu ^\rfloor+\mu ^\rfloor
\nu ^\wedge =\langle \nu ,\mu \rangle 1$ on the first term in the
parenthesis to get
\begin{eqnarray*}
-\Delta _A&=&{\rm I}+{\rm II}+{\rm III},\\
{\rm I}&=&\sum_{j,k}Z_k^{A,*}A_{j,k}Z_j,\\
{\rm II}&=& -\sum_{j,k}Z_k^{A,*}dx_j^\wedge A(dx_k)^\rfloor Z_j\\
{\rm III}&=& \sum_{j,k}dx_j^\wedge Z_jZ_k^{A,*}A(dx_k)^\rfloor,
\end{eqnarray*}
where $A_{j,k}=\langle dx_j,A(dx_k)\rangle $. We have
$$
[Z_j,A(dx_k)^\rfloor ]=h((\partial _{x_j}A)(dx_k))^\rfloor.
$$
Using the identity $(U^{A,*})^{\trans A,*}=U$ (see Subsection
\ref{ssd}), we see that
$$
[Z_k^{A,*},dx_j^\wedge ]^{{\trans A},*}=[({\trans A(dx_j))^\rfloor},Z_k]=-h((\partial
_{x_k}{\trans A})(dx_j))^\rfloor,
$$
hence
$$
[Z_k^{A,*},dx_j^\wedge ]=-h({\trans A}\inv \partial _{x_k}{\trans
A}(dx_j))^\wedge .
$$
Using these commutator relations, we move $dx_j^\wedge $ to the left and
$A(dx_k)^\rfloor$ to the right in ${\rm II}$ and combining with
${\rm III}$, we finally obtain that
\eeekv{ss.11.5}
{
-\Delta _A&=&\sum_{j,k}Z_k^{A,*}A_{j,k}Z_j+\sum_{j,k}dx_j^\wedge
[Z_j,Z_k^{A,*}]A(dx_k)^\rfloor
}
{
&&+h\sum_{j,k} Z_k^{A,*}dx_j^\wedge ((\partial _{x_j}A)(dx_k))^\rfloor
+h\sum_{j,k}({\trans A}\inv (\partial _{x_k}{\trans A})dx_j)^\wedge
A(dx_k)^\rfloor Z_j
}
{&&
+h^2\sum_{j,k}({\trans A}\inv(\partial _{x_k}{\trans A})dx_j)^\wedge
((\partial _{x_j}A)dx_k)^\rfloor .}

\par Modulo ${\cal O}(h)(h\partial _x+\partial _x\phi )+{\cal
O}(h)(-h\partial _x+\partial _x\phi )+{\cal O}(h^2)$, we get
\eekv{ss.12}
{
-\Delta _A&\equiv& \sum_{j,k}(-h\partial _{x_k}+\partial _{x_k}\phi
)A_{j,k}(x)(h\partial _{x_j}+\partial _{x_j}\phi )
}{&&+\sum_{j,k}2h\partial _{x_j}\partial _{x_k}\phi \circ dx_j^\wedge
A(dx_k)^\rfloor ,
}
where the error terms vanish when $A(x)$ and $\mu (x)$
are constant (for the chosen coordinates).

Now write
\ekv{ss.13}
{
A(x)=B(x)+C(x),\ {\trans B}(x)=B(x),\ {\trans C}(x)=-C(x).
}
Then \no{ss.12} gives
\eeekv{ss.14}
{
-\Delta _A&\equiv&\sum_{j,k}(-h\partial _{x_k}+\partial _{x_k}\phi
)B_{j,k}(x)(h\partial _{x_j}+\partial _{x_j}\phi )
}
{
&&+\sum_{j,k}((\partial _{x_k}\phi )C_{j,k}h\partial _{x_j}+h\partial
_{x_j}\circ C_{j,k}\circ (\partial _{x_k}\phi ))-\sum_{j,k}h\partial
_{x_k}(C_{j,k})h\partial _{x_j}
}
{
&&+\sum_{j,k}2h\partial _{x_j}\partial _{x_k}\phi \circ dx_j^\wedge
A(dx_k)^\rfloor .
}
Again, this becomes an equality when $A,\mu $ are constant.
Note that the last term vanishes on $0$-forms, i.e. on scalar \fu{}s. To
recover the \KFP{} \op{} (cf \cite{TaTaKu}), replace $n$ by $2n$, put $M={\bf R}^{2n}_{x,y}$,
$$
A={1\over 2}\pmatrix{0 & 1\cr -1 &\gamma },
$$
and let $\mu =dxdy$ be the Lebesgue measure. Then \no{ss.14} is an
equality and we get for 0-forms:
\eeekv{ss.15}
{
-\Delta _A^{(0)}&=&{\gamma \over 2}\sum_{j=1}^n (-h\partial _{y_j}+\partial
_{y_j}\phi )(h\partial _{y_j}+\partial
_{y_j}\phi )+ }
{&&+{1\over 2}\sum_j [(\partial _{y_j}\phi )h\partial
_{x_j}-(\partial _{x_j}\phi )\, h\partial _{y_j}+h\partial _{x_j}\circ \partial
_{y_j}\phi -h\partial _{y_j}\circ \partial _{x_j}\phi ]
}
{&=&{\gamma \over 2}\sum_{j=1}^n (-h\partial _{y_j}+\partial
_{y_j}\phi )(h\partial _{y_j}+\partial
_{y_j}\phi )+hH_\phi ,
}
where
$$
H_\phi =\sum (\partial _{y_k}\phi \,  \partial _{x_k}-\partial _{x_k}\phi \,
\partial _{y_k})
$$
is the Hamilton field of $\phi $ \wrt{}  the standard symplectic form
$\sum dy_j\wedge dx_j$.

\par If we choose
\ekv{ss.16}
{
\phi (x,y)={1\over 2}y^2+V(x),
}
we get the  \KFP{} \op{}
\ekv{ss.17}
{-\Delta _A^{(0)}=y\cdot h\partial _x-V'(x)\cdot h\partial _y+{\gamma
\over 2}(-h\partial _y+y)\cdot (h\partial _y+y).
}

\subsection{The principal symbol of the Hodge Laplacian}\label{ssb}
\par The principal symbol of $-\Delta _A$ in the
sense of $h$-\dop{}s
is scalar and given by
\eekv{ss.18}
{
p(x,\xi )&=&\sum_{j,k}A_{j,k}(-i\xi _k+\partial _{x_k}\phi )
(i\xi _j+\partial _{x_j}\phi )
}
{&=&\sum_{j,k}B_{j,k}(\xi _j\xi _k+\partial _{x_j}\phi \,  \partial _{x_k}\phi
)+2i\sum_{j,k}C_{j,k}\partial _{x_k}\phi \,  \xi _j.}
The corresponding real symbol $q(x,\xi )=-p(x,i\xi )$ is given by
\eekv{ss.19}
{
q(x,\xi )&=&\sum_{j,k}A_{j,k}(\xi _k+\partial _{x_k}\phi )(\xi _j-\partial
_{x_j}\phi )
}
{&=&\sum_{j,k}B_{j,k}(\xi _j\xi _k-\partial _{x_j}\phi \,  \partial
_{x_k}\phi )+2\sum_{j,k}C_{j,k}\partial _{x_k}\phi \,  \xi _j.
}
It vanishes on the two Lagrangian \mfld{}s $\Lambda _{\pm \phi }$.

\par We define
\ekv{ss.20}
{
\nu _\pm={{H_q}_\vert}_{\Lambda _{\pm \phi }}.
}
Using $x_1,...,x_n$ as coordinates on $\Lambda _{\pm \phi }$, we get
\ekv{ss.21}
{
\nu _+=2\sum_{j,k}A_{j,k}\partial _{x_k}\phi \,  \partial _{x_j}=2A(x)(\phi
'(x))\cdot \partial _x
}
\ekv{ss.22}
{
\nu _-=-2\sum_{j,k}A_{j,k}\partial _{x_j}\phi \, \partial _{x_k}=
-2\,{\trans A(x)}(\phi
'(x))\cdot \partial _x .
}
(Even more radically, we could say that $\nu _+=2A(x)(d\phi (x))$, where
$A(x)$ is viewed as a map $T^*_xM\to T_xM$, and similarly for $\nu _-$.)

\par Let $x_0$ be a \nondeg{} critical point of $\phi $, so that $\Lambda
_\phi $ and $\Lambda _{-\phi }$ intersect transversally at $(x_0,0)$.
The spectrum of the
linearization $F_q$ of $H_q$ at $(x_0,0)$ is equal to the union of the
spectra of the linearizations
\ekv{ss.23}
{\nu _+^0=(2A(x_0)\phi ''(x_0)x)\cdot \partial _x\hbox{ and }
\nu _-^0=-(2{\trans A}(x_0)\phi ''(x_0)x)\cdot \partial _x
}
of $\nu _+$ and $\nu _-$ respectively at $x_0$. Thus we are interested in
the \ev{}s of the matrices $A\phi ''$, ${\trans A}\phi ''$, where we
write $A=A(x_0)$, and $\phi ''=\phi ''(x_0)$ for short. Here, we notice
that ${\trans A}\phi ''={\phi ''}\inv{\trans (A\phi '')\phi ''}$ has the
same \ev{}s as $A\phi ''$ and similarly $\phi ''A$, $\phi ''{\trans A}$ are
isospectral to $A\phi ''$. Thus
\eekv{ss.24}
{
&&\hbox{The \ev{}s of }F_q \hbox{ are given by }\pm 2\lambda _j,}
{&&\hbox{where }\lambda _1,...,\lambda _n\hbox{ are the \ev{}s of }A\phi ''.
}

\par From questions about hypoellipticity (see \cite{Sj}) we would like to
know when all the \ev{}s of $F_p$
avoid the real axis, or equivalently, when all the \ev{}s of  $F_q$ (the
linearization of $H_q$ at $(x_0,0)$) avoid the imaginary axis.

\par We assume from now on that
\ekv{ss.25}
{
B(x)\ge 0,\ x\in M.
}
Then, if $\phi ^0(x)={1\over 2}\phi ''x\cdot x$ is the Hessian
quadratic form of $\phi $ at $x_0$, we
have
\ekv{ss.26}
{
\nu _+^0(\phi ^0)=\langle 2B\phi ''x,\phi ''x\rangle \ge 0.
}

\begin{lemma}\label{ss1}
Let $\mu (x,\partial _x)=Mx\cdot \partial _x$ be a real linear \vf{} on
${\bf R}^n$. Let $n_{\pm}\in{\bf N}$, $n_++n_-=n$. Then the following two
statements are equivalent:
\smallskip
\par\noindent {\rm (A)} M has $n_+$ \ev{}s with real part $>0$ and $n_-$ \ev{}s
with real part $<0$.
\par\noindent {\rm (B)} There exists a quadratic form $G:{\bf R}^n\to {\bf R}$
of signature $(n_+,n_-)$ and a constant $C>0$, such that
\ekv{ss.27}
{
\mu (x,\partial _x)(G)\ge {1\over C}\vert x\vert ^2,\ x\in{\bf R}^n.
}
\end{lemma}
\begin{proof}
Assume first that $(n_+,n_-)=(n,0)$. If (A) holds, we know that
$$
\Vert e^{tM}x\Vert ^2\ge {1\over C}e^{t/C}\Vert x\Vert ^2,\ t\ge 0,
$$
for some constant $C>0$, and we can put (by a classical argument)
$$
G(x)=G_T(x)={1\over T}\int_0^T \Vert e^{tM}x\Vert ^2 dt,\ T\gg 1.
$$
Then
\begin{eqnarray*}
Mx\cdot \partial _x(G_T(x))={1\over T}\int_0^T Mx\cdot \partial _x\Vert
e^{tM}x\Vert ^2dt
={1\over T}\int_0^T {d\over dt}\Vert e^{tM}x\Vert ^2 dt\\
={1\over T}(\Vert e^{TM}x\Vert ^2-\Vert x\Vert ^2)\ge {1\over T}({1\over
C}e^{T/C}-1)\Vert x\Vert ^2\ge {1\over 2T}\Vert x\Vert ^2,
\end{eqnarray*}
if $T$ is large enough. Thus we get (B).

\par Conversely, if (B) holds, we have with a new constant $C>0$, that
$$
\mu (x,\partial _x)G\ge {1\over C}G,
$$
and hence
$$
{d\over dt}G(e^{tM}x)\ge {1\over C}G(e^{tM}x),
$$
so
$$
G(e^{tM}x)\ge e^{t/C}G(x),\ t\ge 0.
$$
Thus $\Vert e^{tM}x\Vert \ge C\inv e^{t/C}\Vert x¯\Vert $ for some new
positive constant $C$ and we conclude that the \ev{}s of $M$ all have
positive real parts.

\par Now consider the case of general $(n_+,n_-)$ and assume first that
(A) holds. Then we have the $M$-invariant decomposition, ${\bf
R}^n=L_+\oplus L_+$ where ${\rm dim\,}(L_{\pm})=n_\pm$ and $\sigma
({M_\vert}_{L_\pm})$ belongs to the open right half plane in the $+$ case
and to the open left half plane in the $-$ case. Hence we have positive
definite quadratic
forms $G_{\pm}$ on $L_{\pm}$
such that
$$
\pm Mx\cdot \partial _x(G_\pm )\ge {1\over C}\vert x\vert ^2,\ x\in L_{\pm}.
$$
Then $G=G_+\oplus (-G_-)$ (defined in the obvious way) has the required
properties in (B).

\par Conversely, assume that (B) holds. Let $\widetilde{L}_+$
be an $n_+$-dimensional subspace on which $G$ is positive definite. By
\no{ss.27}, we have for all $x\in{\bf R}^n$,
$$
\mu (x,\partial _x)G(x)\ge {1\over C}G(x),\ C>0,
$$
so if $x\in\widetilde{L}_+$, we get
$$
G(e^{tM}x)\ge e^{t/C}G(x),\ t\ge 0,
$$
and hence with a new constant $C>0$,
$$
\vert e^{tM}x\vert \ge {1\over C}e^{t/C}\vert x\vert ,\ t\ge 0.
$$
Similarly, if ${\rm dim\,}(\widetilde{L}_-)=n_-$ and $G$ is negative
definite on $\widetilde{L}_-$, we get
$$
\vert e^{tM}x\vert \ge {1\over C}e^{\vert t\vert /C}\vert x\vert ,\
x\in\widetilde{L}_-,\, t\le 0.
$$

\par Now we have the $M$-invariant decomposition
$$ {\bf R}^n=L_+\oplus L_0\oplus L_-, $$ where $L_+^{\bf C}$, $L_0^{\bf
C}$, $L_-^{\bf C}$ are the sums of generalized eigenspaces of $M$
corresponding to the \ev{}s with real parts $>0$, $=0$ and $<0$ respectively.
We see that necessarily, $\widetilde{L}_+\cap (L_0\oplus L_-)=0$, so
$$ n_+={\rm dim\,}\widetilde{L}_+\le {\rm
dim\,}L_+.$$
Similarly,
$$\ n_-={\rm dim\,}\widetilde{L}_-\le {\rm dim\,}L_- ,$$
so $\dim L_+=n_+$, $\dim L_-=n_-$, $\dim L_0=0$, and (A)
follows.
\end{proof}

\par If $B>0$, then \no{ss.26} with strict inequality for $x\ne 0$, and
Lemma \ref{ss1} imply that $\nu _\pm^0$ has $n_\pm$ \ev{}s
with $\pm$ real part $>0$,
where $(n_+,n_-)$ is the signature of $\phi ''(x_0)$. It follows
in that case that $F_p$ has no real \ev{}s. This last conclusion also
follows from \cite{Sj}. Indeed, in that case the quadratic approximation
$p^0(x,\xi )$ of $p$ at $(x_0,0)$ is elliptic in the sense that $\vert
p^0(x,\xi )\vert \backsim \vert x\vert ^2+\vert \xi \vert ^2$ and takes
its values in an angle $-{\pi \over 2}+\epsilon \le \arg p^0 \le {\pi \over
2}-\epsilon $ for some $\epsilon >0$.

\par Now return to the general case, when we only assume \no{ss.25} and
$\phi ''=\phi ''(x_0)$ is \nondeg{} of signature $(n_+,n_-)$.

\par We next make some remarks about the quadratic approximation $p^0$ of $p$
at $(x_0,0)$.
\begin{prop}\label{ss2} {\rm a)} Assume that the matrix $A\phi ''$ of $\nu _+^0$
has $m_{\pm}$ \ev{}s with $\pm$ real part $>0$, $m_++m_-=n$. Then there
exists a real quadratic form ${\cal G}(x,\xi )$ on ${\bf R}^{2n}$ such that
\ekv{ss.31}
{
\Re p^0((x,\xi )+i\epsilon H_{\cal G}(x,\xi ))\ge {\epsilon \over C}\vert
(x,\xi )\vert ^2,\ (x,\xi )\in{\bf R}^{2n},\,\, 0<\epsilon \ll 1.
}
\smallskip
\par\noindent {\rm b)} Conversely, assume that there exists a quadratic form
${\cal G}$ such that {\rm \no{ss.31}} holds. Then $A\phi ''$ has $n_\pm$ \ev{}s
with $\pm$ real part $>0$, where $(n_+,n_-)$ is the signature of $\phi ''(0)$.
\end{prop}

\par Recall that the condition \no{ny.17} implies the existence of ${\cal G}$
as in a) of the proposition. (The converse is not true however.
It is easy to find examples of purely imaginary quadratic forms $p^0$
for which there exist ${\cal G}$ as in a) of the proposition.)

\begin{proof} a)
Choose ${\cal G}(x,\xi )$ of the form $G(x)+\widetilde{G}(\xi )$, so that
$H_{\cal G}(x,\xi )=(\widetilde{G}'_\xi ,-G'_x)$. Recall that $p^0$ is the
quadratic approximation of $p$
at $(x_0,0)$, obtained from \no{ss.18} by freezing $A_{j,k}$
at $x_0$ and replacing $\partial _{x}\phi $ by $\phi ''x=(\phi ^0)'(x)$,
$\phi ''=\phi ''(x_0)$.

\par We get
\begin{eqnarray*}
\Re p_\epsilon ^0(x,\xi )&=&\sum_{j,k}B_{j,k}(\xi _j\xi _k+\partial
_{x_j}\phi ^0\partial _{x_k}\phi ^0)\\ &&+2\epsilon \Big(
\sum_{j,k}C_{j,k}(\partial
_{x_k}\phi ^0)\partial _{x_j}G-\sum_{j,k}C_{j,k}(\phi ''\partial _\xi
\widetilde{G})_k\xi _j\Big) +{\cal O}(\epsilon ^2\vert (x,\xi )\vert ^2),
\end{eqnarray*}
so it suffices to have
\begin{eqnarray}
C\phi ''x\cdot \partial _xG&\ge &{1\over C}\vert x\vert ^2\hbox{ on
}{\phi ''}\inv{\cal N}(B),\label{ss.32}\\
\phi ''C\xi \cdot \partial _\xi \widetilde{G}&\ge &
{1\over C}\vert \xi \vert ^2\hbox{ on
}{\cal N}(B)\label{ss.33},
\end{eqnarray}
where we also used the antisymmetry of $C$ in the last \e{}. This would follow
from
$$
A\phi ''x\cdot \partial _xG\ge {1\over C}\vert x\vert ^2,\ \phi
'' A\xi \cdot \partial _\xi \widetilde{G}\ge {1\over C}\vert \xi
\vert ^2,\ x,\xi \in {\bf R}^n,
$$
and in order to find such \fu{}s $G$ and $\widetilde{G}$,
it suffices to apply Lemma
\ref{ss1} to the isospectral matrices $A\phi ''$ and $\phi '' A$.

\par b) Let ${\cal G}$ be as in \no{ss.31}. The quadratic form $$
p_\epsilon ^0(x,\xi )=p^0(\exp i\epsilon H_{\cal G}(x,\xi ))=p^0(x,\xi
)+i\epsilon H_{\cal G}p^0+{\cal O}(\epsilon ^2)
$$
is elliptic on the real phase space and takes its values in an angle
$$
e^{i[-{\pi \over 2}+{\epsilon \over C},{\pi \over 2}-{\epsilon \over
C}]}[0,+\infty[ ,
$$
so we know from \cite{Sj} that $F_{p_\epsilon ^0}$ has no real \ev{}s.

\par
On the other hand $p_\epsilon ^0$ and $p^0$ are related by a \ctf{},
so $F_{p_\epsilon ^0}$
and $F_{p_0}$ are isospectral. Hence $A\phi ''$ has $m_\pm$ \ev{}s with $\pm$
real part $>0$, where $m_++m_-=n$. To see that $m_\pm =n_\pm$, we just
replace $B$ by $B+\delta 1$, $0<\delta \ll 1$, to reduce ourselves to the
elliptic case, and apply the observation after the proof of Lemma \ref{ss1}.
\end{proof}

\subsection{The subprincipal symbol}\label{ssc}

\par We next look at the subprincipal term in \no{ss.12}. Write
$$
A(dx_k)=\sum_\nu  A_{\nu ,k}\partial _{x_\nu },\
A(dx_k)^\rfloor =\sum_\nu A_{\nu ,k}\partial _{x_\nu }^\rfloor ,
$$
so the second sum in \no{ss.12} becomes
$$
2h\sum_{j,k,\nu }\phi ''_{j,k}A_{\nu ,k}\circ dx_j^\wedge  \partial
_{x_\nu }^\rfloor = 2h\sum_{j,\nu}(\phi ''\circ {\trans A})_{j,\nu }\circ
dx_j^\wedge \partial _{x_\nu }^\rfloor
$$
which simplifies further to
\ekv{ss.35}
{
2h\sum _j (\phi ''\circ {\trans A})(dx_j ) ^\wedge \partial _{x_j
}^\rfloor .
}

\par Now we restrict the attention to a \nondeg{} critical point $x_0$
of $\phi $ and we shall compute the subprincipal symbol of $-\Delta _A$
at the corresponding doubly characteristic point $(x_0,0)$. At that
point $\phi ''\circ {\trans A}:T_{x_0}^*M\to T_{x_0}^*M$ is invariantly
defined and it is easy to check that \no{ss.35} is also invariantly
defined: we get the same quantity if we replace $dx_1,..,dx_n$, $\partial
_{x_1},..,\partial _{x_n}$, by $\omega _1,..,\omega _n$, $\omega
_1^*,..,\omega _n^*$, where $\omega _1,..,\omega _n$ is any basis in the
complexified cotangent space and $\omega _1^*,..,\omega _n^*$ is the dual
basis of tangent vectors for the natural bilinear pairing.

\par Assume that the equivalent conditions of Proposition \ref{ss2} hold
and denote the corresponding
\ev{}s (that are also the \ev{}s of $\phi ''\circ {\trans A}$) by $\lambda
_1,..,\lambda _n$ with $\Re \lambda _j>0$
for $1\le j\le n_+$ and with $\Re \lambda _j <0$ for $n_++1\le j\le
n=n_++n_-$. The \ev{}s of $F_p$ are then $\pm 2i\lambda _j$
(in view of \no{ss.24}
and the isospectrality of $F_p$ and $iF_q$ reviewed prior to Lemma
\ref{as3})
, so
\ekv{ss.36}
{
\widetilde{{\rm tr\,}}F_p:={1\over i}\sum_{\mu \in \sigma (F_p)\atop \Im
\mu >0}\mu =\sum_1^{n_+}2\lambda _j-\sum_{n_++1}^n2\lambda _j .
}

\par The subprincipal symbol of the first term in \no{ss.12} (at $(x_0,0)$)
is equal to
\ekv{ss.37}
{
\sum_{j,k}A_{j,k}{1\over 2i}\{ -i\xi _k+\partial _{x_k}\phi ,i\xi _j+\partial
_{x_j}\phi \}
=-\sum_{j,k}A_{j,k}\phi ''_{j,k}=-\tr (A\phi '')=-\sum_1^n \lambda _j.
}

\par The \ev{}s of $\sum_j (\phi ''\circ {\trans A})(dx_j)^\wedge \partial
_{x_j}^\rfloor$ on the space of $m$-forms are easily calculated, if we
replace $dx_1,..,dx_n$ by a basis of eigenvectors $\omega _1,..,\omega _n$
of $\phi ''{\trans A}$, so that
$$
(\phi ''\circ {\trans A})(\omega _j)=\lambda _j\omega _j,
$$
and $\partial _{x_j}$ by the corresponding dual basis vectors $\omega _j^*$. (Here we assume to
start with that there are no Jordan blocks. This can be achieved by an
\ably{} small perturbation of $A$, and we can extend the end result of our
calculation to the general case by continuity.) We get
\ekv{ss.38}
{
\sum _j (\phi ''\circ {\trans A})(dx_j )^\wedge \partial _{x_j
}^\rfloor=\sum_j \lambda _j\omega _j^\wedge  {\omega _j^*}^\rfloor .
}
A basis of eigenforms of this \op{} is given by $\omega _{j_1}\wedge
..\wedge \omega _{j_m}$, $1\le j_1<j_2<..<j_m\le n$ and the corresponding
\ev{}s are $\lambda _{j_1}+..+\lambda _{j_m}$.

\par Let $S_P$ be the subprincipal symbol of \no{ss.12} at $(x_0,0)$. Then
the \ev{}s of
$$
{1\over 2}\widetilde{\tr}F_p+S_P,\hbox{ acting on }m\hbox{ forms }
$$
are
\eekv{ss.39}
{
\sum_1^{n_+}\lambda _j-\sum_{n_++1}^n\lambda _j-\sum_1^n\lambda
_j+2(\lambda _{j_1}+..+\lambda _{j_m})
}
{=2(\lambda _{j_1}+..+\lambda _{j_m}-\sum_{n_++1}^n\lambda _j),\ 1\le
j_1<..<j_m\le n.}
We conclude that if $m\ne n_-$, then all the \ev{}s have a real part $>0$
and if $m=n_-$, then precisely one \ev{} is equal to 0, while the others
have positive real part.

\subsection{A symmetry for adjoints}\label{ssd}
\par Our last remark in this section concerns symmetry relations for
the $A,*$ adjoints. If $D:L^2(\Omega ;\wedge ^kT^*\Omega )\to L^2(\Omega
;\wedge ^jT^*\Omega )$, then a simple calculation shows that
$$
D^{A,*}=(\wedge ^k(\trans{A}))\inv D^*(\wedge ^j(\trans{A})),
$$
where $D^*$ denotes the adjoint \wrt{} the measure $\mu $.

\par We also have
$$
(u\vert v)_A=\overline{(v\vert u)}_{\trans{A}}.
$$
Playing with these relations we see that
$$
D=(D^{A,*})^{\trans{A},*}.
$$

This can be applied to $-\Delta _A$ and we get
$$
(-\Delta _A)^{\trans{A},*}=-\Delta _{\trans{A}}.
$$

\section{The double well case}\label{2w}
\setcounter{equation}0

\par In this section we assume that
\ekv{dw.1}
{
M=\R ^n,\  A=A(x)\hbox{ is \indep{} of  }x,\hbox{ and  invertible}.
}
We decompose $A$ as in \no{ss.13} and assume \no{ss.25}. Let $\phi\in
C^{\infty}(\R ^n;\R)$ be such that 
\ekv{dw.2}
{
\partial ^\alpha_x \phi(x) ={\cal O}(1),\quad
\partial ^{\alpha}_x \left(\langle B\dd_x\phi ,\dd_x\phi\rangle\right)
= {\cal O}(1),\quad \abs{\alpha}\geq 2. 
}

\par Consider $P^{(q)}=-\Delta _A^{(q)}$ which according to \no{ss.14} becomes
\eekv{dw.3}
{
P^{(q)}&=&\sum_{j,k}hD_{x_j}B_{j,k}hD_{x_k}+\sum_{j,k}(\partial_{x_j}\phi)B_{j,k}(\partial_{x_k}\phi)-h{\rm tr\,}(B\phi '')
}{&&+
\sum_{j,k}((\partial_{x_k}\phi )C_{j,k}h\partial_{x_j}+h\partial_{x_j}\circ C_{j,k}(\partial_{x_k}\phi ))+2h\sum_{j,k}(\partial _{x_j}\partial _{x_k}\phi )dx_j^{\wedge} A(dx_k)^\rfloor .
}
Apart from the third and the  last terms which are ${\cal O} (h)$ with
all their derivatives, this is of the form \no{ap.1} with $b_{j,k}=B_{j,k}$,
$c_j(x)=\sum_k C_{j,k}\dd _{x_k}\phi$, $p_0(x)=\langle
B\dd_x\phi ,\dd_x\phi\rangle $.
We define $p_2,p_1,p_0$ as in Section \ref{ap} and see that
\no{ap.4}--\no{ap.7} hold.

\par  Assume that
\ekv{dw.4}
{
\phi\hbox{ is a Morse \fu{} with critical points }x_1,...,x_N\in\R ^n,
}
\ekv{dw.5}
{
\vert \phi ' (x)\vert\ge 1/C,\ \vert x\vert\ge C.
}
Let ${\cal C}=\{ \rho _j;\, j=1,...,N\} $ where $\rho _j=(x_j,0)\in T^* \R^n$.
Then $\rho _j$ are critical points with critical value 0 for $p_2,p_1,p_0$.
Since $A$ is invertible, ${\cal N}(B)\cap {\cal N}(C)=0$, and hence
\no{ny.1} holds. We adopt the dynamical assumptions \no{ny.17},
\no{ny.18} (or equivalently \no{ny.17}, \no{ny.19}) and \no{ny.20}).
Then we can apply the results of Sections \ref{as}, \ref{exp} to
$P^{(q)}$ since neither the presence of the \bdd{} subprincipal
symbol in \no{dw.3} nor the non-scalar nature of the \op{}s has any
serious influence. In the preceding section we saw that we are in the
case when the conditions of Proposition \ref{ss2} are fulfilled.
The only change in Theorem \ref{as1} is
that the $\mu _{j,k}$ in \no{as.8} are of the form
\ekv{2w.0}
{
\mu _{j,k}={1\over i}\sum_{k=1}^n \nu _{j,k,\ell}\lambda _\ell +{1\over
2}{\widetilde{\rm tr}\,}(F_p(\rho _j))+\gamma _{j,k},
}
where $\gamma _{j,k}$ is any \ev{} of the subprincipal symbol
$S_{P^{(q)}}$ at $(x_j,0)$. From the calculations in Subsection \ref{ssc}
we notice that the $\mu _{j,k}$ will be confined to a sector $\{
0\}\cup \{ |{\rm arg\, }z\vert <{\pi /2}-1/C\}$ around
$[0,+\infty [$ and it is precisely when $x_j$ is of index $q$ (i.e. when
the Hessian of $\phi $ at $x_j$ has precisely $q$ negative \ev{}s) that one
of the $\mu _{j,k}$ may be equal to $0$.

\par We now add more specific conditions for the double well case.
Assume that
\eekv{2w.0.5}
{&&\phi  \hbox{ has precisely three critical points, two local}}
{&&\hbox{minima }U_{\pm 1}, \hbox{ and a "saddle point" }U_0\hbox { of index one.}} Then $\phi
(x)\to +\infty $ with $\phi (x)\ge {1\over C}\vert x\vert $, for $\vert
x\vert \ge C$.

\par Put $S_j=\phi (U_0)-\phi (U_j)$, $j=\pm 1$, so that $S_j>0$. The set
$\phi \inv (]-\infty ,\phi (U_0)[)$ has precisely two connected
components $D_j$, $j=\pm 1$, determined by the condition $U_j\in D_j$.
Under these assumptions we know that $P^{(0)}=-\Delta _A^{(0)}$ has precisely two
\ev{}s $\mu _0,\mu _1=o(h)$ spanning a corresponding 2-dimensional
spectral subspace $E^{(0)}$. Actually one of these two \ev{}s, say $\mu _0$,
is equal to $0$ with $e^{-\phi /h}$ as the corresponding \ef{} and since a
truncation of this \fu{} can be used as a quasimode near each of $U_{\pm 1}$
we also know that $\mu _{1}={\cal O}(h^\infty )$ (cf \no{2w.4}).
Moreover, $-\Delta _A^{(1)}$ has
precisely one \ev{} $\widetilde{\mu }_1=o(h)$ and $-\Delta _A^{(k)}$
has no \ev{}s $=o(h)$ for $k\ge 2$.
Since our \op{}s are real we know that the spectra are symmetric around the real axis, hence $\mu _0,\mu _1,\widetilde{\mu }_1$ are real.
From the intertwining relations
$$
-\Delta _A^{(1)}d_\phi =d_\phi (-\Delta _A^{(0)}),\ -\Delta _A^{(0)}d_\phi
^{A,*}=-d_\phi ^{A,*}\Delta _A^{(1)},
$$
we then also know {}that $\widetilde{\mu }_1=\mu _1$.

\par In fact, when $B>0$ it follows from the ellipticity and the estimates
in Section \ref{ap} that all eigenforms and generalized eigenforms
corresponding to an \ev{} in $D(0,Ch)$ belong to ${\cal S}(\R ^n)$, so
if $\mu_1\ne 0$ and $(-\Delta _A^{(0)}-\mu_1 )u=0$, $u\in L^2$, then
$u\in{\cal S}(\R ^n)$ and $0\ne d_\phi u\in{\cal S}$ is an eigenform for
$-\Delta_A^{(1)}$ with the same \ev{}. A priori we cannot exclude that
$\mu_1=0$ and that $-\Delta_A^{(0)}u={\rm Const\,} e^{-\phi /h}$. Then
again, $0\ne d_\phi u\in\cS$ is a corresponding eigenvector of
$-\D_A^{(1)}$. In the general case, we let $0<B_\epsilon \to B=:B_0$
when $\epsilon \searrow 0$. (Take for instance $B_\epsilon =B+\epsilon I$.) On a small circle $D(0,h/C)$ we know that
$(P_\epsilon ^{(q)}-z)\inv =\OO({1\over h})$ \ufly{} for $0\le \epsilon
\ll 1$. If $(P_0^{(q)}-z)u=v$, $u,v\in\cS $ (for $\epsilon =0$), then
$(P_\epsilon ^{(q)}-z)u=v+r_\epsilon $, $r_\epsilon \to 0$, so
$(P_\epsilon ^{(q)}-z)
\inv v=u-(P_\epsilon ^{(q)}-z)\inv r_\epsilon \to u$, $\epsilon \to 0$.
Since $(P^{(q)}-z)(\cS )$ is dense in $L^2$ we conclude that
$(P_\epsilon ^{(q)}-z)\inv\to (P^{(q)}-z)\inv$ strongly for
$\vert z\vert =h/C$. We have then the corresponding fact for the finite
rank spectral projections and their compositions with $P_\epsilon ^{(q)}$:
$$
{1\over 2\pi i}\int_{\vert z\vert =h/C}(z-P_\epsilon )\inv dz \hbox{ and }
{1\over 2\pi i}\int_{\vert z\vert =h/C}z(z-P_\epsilon )\inv dz.
$$
These are finite rank \op{}s and converge therefore in norm. It follows that
\ekv{2w.0.6}{
\sigma (P_\epsilon ^{(q)})\cap D(0,h/C)\to \sigma (P^{(q)})\cap D(0,h/C),\hbox{ when }\epsilon \to 0,
}
so we get $\widetilde{\mu }_1=\mu _1$ also in the general case.

\par Let $\chi _j\in C_0^\infty (D_j)$ be equal to 1 on $D_j\cap\phi \inv
(]-\infty ,\phi (U_0)-\epsilon _0])$ for $\epsilon _0>0$ fixed but \ably{}
small. Consider
\ekv{2w.1}
{
f_j=h^{-n/4}c_j(h)e^{-{1\over h}(\phi (x)-\phi (U_j))}\chi _j(x),\ j=\pm 1,
}
where $c_j\sim c_{j,0}+hc_{j,1}+...>0$ is a normalization constant with
$c_{j,0}>0$, such that
\ekv{2w.2}
{
\Vert f_j\Vert =1.
}
We also have
\ekv{2w.3}
{
P^{(0)}(f_j)=[P^{(0)},\chi _j](c_je^{-{1\over h}(\phi (x)-\phi
(U_j))})={\cal O}(h^{-N_0}e^{-{1\over h}(S_j-\epsilon _0)}),
}
for some $N_0>0$.

\par If
\ekv{2w.3.5}
{
\Pi ^{(0)}={1\over 2\pi i}\int_\gamma  (z-P^{(0)})\inv dz,\ \gamma
=\partial D(0,{h\over C})
}
is the spectral projection of $P^{(0)}$ onto $E^{(0)}$ we know from
Theorem \ref{as2} that $\Pi ^{(0)}={\cal O}(1)$. It follows from \no{2w.3} that
\ekv{2w.4}
{
e_j:=\Pi ^{(0)}f_j=f_j+{\cal O}(h^{-N_1}e^{-{1\over
h}(S_j-\epsilon _0)})\hbox{ in }L^2.
}
In fact, we write \no{2w.3} as $P^{(0)}f_j=r_j,$
$$
(z-P^{(0)})(f_j)=zf_j-r_j,
$$
$$
(z-P^{(0)})\inv f_j={1\over z}f_j+(z-P^{(0)})^{-1}z^{-1}r_j
$$
and integrate, using the bounds on the resolvent provided by Theorem \ref{as2}.

\par From \no{2w.4} we see that
\ekv{2w.5}
{
\Vert e_j\Vert ^2=1+{\cal O}(h^{-N_2}e^{-{1\over
h}(S_j-\epsilon _0)}),
}
\ekv{2w.6}
{
(e_1\vert e_{-1})={\cal O}(h^{-N_2}e^{-{1\over
h}(S_{\rm min}-\epsilon _0)}),
}
where
\ekv{2w.9}
{
S_{\rm min}=\min (S_{-1},S_1).
}

\par Let $E^{(1)}$ be the one-dimensional eigenspace of $P^{(1)}$
corresponding to $\mu _1$. From an easy extension of Theorem \ref{exp1} to
the non-scalar case with the presence of other non-resonant wells ($U_{\pm
1}$) as in Remark \ref{exp2}, we know that $E^{(1)}$ is generated by an eigenform
\ekv{2w.10}
{
e_0(x;h)=\chi _0(x)e^{-{1\over h}\phi _+(x)}h^{-{n\over 4}}a_0(x;h)+{\cal
O}(e^{-S_0/h}),
}
where $\chi _0\in C_0^\infty ({\rm neigh\,}(U_0))$ is equal to one near
$U_0$, $S_0>0$,
$$a_0(x;h)\sim \sum_0^\infty  a_{0,k}(x)h^k$$
is a symbol as in Theorem \ref{exp1} with $a_{0,0}(U_0)\ne 0$, and $\phi
_+\in C^\infty ({\rm neigh\,}(U_0);[0,\infty [)$ satisfies
\ekv{2w.10.1}
{
\phi _+(x)\backsim \vert x-U_0\vert ^2
}
and solves the eikonal \e{}
\ekv{2w.10.2}
{
q(x,\phi _+'(x))=0,
}
with $q=p_2+p_1-p_0$,
$$p_2=\langle B(x)\xi ,\xi \rangle ,\ p_1(x,\xi
)=2\langle C(x)\phi '(x),\xi \rangle ,\ p_0(x)=\langle B(x)\phi
'(x),\phi '(x)\rangle .$$

\par $\Lambda _{\phi _+}$ is the stable outgoing
\mfld{} through $(U_0,0)$ for the $H_q$-flow and recall that $\phi
_+''(U_0)>0$ by Proposition \ref{as4}. (Similarly we have a stable
incoming \mfld{} $\Lambda _{\phi _-}$.) Let $k_\pm$ be the number of \ev{}s
of the linearization of ${{H_q}_\vert }_{\Lambda _\phi }$ at that point
with $\pm$ real part $>0$,
so that $k_++k_-=n$. Let $K_+, K_-\subset \Lambda _\phi $ be the
corresponding stable outgoing and incoming sub\mfld{}s of dimension $k_+$
and $k_-$ respectively. Then $K_+\subset \Lambda _{\phi _+}$, $K_-\subset
\Lambda
_{\phi _-}$ and $\phi -\phi (U_0)-\phi _\pm$ vanishes to the second order
on $\pi
_x(K_\pm )$. Since $\phi ''(U_0)$ has signature $(n-1,1)$, we conclude
that ${\rm dim\,}K_+=n-1$, ${\rm dim\,}K_-=1$. (This also follows from
Proposition \ref{ss2}.) It is also clear that
$\Lambda _\phi ,\Lambda _{\phi _\pm}$ intersect cleanly along $K_{\pm}$,
so we get
\ekv{2w.10.3}
{
\phi _+-(\phi -\phi (U_0))\backsim {\rm dist\,}(x,\pi _x(K_+))^2,\ \phi -\phi (U_0)-\phi
_-\backsim {\rm dist\,}(x,\pi _x(K_-))^2.
}

\par We next make some remarks about the adjoint \op{} $-\Delta _{\trans
A}=(-\Delta _A)^{\trans A,*}$ (cf Subsection \ref{ssd}). The principal
symbol is $p_2-ip_1+p_0=p(x,-\xi )=\check{p}(x,\xi )$ and the
corresponding real "$q$"-symbol is $\check{q}(x,\xi )=q(x,-\xi )$. Since
our dynamical conditions are invariant under a change of sign of the
$H_{p_1}$-direction, all our assumptions are equally valid for $-\Delta
_{\trans A}$. This also holds for the geometric discussion above, so if
$\Lambda _{\phi _+^*}$, $\Lambda _{\phi _-^*}$ denote the outgoing and
incoming $H_{\check{q}}$-invariant Lagrangian \mfld{}s through $(U_0,0)$
and $K_\pm^*\subset \Lambda _\phi $ the outgoing/incoming \mfld{}s for
${{H_{\check{q}}}_\vert}_{\Lambda _\phi }$ (noting that $\check{q}=0$ on
$\Lambda _\phi $), then ${\rm dim\,}K_+^*=n-1$, ${\rm dim\,}K_-^*=1$ and
\ekv{2w.10.4}
{
\phi _+^*-(\phi -\phi (U_0))\backsim {\rm dist\,}(x,\pi _x(K_+^*))^2,\ \phi -\phi (U_0)-\phi
_-^*\backsim {\rm dist\,}(x,\pi _x(K_-^*))^2.
}

\par In view of the general relation
\ekv{2w.10.5}
{
J_*(H_q)=-H_{\check{q}},\hbox{ where }J:(x,\xi )\mapsto (x,-\xi ),
}
we see that $\Lambda _{\phi _-^*}=J(\Lambda _{\phi _+})$, $\Lambda _{\phi
_+^*}=J(\Lambda _{\phi _-})$, or more simply
\ekv{2w.10.6}
{
\phi _-^*=-\phi_+,\ \phi _+^*=-\phi_- ,
}
giving in particular from \no{2w.10.3}, \no{2w.10.4},
\ekv{2w.10.7}
{
\phi -\phi (U_0) +\phi _+^*\backsim {\rm dist\,}(x,\pi _x(K_-))^2,\ \phi -\phi (U_0)+\phi
_+\backsim {\rm dist\,}(x,\pi _x(K_-^*))^2.
}

\par Let $\mu _0^*=0$, $\mu _1^*$ be the two \ev{}s of $P_*^{(0)}:=-\Delta
_{\trans A}^{(0)}$ that are $o(h)$ and let $\Pi _*^{(0)}$ be the
 spectral projection onto the corresponding spectral subspace $E_*^{(0)}$.
 Then $e_j^*=\Pi _*^{(0)}f_j$, $j=\pm 1$, span $E_*^{(0)}$, and satisfy
 \no{2w.4}. Similarly for $P_*^{(1)}=-\Delta _{\trans A}^{(1)}$ we have the
 generating eigenform
\ekv{2w.10.8}
{
e_0^*(x;h)=\chi _0(x)e^{-{1\over h}\phi _+^*(x)}h^{-{n\over 4}}a_0^*(x;h)+{\cal
O}(e^{-S_0/h}) \hbox{ in }L^2
}
for the one dimensional eigenspace $E_*^{(1)}$ corresponding to $\mu _1^*$.

\par Now, using that our eigenvalues and \op{}s are real, we know by
duality that
\ekv{2w.10.9}
{
\mu _1^*=\mu _1,
}
and that $(E_*^{(0)},E^{(0)})$ and $(E_*^{(1)},E^{(1)})$ are dual pairs for
the scalar products $(u\vert v)_{L^2}$ and $(u\vert v)_A$ respectively. In fact, $((z-P_*^{(q)})\inv)^{A,*}=(\overline{z}-P^{(q)})\inv$.

\par From Subsection \ref{ssc} we know that $a_{0,0}(U_0)$ is an
eigenvector corresponding to the negative \ev{} of $\phi ''\circ {\trans
A}$ at $U_0$, and $a_{0,0}^*(U_0)$ is an
eigenvector corresponding to the negative \ev{} of $\phi ''\circ
A$. Since $(\phi ''\circ A)^{A,*}=\phi ''\circ {\trans A}$, we know that
the two \ev{}s are equal and that the $A$-product of the two eigenvectors
is $\ne 0$;
\ekv{2w.11}
{
(a_{0,0}^*(U_0)\vert a_{0,0}(U_0))_A\ne 0.
}
It follows that $(e_0^*\vert e_0)_A\backsim 1$ and after renormalization
of $e_0^*$ we may assume that
\ekv{2w.12}
{
(e_0^*\vert e_0)_A=1.
}
Similarly, using \no{2w.2}
\ekv{2w.13}
{
(e_j^*\vert e_k)=\delta _{j,k}+{\cal O}(e^{-{1\over Ch}}), \ j,k=\pm 1.
}

\par Let $\pmatrix{\lambda _{-1} &\lambda _1}$ be the matrix of $d_\phi
:E^{(0)}\to E^{(1)}$ \wrt{} the bases $\pmatrix{e_{-1}, &e_1}$ and $(e_0)$.
(Strictly speaking, we approximate our operators by elliptic ones as in
\no{2w.0.6} and pass to the limit.)
Let
$$
\pmatrix{\lambda _{-1}^*\cr \lambda _1^*}
$$
be the matrix of $d_\phi ^{A,*}$ for the same bases. The \ev{} $\mu _1$ can
be viewed as the second \ev{} of $d_\phi ^{A,*}d_\phi :E^{(0)}\to E^{(0)}$
or equivalently as the scalar $d_\phi d_\phi ^{A,*}:E^{(1)}\to E^{(1)}$ (using also that $P^{(2)}$ has no \ev{} $=o(h)$).
Either way, we get
\ekv{2w.14}
{
\mu _1=\lambda _{-1}^*\lambda _{-1}+\lambda _{1}^*\lambda _{1}.
}

\par We get
\ekv{2w.15}
{
\overline{\lambda }_k=(e_0^*\vert d_\phi e_k)_A,\ k=\pm 1,
}
\ekv{2w.16}
{
\overline{\lambda }_j^*=(g_j\vert d_\phi^{A,*} e_0)_A,\ j=\pm 1,
}
where
\ekv{2w.16.5}
{
\pmatrix{g_{-1} &g_1}=\pmatrix{e_{-1}^* &e_1^*}(1+{\cal O}(e^{-{1\over
Ch}}))
}
is the base in $E_*^{(0)}$ that is dual to $\pmatrix{e_{-1} &e_1}$. Here
the complex conjugate signs are superfluous since we work with real
\op{}s, \ev{}s and functions.

\par Let $\chi \in C_0^\infty ({\rm neigh\,}(U_0);[0,1])$ be equal to 1
near $U_0$. Using that $d_\phi e_{-1}=\lambda _{-1}e_0$, we get, dropping
the bars from now on,
\eeekv{2w.17}
{\lambda _{-1}&=&(e_0^*\vert d_\phi e_{-1})_A}
{&=& (e_0^*\vert \chi d_\phi e_{-1})_A+\lambda _{-1}(e_0^*\vert (1-\chi
)e_0)_A}
{&=& (e_0^*\vert \chi d_\phi e_{-1})_A+{\cal O}(e^{-{1\over Ch}})\lambda
_{-1}.}

\par Here
\ekv{2w.18}
{
(e_0^*\vert \chi d_\phi e_{-1})_A=(e_0^*\vert [\chi ,d_\phi
]e_{-1})_A+(d_\phi ^{{\trans A},*}e_0^*\vert \chi e_{-1})_A.
}
Now the matrix of $d_\phi ^{{\trans A},*}:E_*^{(1)}\to E_*^{(0)}$ \wrt{}
the dual bases is the adjoint of the one of $d_\phi :E^{(0)}\to E^{(1)}$,
so $d_\phi ^{{\trans A},*}e_0^*=\lambda _{-1}g_{-1}+\lambda _1g_1$, and
expressing $g_j$ as linear combinations of the $e_j^*$ by means of
\no{2w.16.5} and using \no{2w.4} for the $e_{\pm 1}^*$ we see that the
last term in \no{2w.10.7} is of the form
\ekv{2w.18.5}
{(d_\phi ^{{\trans A},*}e_0^*\vert \chi e_{-1})_A={\cal O}(e^{-{1\over
Ch}})\lambda
_{-1}+{\cal O}(e^{-{1\over Ch}})\lambda
_{1}.}
Thus we have obtained
\ekv{2w.19}
{
(1+{\cal O}(e^{-{1\over Ch}}))\lambda _{-1}+{\cal O}(e^{-{1\over
Ch}})\lambda _1=(e_0^*\vert [\chi ,d_\phi ]e_{-1})_A=-h(e_0^*\vert
(d\chi )\wedge e_{-1})_A,
}
and we shall study the last expression. The contribution from the
remainder in \no{2w.4} is ${\cal O}(h^{-N_2})\exp {1\over
h}(-S_{-1}+\epsilon _0-{1\over C})={\cal O}(1)e^{-{1\over
h}(S_{-1}+{1\over 2C})}$ if we choose $\epsilon _0$ small enough.
A similar estimate holds for the contribution from the remainder term in
\no{2w.10.8}. As we
shall see, the contribution from the leading terms in \no{2w.4}, \no{2w.10.8}
 will be
larger. It is equal to
\ekv{2w.20}
{
-c_{-1}(h)h^{1-{n\over 2}}\int \chi _{-1}(x)\langle A(x)a_0^*(x;h)\vert d\chi
(x)\rangle e^{-{1\over h}(\phi _+^*(x)+\phi (x)-\phi (U_{-1}))}dx.
}
Here by \no{2w.10.7},
\ekv{2w.21}
{
\phi _+^*(x)+\phi (x)-\phi (U_{-1})
=\phi _+^*+(\phi (x)-\phi (U_0))+S_{-1}
\backsim S_{-1}+{\rm dist\,}(x,\pi
_x(K_{-}))^2,
}
so we expect \no{2w.20} to behave like some power of $h$ times $e^{-{1\over
h}S_{-1}}$ with the main contribution coming from a \neigh{} of ${\rm
supp\,}(d\chi )\cap {\rm supp\,}(\chi _{-1})\cap\pi _x(K_{-})$. Now $\phi
(x)-\phi (U_0)\backsim -\vert x-U_0\vert ^2$ on $\pi _x(K_{-})$ while
$\chi _{-1}$ has its support in $D_{-1}$ and equals 1 in the subset of
$D_{-1}$ where $\phi (x)-\phi (U_0)\le -\epsilon _0$
with $\epsilon _0$
\ably{} small. Because of the presence of $d\chi $ which has its support in
an annular region around $U_0$, we see that $\chi _{-1}(x)=1$ in ${\rm
neigh\,}(D_{-1}\cap {\rm supp\,}(d\chi )\cap \pi _x(K_-))$, so we can
forget about $\chi _{-1}$ in \no{2w.20} and just integrate over ${\rm
neigh\,}({\rm supp\,}(d\chi )\cap \pi _x(K_{-1}),D_{-1})$.

\par Let us look at $\langle A(x)a_{0,0}^*(x;h)\vert d\chi (x_1)\rangle $
at a point $x_1\in {\rm supp\,}(d\chi )\cap \pi _x(K_-)$. At $U_0$ we know
that $a_{0,0}^*$ is an eigenvector of $\phi ''\circ A$ associated to the
negative \ev{}, so $Aa_{0,0}^*$
is a corresponding eigenvector for $A\phi ''$ which is the linearization
of ${1\over 2}{{H_q}_\vert}_{\Lambda _\phi }$. Hence $Aa_{0,0}^*$ at $U_0$ is tangent
to $\pi _x(K_-)$. This will remain approximately true at $x_1$ since
the latter point is close to $U_0$. Choosing $\chi $ to be a "circular"
standard cut-off, we see that $\langle Aa_{0,0}^*,d\chi \rangle $ is
non-vanishing of constant sign at every point $x_1\in \pi _x(K_-)\cap
D_{-1}$ where $d\chi \ne 0$.

\par By stationary phase it is now clear that the integral \no{2w.20} is
equal to
\ekv{2w.22}
{
h^{1\over 2}\ell_{-1}(h)e^{-{1\over h}S_{-1}},\ \ell_{-1}\sim
\ell_{-1,0}+h\ell_{-1,1}+...,\quad \ell_{-1,0}\ne 0.
}
Returning to \no{2w.19} and modifying $\ell_{-1}$ by an exponentially small
term, we get
\ekv{2w.23}
{
(1+{\cal O}(e^{-{1\over Ch}}))\lambda _{-1}+{\cal O}(e^{-{1\over
Ch}})\lambda _{1}=h^{1\over 2}\ell_{-1}(h)e^{-{1\over h}S_{-1}}.
}
Similarly,
\ekv{2w.24}
{
{\cal O}(e^{-{1\over Ch}})\lambda _{-1}+(1+{\cal O}(e^{-{1\over
Ch}}))\lambda _{1}=h^{1\over 2}\ell_{1}(h)e^{-{1\over h}S_1},
\ \ell_{1}\sim
\ell_{1,0}+h\ell_{1,1}+...,\ \ell_{1,0}\ne 0.
}
Inverting the system, we get
\ekv{2w.25}
{
\pmatrix{\lambda _{-1}\cr \lambda _1}=(1+{\cal O}(e^{-{1\over
Ch}}))\pmatrix{h^{1\over 2}\ell_{-1}(h)e^{-{1\over h}S_{-1}}\cr
h^{1\over 2}\ell_{1}(h)e^{-{1\over h}S_{1}}}.
}

\par Now turn to $\lambda _j^*$ in \no{2w.16}. In view of \no{2w.16.5}, we
have
\ekv{2w.26}
{
\pmatrix{\lambda _{-1}^*\cr \lambda _1^*}=(1+{\cal O}(e^{-{1\over
Ch}}))\pmatrix{\alpha _{-1}\cr \alpha _1},
}
where
$$
\alpha _j=(d_\phi e_j^*\vert e_0)_A=(e_0\vert d_\phi e_j^*)_{{\trans A}}
$$
which can be identified with the expression \no{2w.15} after replacing
$A$ by $\trans A$ and making the corresponding substitutions, $e_0^*\to
e_0$, $e_j\to e_j^*$. Hence we have the analogue of \no{2w.25},
\eekv{2w.27}
{
&\hskip -1truecm \pmatrix{\lambda _{-1}^*\cr \lambda _1^*}=(1+{\cal O}(e^{-{1\over
Ch}}))
\pmatrix{\alpha _{-1}\cr \alpha _1}
=(1+{\cal O}(e^{-{1\over Ch}}))
\pmatrix{h^{1\over 2}\ell_{-1}^*(h)e^{-{1\over h}S_{-1}}\cr
h^{1\over 2}\ell_{1}^*(h)e^{-{1\over h}S_{1}}},& }
{&\ell_j^*(h)\sim \ell_{j,0}^*+h\ell_{j,1}^*+...,\ \ell_{j,0}^*\ne 0.&}

We finally claim that $\ell_{j,0}\ell_{j,0}^*>0$. Indeed, this number is
real and different form zero and if we deform our matrices to reach the
selfadjoint case (with $A>0$) we see that we have a positive sign).

\par Combining this with \no{2w.14} we get the main result of this work:
\begin{theo}\label{2w1}
Let $P=-\D_A^{(0)}$ where we assume {\rm \no{dw.1}}, {\rm \no{dw.2}}, {\rm \no{2w.0.5}}.
We also assume that $P$ satisfies the additional dynamical conditions
{\rm \no{ny.17}}, {\rm \no{ny.18}} (or equivalently {\rm \no{ny.17}}, {\rm \no{ny.19}}) and
{\rm \no{ny.20}}. Then for $C>0$ large enough,  $P$ has precisely {\rm 2} \ev{}s, $0$
and $\mu _1$  in the disc $D(0,h/C)$ when $h>0$ is small enough. Here
$\mu _1$ is real and of the form
\ekv{2w.28}
{
\mu _1=h(a_1(h)e^{-2S_1/h}+a_{-1}(h)e^{-2S_{-1}/h}),
}
where $a_j(h)$ are real, $a_j(h)\sim a_{j,0}+a_{j,1}h+...$, $a_{j,0}>0$, $S_j=\phi (U_0)-\phi (U_j)$.
\end{theo}

\par


\begin{thebibliography}{30}

\bibitem{BaChSh} D.~Bao, S.-S.~Chern, Z.~Shen, \it An introduction to
Riemann-Finsler geometry, \rm Graduate texts in Mathematics, 200. Springer--Verlag,
New York, 2000.

\bibitem{Bi} J. M.~Bismut, \it The hypoelliptic Laplacian on the
cotangent bundle, \rm J. Amer. Math. Soc. 18(2005), 379-476.

\bibitem{BiLe} J. M.~Bismut, G.~Lebeau, \it The hypoelliptic Laplacian
  and Ray-Singer metrics, \rm preprint (2006).

\bibitem{DeVi}
L.~Desvillettes and C.~Villani, \it On the trend to global equilibrium in
spatially inhomogeneous entropy-dissipating systems: the linear
Fokker-Planck equation, \rm Comm. Pure Appl. Math.,
54(1)(2001), 1--42.

\bibitem{DiSj} J.~Sj\"ostrand, M.~Dimassi, \it
Spectral asymptotics in the semi-classical limit, \rm
London Math. Soc. Lecture Notes Series 269, Cambridge University Press 1999.

\bibitem{EcHa}
J.-P. Eckmann, M. Hairer, \it
Spectral properties of hypoelliptic operators, \rm
Comm. Math. Phys. 235(2)(2003), 233--253.

\bibitem{FrWe} M. I.~Freidlin, A. D.~Wentzell, \it
Random perturbations of dynamical systems, \rm Springer--Verlag, New York, 1984.

\bibitem{HeKlNi}
B.~Helffer, M.~Klein, F.~Nier, \it
Quantitative analysis of metastability in reversible diffusion
processes via a Witten complex approach, \rm Mat. Contemp. 26 (2004), 41--85.

\bibitem{HelNi} B.~Helffer, F.~Nier, \it
Hypoelliptic estimates and spectral theory for Fokker-Planck operators
and Witten Laplacians, \rm Lecture Notes in Mathematics, 1862.
Springer--Verlag, New York, 2005

\bibitem{HeSj} B.~Helffer, J.~Sj\"ostrand, {\it Multiple wells in the
semiclassical limit. {\rm I}},  Comm. Partial Differential Equations
9(4)(1984), 337--408, {\it Puits multiples en limite
semi-classique. {\rm II}.
Interaction mol\'eculaire. Sym\'etries. Perturbation}, \rm
 Ann. Inst. H. Poincar\'e Phys. Th\'eor.  42(2)(1985), 127--212.

\bibitem{HeSj3} B.~Helffer, J.~Sj\"ostrand, \it Multiple wells in the
semiclassical limit. {\rm III}. Interaction through nonresonant wells, \rm
Math. Nachr.  124(1985), 263--313.

\bibitem{HeSj4}
B.~Helffer, J.~Sj\"ostrand, \it Puits multiples en m\'ecanique
semi-classique. {\rm IV}. Etude du complexe de Witten, \rm
Comm. Partial Differential Equations  10(3)(1985), 245--340.

\bibitem{HeSj2} B.~Helffer, J.~Sj\"ostrand, \it R\'esonances en limite
semi-classique, \rm Bull. de la S.M.F., M\'emoire 24/25, Suppl. du Tome
114(3)(1986).

\bibitem{HeNi} F.~H\'erau, F.~Nier, \it
Isotropic hypoellipticity and trend to equilibrium for the
Fokker-Planck equation with a high-degree potential, \rm
Arch. Ration. Mech. Anal. 171(2)(2004), 151--218.

\bibitem{HeSjSt} F.~H\'erau, J.~Sj\"ostrand, C.~Stolk, \it Semiclassical
analysis for the Kramers-Fokker-Planck equation, \rm  Comm. Partial
Differential Equations 30(4-6)(2005), 689--760.

\bibitem{Hi} M.~Hitrik, \it
Boundary spectral behavior for semiclassical operators in dimension
one, \rm  Int. Math. Res. Not. 64(2004), 3417--3438.

\bibitem{Ko} V. N.~Kolokoltsov, \it Semiclassical analysis for
diffusions and stochastic processes, \rm Lecture Notes in Mathematics, 1724,
Springer-Verlag, Berlin, 2000.

\bibitem{Le} G.~Lebeau,  \it Le bismutien, \rm S\'eminaire \'equations aux
d\'eriv\'ees partielles, Ecole Polytechnique 2004--05, I.1--I.15

\bibitem{LiNi} Y.~Li, L.~Nirenberg, \it The distance \fu{} to the \bdy{},
Finsler's geometry and the singular set of viscosity solutions of some
Hamilton-Jacobi \e{}s, \rm Comm. Pure Appl. Math, 58(2005), 85--146.

\bibitem{MeSj} A.~Melin, J.~Sj\"ostrand, {\it \fop{}s with complex-valued
phase \fu{}s}, Springer Lect. Notes in Math., 459.

\bibitem{MeSj2} A. Melin, J. Sj\"ostrand, {\it Determinants of pseudodifferential operators and complex deformations
of phase space}, Meth. Appl. Analysis 9(2002), 177--238.

\bibitem{TaTaKu}
J.~Tailleur, S.~Tanase-Nicola, J.~Kurchan, \it
Kramers equation and supersymmetry, \rm  J. Stat. Phys.  122  (2006),
no. 4, 557--595. (preprint: arxiv.org/abs/cond-mat/0503545).

\bibitem{Sj} J.~Sj{\"o}strand, \it Parametrices for \pop{}s with multiple
characteristics, \rm Ark. f. Mat., 12(1)(1974), 85--130.

\bibitem{Sj2} J.~Sj{\"o}strand, \it Density of \res{}s for \st{} convex
obstacles, \rm Can. J. Math., 48(2)(1996), 397--447.

\bibitem{Vi} C.~Villani, {\it Hypocoercivity}, preprint, 2006 (arxiv.org/abs/math.AP/0609050).

\end{thebibliography}
\end{document}